%% file: _0DANDER.tex
\documentclass[12pt,a4paper,fleqn]{article}
\pdfoutput=1

\usepackage[english]{babel}
\RequirePackage{ucs}
\RequirePackage[utf8x]{inputenc}
\RequirePackage[T1]{fontenc}
\RequirePackage{amsthm,amsmath}
\RequirePackage{amssymb,amstext}
\RequirePackage{amsfonts,stmaryrd}
\RequirePackage{amsbsy} 
\RequirePackage{mathrsfs}
\RequirePackage{bm}
\RequirePackage{bbm}
\RequirePackage{aliascnt,ifthen}
\RequirePackage{graphicx}
\RequirePackage{color,calc}
\RequirePackage{times}
\RequirePackage{booktabs}
\usepackage{url}
\usepackage{verbatim}
\usepackage{subcaption}
\usepackage[pdfpagemode=UseOutlines ,plainpages=false,
hyperindex=true,colorlinks=true]{hyperref}
	\makeatletter
	\Hy@breaklinkstrue
	\makeatother

\definecolor{darkred}{rgb}{0.6,0,0.1}
\definecolor{darkgreen}{rgb}{0,0.5,0}
\definecolor{darkblue}{rgb}{0,0,0.5}

\hypersetup{colorlinks ,linkcolor=gray ,filecolor=darkgreen, urlcolor=darkblue ,citecolor=black}

\RequirePackage{paralist}
\RequirePackage{relsize}
\usepackage{enumitem}

\usepackage{natbib}
\usepackage{hypernat}
\renewcommand{\cite}{\citet}
\usepackage{DANDER-abrev-package}

\definecolor{dgreen}{rgb}{0,0.5,0}
\definecolor{dblue}{rgb}{0,0,0.5}
\definecolor{dred}{rgb}{0.6,0.0,0.1}
\definecolor{dgold}{rgb}{0.5,0.3,0.0}
\definecolor{dvio}{rgb}{0.6,0.3,0.5}
\definecolor{gray}{rgb}{0.5,0.5,0.5}
\definecolor{dbraun}{rgb}{.5,0.2,0}

\newcommand{\dr}{\color{dred}}
\newcommand{\db}{\color{dblue}}

\newcommand{\dgrau}{\color{gray}}

\newcommand{\colre}{dred}
\newcommand{\colas}{dblue}
\newcommand{\colrem}{dgold}
\newcommand{\colil}{dgreen}

\newtheoremstyle{styre}
  {1.1\topsep}
  {\topsep}
  {\normalfont\itshape}
  {}
  {\color{\colre}}
  {.}
  {.5em}
  {\thmname{\textbf{#1}\xspace}{\xspace\dgrau\thmnumber{#2}}\thmnote{\xspace\textit{\small(#3)}}}

\newtheoremstyle{styas}
  {1.1\topsep}
  {\topsep}
  {\normalfont\itshape}
  {}
  {\color{\colas}}
  {.}
  {.5em}
  {}

\newtheoremstyle{styrem}
  {1.1\topsep}
  {\topsep}
  {\normalfont\itshape}
  {}
  {\color{\colrem}}
  {.}
  {.5em}
  {\thmname{\textbf{#1}\xspace}{\xspace\dgrau\thmnumber{#2}}\thmnote{\xspace\textit{\small(#3)}}}

\newtheoremstyle{styil}
  {1.1\topsep}
  {\topsep}
  {\normalfont\rmfamily}
  {}
  {\color{\colil}}
  {.}
  {.5em}
  {\thmname{\textbf{#1}\xspace}{\xspace\dgrau\thmnumber{#2}}\thmnote{\xspace\textit{\small(#3)}}}

\newtheoremstyle{stypro}%
	{0.5\topsep}
	{1.1\topsep}
	{\upshape}
	{}
	{}
	{.}
	{.5em}
	{\thmnote{\textit{#3}}}

\theoremstyle{styre}
\newaliascnt{co}{pr}
\theoremstyle{styre}\newtheorem{co}[co]{Corollary}
\aliascntresetthe{co}
\newaliascnt{thm}{pr}
\theoremstyle{styre}\newtheorem{thm}[thm]{Theorem}
\aliascntresetthe{thm}
\newaliascnt{lem}{pr}
\theoremstyle{styre}\newtheorem{lem}[lem]{Lemma}
\aliascntresetthe{lem}
\newaliascnt{rem}{pr}
\theoremstyle{styrem}\newtheorem{rem}[rem]{Remark}
\aliascntresetthe{rem}
\newaliascnt{il}{pr}
\theoremstyle{styil}
\aliascntresetthe{il}
\theoremstyle{styas}\newtheorem{ass}{Assumption}

\theoremstyle{styas}

\theoremstyle{stypro}\newtheorem*{pro}{}

\newcommand{\remEnd}{{\scriptsize\textcolor{\colrem}{\qed}}}
\newcommand{\proEnd}{{\scriptsize\textcolor{\colre}{\qed}}}
\newcommand{\ilEnd}{{\scriptsize\textcolor{\colil}{\qed}}}

\usepackage{cleveref}
\crefname{pr}{\color{\colre}Proposition}{\color{\colre}Propositions}
\crefname{co}{\color{\colre}Corollary}{\color{\colre}Corollaries}
\crefname{thm}{\color{\colre}Theorem}{\color{\colre}Theorems}
\crefname{lem}{\color{\colre}Lemma}{\color{\colre}Lemmata}
\crefname{ass}{\color{\colas}Assumption}{\color{\colas}Assumptions}
\crefname{de}{\color{\colas}Definition}{\color{\colas}Definitions}
\crefname{rem}{\color{\colrem}Remark}{\color{\colrm}Remarks}
\crefname{il}{\color{\colil}Illustration}{\color{\colil}Illustrations}

\usepackage{environ}
\NewEnviron{te}{%
{\BODY}%
}

\numberwithin{equation}{section} 

\makeatletter
\newcommand{\mylabel}[2]{#2\def\@currentlabel{#2}\label{#1}}
\makeatother

\newcommand{\setListe}[5][3ex]{\setlength{\itemsep}{#2}\setlength{\topsep}{#3}\setlength{\leftmargin}{#4}\setlength{\rightmargin}{#5}\setlength{\labelwidth}{#1}}
\newcommand{\setListeStandard}{\setListe{0ex}{.5ex}{4ex}{0ex}}
\newcounter{ListeN}
\renewcommand{\theListeN}{(\roman{ListeN})}

\newenvironment{resListeN}[1][~]%
{\setcounter{ListeN}{0}\renewcommand{\theListeN}{\normalfont\rmfamily\color{\colre}(\roman{ListeN})}\begin{list}{\theListeN}%
{\usecounter{ListeN}\setListeStandard #1}}
{\end{list}}

\newenvironment{Liste}[1][~]%
{\begin{list}{}%
{\setListeStandard #1}}
{\end{list}}

\makeatletter%
\def\@fnsymbol#1{\ensuremath{\ifcase#1\or * \or 1 \or 2 \or 3 \or  *\or  \star \or 4\or  , \or 
g\or h\or i\else\@ctrerr\fi}}%
\makeatother%

\author{{\sc Sergio Brenner Miguel}\;\thanks{Institut f\"ur Angewandte
    Mathematik, M$\Lambda$THEM$\Lambda$TIKON, Im Neuenheimer Feld 205,
  D-69120 Heidelberg, Germany, e-mail:
  \url{{brennermiguel|johannes}@math.uni-heidelberg.de}} \and {\sc 
  Jan Johannes}$\;^*$}

\date{Ruprecht-Karls-Universität Heidelberg} 
\title{Data-driven aggregation in non-parametric density estimation on the real line} 
\begin{document} 
\maketitle 
\begin{abstract}
  We study non-parametric estimation of an unknown density with support in $\Rz$ (respectively  $\pRz$). 
  The proposed estimation procedure is based on the projection on finite dimensional subspaces spanned by the Hermite (respectively the Laguerre) functions. The focus of this paper is to introduce a data-driven aggregation approach in order to deal with the upcoming bias-variance trade-off. Our novel procedure integrates the usual model selection method as a limit case. 	We show the oracle- and the minimax-optimality of the
	data-driven aggregated density estimator and hence its adaptivity. 
	We present results of a simulation study which allow to compare the finite sample performance of the data-driven estimators using model selection compared to the new aggregation. 
\end{abstract} 
{\footnotesize
\begin{tabbing} 
aggregation, adaptation\\[.2ex] 
\end{tabbing}}
\input{_1Intro}
\input{_2Minimax}
\input{_3Aggreg}
\input{_4Numeric}
\appendix
\setcounter{subsection}{0}
\section*{Appendix}
\numberwithin{equation}{subsection}  
\renewcommand{\thesubsection}{\Alph{subsection}}
\renewcommand{\theco}{\Alph{subsection}.\arabic{co}}
\numberwithin{co}{subsection}
\renewcommand{\thelem}{\Alph{subsection}.\arabic{lem}}
\numberwithin{lem}{subsection}
\renewcommand{\therem}{\Alph{subsection}.\arabic{rem}}
\numberwithin{rem}{subsection}
\renewcommand{\thepr}{\Alph{subsection}.\arabic{pr}}
\numberwithin{pr}{subsection}
\input{_app1.tex}
\input{_app2.tex}
\input{_app3.tex} 
\bibliography{DANDER} 
\end{document}

%% file: _1Intro.tex
%
%
%
\section{Introduction}\label{i}
\begin{te}
In this paper we consider the data-driven estimation of an unknown
density $f$ with non-compact support in the real line given an
independent and identically distributed (i.i.d.) sample $X_1, \dots,
X_n$ from $\So$. In the literature, non-parametric density estimation is a well-discussed problem and many estimation strategies based on splines, kernels or wavelets, to name but a few,   are considered. For an overview of various methods we refer to  \cite{Comte2017}, \cite{Efromovich1999}, \cite{Silverman2018}   and  \cite{Tsybakov2008}. 
Here, we will focus on the projection of the density $\So$ on an orthonormal basis and therefore assume its square integrability. This has been studied  for densities  with compact support  (e.g.  \cite{Massart2007} and \cite{Efromovich1999}), with support in $\Rz$  using wavelets and Hermite functions (e.g. \cite{JuditskyLambert-Lacroixothers2004} and \cite{BelomestnyComteGenon-Catalot2019}, respectively) or with support in $\pRz$ using Laguerre functions  (e.g. \cite{ComteGenon-Catalot2018}).

Here, we cover  the projection on the Hermite functions or the Laguerre functions, i.e. the estimation over a set $A\subseteq \Rz$ where in the Hermite case \HerCase $A= \Rz$ and in the Laguerre case \LagCase $A=\pRz$.
In the sequel $\LpA$ denotes the set of all square integrable functions over $A$ endowed with its usual inner product $\skalarA$ 
	and norm $\normA$. Furthermore, let $\Kset[j\in \Nz_0:= \Nz \cup \{0\}]{\varphi_j}$ be in case \HerCase  and \LagCase the Hermite and Laguerre basis, respectively.
Therewith, for each $\So\in\LpA[2]$ we define the family
$\Ksuite[\Di\in\Nz]{\SoPr}$ of projections of $\So$ onto the subspaces
$\Ksuite[\Di\in\Nz]{S_{\Di}}$ where $S_k$ is the linear subspace
spanned by the first $k$ basis functions. 
	By replacing 
	the unknown coefficients by their empirical counterparts, we consider
	for $\SoPr$, $\Di\in\Nz$,  the unbiased orthogonal series estimator (OSE) 
	\begin{align}\label{equation:proj_est_const}
\hSoPr := \sum_{j=0}^{\Di-1} \hfSo{j}\bas_j\quad\text{ with }\quad\hfSo{j}:= n^{-1} \sum_{i=1}^n \bas_j(X_i).
\end{align}

To measure the performance of the estimator we discuss its mean
integrated squared error (MISE) as risk,  we state oracle rates and we derive upper bounds for its maximal risk over Sobolev classes. Further we show that the projection estimator  with optimal choice of the dimension parameter is minimax-optimal  over Sobolev classes.  The proof of the lower bound borrows ideas from  \cite{BelomestnyComteGenon-Catalotothers2017} and \cite{ComteDuvalSacko2019}.
In practice, however, the optimal choice of the dimension parameter is not feasible since it depends on characteristics of the unknown density $f$. Therefore, \cite{ComteGenon-Catalot2018} consider a model selection approach, inspired by the work of \cite{BarronBirgeMassart1999} and extensively described in  \cite{Massart2007}, to select fully data-driven the dimension parameter in such a way that the bias and variance compromise is automatically reached by the resulting estimator. 
More precisely, the authors choose the random dimension $\widehat k$ as  a minimum of $-\VnormA{\hSoPr}^2+\mathrm{pen}_k$ over the set of admissible parameters $\nset{1,\maxDi}:= [1, \maxDi] \cap \Nz$ for a given upper bound $\maxDi$ and sequence of penalty terms $(\mathrm{pen}_k)_{k\in \Nz}$.

In this work, we study a different data-driven procedure. Introducing $\maxDi\in\Nz$ and \textit{aggregation weights} $\We[]:=(\We[m])_{m\in \nset{1, \maxDi}} \in [0,1]^{\maxDi}$ with $\sum_{i=1}^{\maxDi} \We[i]=1$  we define the \textit{aggregated estimator} $\hSoPr[{\We[]}]:= \sum_{k=1}^{\maxDi} \We[k]\hSoPr$. Note that the \textit{aggregation} is called fully data-driven  if the \textit{aggregation weights} depend on the data only. The model selection procedure can be integrated into this aggregation framework via the \textit{\mwn} 
\begin{equation}\label{ag:de:msWe}\msWe:=\dirac[\hDi](\{\Di\}),\quad
  \Di\in\nset{1,\maxDi}
\end{equation}
where $\delta_x$ denotes the usual Dirac measure in $x\in \Rz$. 
However we suggest as a new data-driven choice the  weights defined by
\begin{equation}\label{ag:de:rWe}
\rWe:=\frac{\exp(-\rWn\{-\VnormA{\hSoPr}^2+\hpen\})}
{\sum_{l=1}^{\maxDi}\exp(-\rWn\{-\VnormA{\hSoPr[l]}^2+\hpen[l]\})},\quad \Di\in\nset{1,\maxDi}
\end{equation}
where the choice of the penalties $(\hpen)_{k\in \Nz_0}$, the numerical constant $\kappa\geq1$ and the upper bound  $\maxDi$ will be further discussed in Section \ref{ag}. 
We refer to them as \textit{\bwn} since their particular form   takes its inspiration from  a-posteriori weights in a Bayesian sequence space model (c.f. \cite{JohannesSimoniSchenk2015}).  In this paper we derive upper bounds for the (maximal) risk of the aggregated estimator using either \bwn or \mwn where throughout the paper $\hDi$ is chosen as  a minimum of $-\VnormA{\hSoPr}^2+\hpen$ over $\nset{1,\maxDi}$. Note that in this situation the \bwn  converge to  the corresponding \mwn as $\rWc$ tends to infinity, i.e., $\lim_{\rWc\to\infty} \rWe=\msWe$ for each $\Di\in\nset{1,\maxDi}$. 
\end{te}
\begin{te}The paper is organised as follows: in Section \ref{mm} we introduce
	our basic assumptions, recall the oracle inequalities and develop the minimax
	theory. 
	We show, in Section \ref{ag}, the oracle- and the minimax-optimality of the
	data-driven aggregated density estimator and hence its adaptivity. In Section \ref{mm} and \ref{ag}  we only present key arguments of the proofs while more technical details are postponed to the  \cref{a:mm,a:ag}, respectively.
	Finally,  results of a simulation study are
	reported in Section \ref{si} which allow to compare the finite sample performance of the aggregated estimator with \mwn and \bwn  of a density given independent observations. Further we introduce the Laguerre and Hermite functions and recall some of their properties in the \cref{a:prel}.
\end{te}


%% file: _2Minimax.tex
%
%
%
\section{Minimax theory}\label{mm}
Given an orthonormal basis $\Ksuite[j\in\Nz_0]{\bas_j},$ we
consider for any function $f\in\LpA[2]$ its expansion
$f =\sum_{j\in\Nz_0}\fSo{j}\bas_j$ with $\fSo{j} := \VskalarA{\So,
  \bas_j}$ and for each $\Di\in\Nz$ the subspace  $S_{\Di}$ spanned by the first
$\Di$ basis functions $\Kset[j\in\nsetro{0,\Di}]{\bas_j}$, where here and subsequently for real numbers  $a\leq b$ we write shortly $\nsetro{a,b}:=[a,b)\cap
\Zz$,  $\nsetlo{a,b}:=(a,b]\cap \Zz$,  and so forth. Consequently, the
projection of $f\in\LpA[2]$ onto $S_{\Di}$ is given by $\SoPr =
\sum_{j=0}^{\Di-1}\fSo{j}\bas_j$. For each $\Di\in\Nz_0$ and density
$\So\in\LpA[2]$ we define 
  $\bias^2(\So)\in[0,1]$ as follows
  $\VnormA{\So}^2\bias^2(\So)=\VnormA{\SoPr-\So}^2=\sum_{j\geq\Di}|\fSo{j}|^2$,
  where   we agree on $\SoPr[0]:=0$ and hence
  $\bias[0]^2(\So)=1$.  Let $\FuEx[]{\So}$ and
$\FuEx[n]{\So}$ denote, respectively, the
expectation with respect to the marginal and joint distribution $\FuVg[n]{\So}$ of
the i.i.d. $n$-sample
$\Ksuite[i\in\nset{1,n}]{X_i}$.

\paragraph{Oracle optimality.}   Elementary calculations show for each $\Di\in\Nz$ the  identity 
\begin{multline}\label{oo:rdec}
\FuEx[n]{\So}\big(\VnormA{\hSoPr-\So}^2\big) +
  \tfrac{1}{n}\VnormA{\So}^2=
  \tfrac{n+1}{n}\VnormA{\So-\SoPr}^2 + \tfrac{1}{n} V_{\Di}  \text{ with } V_{\Di}:= \sum_{j=0}^{\Di-1} \FuEx{\So}\big(\bas_j^2(X_1)\big).
\end{multline}
In \cref{a:prel} we briefly recall elementary properties of Laguerre
and Hermite functions. As for example, that they are bounded in the usual
uniform norm, precisely, $\sup_{j\in\Nz_0}\VnormInf{\bas_j}\leq C$  and,
hence $\sup_{\Di\in\Nz}\tfrac{1}{\Di}V_{\Di} \leq C$ in case \LagCase and 
\HerCase  with $C= \sqrt{2}$ and $C= 1$, respectively. Moreover,  \cite{ComteGenon-Catalot2018} and
\cite{BelomestnyComteGenon-Catalot2019} have  shown sharper upper and lower
bounds for the term $V_k$.
Precisely, setting
\begin{multline}\label{de:bas:vFu}
\dSob:=1,\quad \aSob:=-1/2,\quad  \vFu:= \FuEx{\So}(X^{\aSob})+1 \text{ in case \LagCase and }\\
\dSob:=10/12,\quad \aSob:=2/3,\quad \vFu:=\FuEx{\So}(|X|^{\aSob})+1\text{ in case \HerCase} 
\end{multline}
there exists a numerical
constant $\Vcst\geq1$ such that for each $\Di\in\Nz$ hold
\begin{equation}\label{de:bas:cst}
  \VnormInf{\sum_{j=0}^{\Di-1} \bas_j^2}\leq\Vcst\,\Di^{\dSob}\quad\mbox{and}\quad  \sum_{j=0}^{\Di-1}\FuEx{\So}\big(
  \bas_j^2(X)\big)\leq \Vcst\, \vFu\, \Di^{1/2}. 
\end{equation}
For a sequence $\Nsuite[n]{a_n}$ of real numbers with minimal value in a set
$B\subset{\Nz}$ we define
$\argmin\set{a_n,n\in B}:=\min\{m\in B:a_m\leq a_n,\;\forall n\in B
\}$. For $n,\Di\in\Nz$ we set 
\begin{multline}\label{de:oRa}
\dRaSo:=[\bias^2(\So)\vee n^{-1}\Di^{1/2}],\quad 
 \oDiSo:=\argmin\Nset[\Di\in\Nz]{\dRaSo}\quad\text{and}\\
 \oRaSo :=\min\Nset[\Di\in\Nz]{\dRaSo}.\hfill\end{multline}
 Here for two real numbers $a,b\in \Rz$ we define $a\vee b:= \max\{a,b\}$ and $a\wedge b:=\min\{a,b\}$.
\begin{rem}\label{ag:rem:pen:oo} Note that by construction
  $\oRaSo=\dRaSo[\oDiSo]$ and  $\oDiSo\in\nset{1,n^2}$, since  $\bias[n^2]^2(\So)\leq 1<(n^2+1)^{1/2}n^{-1}$, and hence
  $\dRaSo[n^2]<\dRaSo$ for all $\Di\in \nsetro{n^2+1,\infty}$.
 It is worth stressing
  out that in compact density estimation the oracle dimension
  typically satisfies  $\oDiSo\in\nset{1,n}$. 
  Obviously, it follows thus
  $\oRaSo=\min\set{ \dRaSo,\Di\in\nset{1,n^2}}$ for all $n\in\Nz$. 
  Moreover,  we shall emphasise that $\oRa\geq n^{-1}$ for all
  $n\in\Nz$, and $\oRaSo=o(1)$ as $n\to\infty$. We eventually use those
  elementary findings in the sequel without further
  reference.\remEnd
\end{rem}
Combining \eqref{oo:rdec} and \eqref{de:bas:cst} we immediately obtain  
\begin{equation}\label{oo:oe1}
  \inf\NsetB[\Di\in\Nz]{\nEx\VnormA{\hSoPr-\So}^2}\leq \nEx\VnormA{\hSoPr[\oDiSo]-\So}^2
  \leq [\tfrac{n+1}{n}\VnormA{\So}^2+\Vcst \, \vFu]\oRaSo.
\end{equation}
The upper bound \eqref{de:bas:cst} for the variance term cannot be
improved, since under the additional assumption $c_{\So}:=\inf_{x\in[a, b]}
\So(x)>0$ for some $a,b\in\Rz$, $a<b$, there exists a
constant $\cst{\So,a,b}>0$ depending on $a,b$ and $\So$ such that for
each $\Di$ holds $V_{\Di} \geq
\cst{\So,a,b} \sqrt{\Di}$ (e.g. \cite{ComteGenon-Catalot2018}). In this situation, from \eqref{oo:rdec} it follows also 
\begin{equation}\label{oo:oe2}
  \inf\NsetB[\Di\in\Nz]{\nEx\VnormA{\hSoPr-\So}^2}\geq [\tfrac{n+1}{n}\VnormA{\So}^2 \wedge \cst{\So,a,b}-\tfrac{\VnormA{\So}^2}{n\oRaSo}]\oRaSo.
\end{equation}
Consequently, the  rate $\Nsuite[n]{\oRaSo}$, the dimension parameters $\Nsuite[n]{\oDiSo}$  and  the OSE's  $\Nsuite[n]{\hSoPr[\oDiSo]}$, respectively, is an oracle
rate, an oracle dimension and oracle optimal (up to a constant) as soon as the leading factor on the right hand side is positive.
\begin{te}
  Throughout the paper we shall distinguish for the 
  density $\So$ and hence it's associated sequence
  $\Nsuite[\Di]{\bias(\So)}$ the following two cases
  \begin{Liste}[]
  \item[\mylabel{oo:p}{\dgrau\upshape\bfseries{(p)}}] there is $K\in\Nz$
    with $\bias[K-1](\So)>0$ (with $\bias[0](\So)=1$) and
    $\bias[K](\So)=0$,
  \item[\mylabel{oo:np}{\dgrau\upshape\bfseries{(np)}}] for all $K\in\Nz$
    holds $\bias[K](\So)>0$.
    \end{Liste}
  \end{te}
\begin{rem}\label{oo:rem:ora}
Note that the expansion of $\So$
is in case \ref{oo:p} \textit{finite}, i.e., $\So=\sum_{j=0}^{K-1}\fSo{j}\bas_j$
for some $K\in\Nz$ while in the opposite case
\ref{oo:np}, it isn't. Interestingly,  in case \ref{oo:p} the oracle
rate is parametric, that is, 
$\oRaSo$ is of order $n^{-1}$. More precisely, 
if there is $K\in\Nz$  with $\bias[K-1](\So)>0$ and
$\bias[K](\So)=0$, then setting
$n_{\So}:=\tfrac{K^{1/2}}{\bias[K-1]^2(\So)}$, for all
$n\geq n_{\So}$ holds
$\bias[K-1]^2(\So)>K^{1/2}n^{-1}$, and hence  $\oDiSo=K$ and
$\oRaSo= K^{1/2}n^{-1}$.
On the other hand side, in case \ref{oo:np} the oracle rate is
non-parametric, more precisely, it holds
$\lim_{n\to\infty}n\oRaSo=\infty$. Indeed, since
$\bias[\oDiSo]^2(\So)\leq\dRaSo[\oDiSo]=\oRaSo=o(1)$ as $n\to\infty$
follows $\oDiSo\to\infty$, which implies the claim because
$n\oRaSo\geq(\oDiSo)^{1/2}$.\remEnd
\end{rem}
\begin{te}
  Let us first briefly illustrate the last definitions by stating the
  order of $\oDiSo$ and $\oRaSo$ in case \ref{oo:np} for an often
  considered behaviour of the sequence $\Nsuite[\Di]{\bias^2(\So)}$.
  Here and subsequently, we use for two strictly positive sequences
  $\Nsuite[n]{a_{n}}$ and $\Nsuite[n]{b_{n}}$ the notation
  $a_{n}\sim b_{n}$ if the sequence $\Nsuite[n]{a_{n}/b_{n}}$ is
  bounded away both from zero and infinity.  Let us use
  $\bias^2(\So)\sim \Di^{-s}$, $s>0$, as a particular specification.
  In this situation $\oDiSo\sim n^{2/(2s+1)}$ is the oracle dimension
  and $\oRaSo\sim n^{-2s/(2s+1)}$ is the oracle rate.
  \paragraph{Minimax optimality.} For each $\Di\in\Nz$ let us measure
  the accuracy of $\hSoPr$ by its maximal risk over classical Sobolev
  spaces or ellipsoids, respectively, defined for
  $\wSob, \rSob\in\pRz$ by
  \begin{align*}
    \rwcSob{\wSob}&:=\{f \in \LpA[2] : |f|_{\wSob}^2 :=
                    \sum_{\Di\in \Nz_0} \Di^{\wSob} |\fSo{\Di}|^2<
                    \infty\}
                    \text{ and } \rwcSob{\wSob,\rSob}:= \{f \in \rwcSob{\wSob}: |f|_{\wSob}^2\leq \rSob\}.
  \end{align*}
  For a more detailed discussion on Sobolev-Laguerre spaces
  $\rwcSob[\pRz]{s}$ and Sobolev-Hermite space $\rwcSob[\Rz]{s}$ we
  refer to \cite{BongioanniTorrea2009} and
  \cite{BongioanniTorrea2006}, respectively. For $\aSob$ as in
  \eqref{de:bas:vFu} we denote further a corresponding subset of
  densities with finite $\aSob\mSob$-th moment, $\mSob\in\pRz$, by
  \begin{equation}\label{equation:dens_sobol}
    \rwcSobD{\wSob,\rSob,\mSob}:= \{f\in \rwcSob{\wSob,\rSo}: f \text{ is
      a density and } \FuEx{\So}\big(|X|^{\aSob\mSob}\big) \leq L\}
  \end{equation}
  and let
  $\rwcSobD{\wSob,\mSob} := \bigcup_{\rSob>0}
  \rwcSobD{\wSob,\rSob,\mSob}$. We note, that for each
  $\So\in \rwcSobD{\wSob,\rSob,\mSob}$ with $\mSob\geq1$ and $ \vFu$
  as in \eqref{de:bas:vFu} we have $\vFu\leq \rSob +1$. Moreover, for
  each $\Di\in\Nz$ holds
  $\VnormA{\So}^2\bias^2(\So)=\VnormA{\So-\SoPr}^2\leq
  \rSob\Di^{-\wSob}$. Exploiting the upper bound \eqref{oo:oe1} there
  is a constant $\cst{\wSob,\rSob}$ depending on the class
  $\rwcSobD{\wSob,\rSob,\mSob}$ only, such that for each $n\in\Nz$
  with $\oDi\sim n^{2/(2s+1)}$ it holds
  \begin{equation}\label{theorem:directupper}
    \mmRi[\oDi]{\wSob}{\mSob}
    \leq \cst{\wSob,\rSob} n^{-2s/(2s+1)}.
  \end{equation}
  Now we provide a lower bound to show that 
  the upper bound in  \eqref{theorem:directupper} is minimax-optimal
  over the ellipsoid $\rwcSobD{\wSob,\rSob,\mSob}$. Further the
  following theorem is formulated for the more general ellipsoids
  $\rwcSobD{\wSob,\rSob,\mSob}$ because the data-driven
  aggregation requires stronger moment assumptions.
\end{te}
\begin{thm}\label{theorem:lower_bound}
  Let $n,\mSob\in\Nz$ with $n\geq n_s=8^{2s+1}$, then there exist
  constants $\cst{\mSob},L_{\wSob,\mSob}>0$ such that for all
  $L\geq L_{\wSob,\mSob}$ and for any estimator $\hSo$ of $\So$ based
  on an i.i.d. sample $\Nsample{X_j}$ 
  \begin{align*}
   \mmRi[]{\wSob}{\mSob}
    \geq \cst{\mSob} n^{-2s/(2s+1)}.
  \end{align*}
\end{thm}
\begin{pro}[Proof of \cref{theorem:lower_bound}]
  We outline here the main steps of the proof, while more technical
  details are deferred to the appendix.  We will construct a family of
  functions in $\rwcSobD{\wSob,\rSob,\mSob}$ by a perturbation of a
  density $\SoPr[o]$ with small bumps, such that their $\LpA$-distance
  and their Kullback-Leibler divergence can be bounded from below and
  above, respectively. The claim follows then by applying Theorem 2.5
  in \cite{Tsybakov2008}. In both cases \LagCase and \HerCase we use
  the same construction, which we present first.\\ 
   Given a function $\psi\in\LpA$ for each $K\in\Nz$ (to be selected below) and
  $k\in\nsetro{0,K}$ we define the bump-functions
  $\psi_{k, K}(x):= \psi(xK-K-k),$ $x\in\Rz$. For a density
  $\SoPr[o]\in\LpA$ (specified in \cref{lemma:Lag_SobDen,lemma:Her_SobDen} in the \cref{a:mm}), a bump-amplitude $\delta>0$ and a vector
  $\bm{\theta}=(\theta_1,\dots,\theta_K)\in \{0,1\}^K$ we define
  \begin{equation}\label{equation:lobodens}
    \SoPr[\bm{\theta}](x)=\SoPr[o](x)+ \delta K^{-s} \sum_{k=0}^{K-1}
    \theta_{k+1} \psi_{k, K}(x).
  \end{equation}
  The choice of $\psi$ is discussed in \cref{a:mm} too, however, it
  ensures that $\int_A\psi(x)dx=0$, and hence,
  $\SoPr[\bm{\theta}]$ integrates to one, and that the support of
  $\psi$ is contained in $[0,1]$. Moreover, $\SoPr[o]$ satisfies $c_{\SoPr[o]}:=\inf_{x \in [1,2]} \SoPr[o] >0$ which in turn for any $\delta \in (0,\delta_{\SoPr[o],\psi}]$ with
  $\delta_{\SoPr[o], \psi}:=c_{\SoPr[o]}/\VnormInf{\psi}$ implies
  $\SoPr[\bm\theta](x) \geq 0$ for all $x\in A$. Indeed, on $[1,2]^c $
  holds $\SoPr[\bm\theta] = \SoPr[o]$ and for the non trivial case $x\in [1, 2]$ there is $k_o\in\nsetro{0,K}$ such that $x \in [1+k_o/K,1+ (k_o+1)/K]$ and hence
  \begin{equation*}
   \SoPr[\bm{\theta}](x)= \SoPr[o](x) + \theta_{k_o+1}\delta K^{-\wSob} \psi(xK-K-k_o) \geq c_{\SoPr[o]}  - \delta \| \psi\|_{\infty}K^{-\wSob}\geq0.
 \end{equation*}
  Moreover, due to \cref{lemma:Lag_SobDen} and
  \cref{lemma:Her_SobDen} for the cases \LagCase and \HerCase,
  respectively,  $\SoPr[o]$ and  the family $\{\SoPr[\bm\theta] : \bm\theta \in \{0,1\}^K\}$ belong to
  $\rwcSobD{\wSob,L_{\wSob,\mSob},\mSob}$ for some $L_{\wSob,\mSob}>0$.
  Exploiting Varshamov-Gilbert's
  Lemma (see \cite{Tsybakov2008}) in \cref{lemma:tsyb_vorb} we show
  further that there is $M\in\Nz$ with $M\geq 2^{K/8}$ and a subset
  $\{\bm \theta^{(0)}, \dots, \bm \theta^{(M)}\}$ of $\{0,1\}^K$ with
  $\bm \theta^{(0)}=(0, \dots, 0)$ such that for all
  $j, l \in \nset{0, M}$, $j \neq l$ the $\LpA$-distance and the
  Kullback-Leibler divergence are bounded: 
  \begin{equation}\label{vg:e1}
    \VnormA{\SoPr[\bm{\theta^{(j)}}]-\SoPr[\bm{\theta^{(l)}}]}^2 \geq
    \cst[(1)]{\psi, \delta} K^{-2s}\quad\text{ and }\quad
    \text{KL}(\SoPr[\bm{\theta^{(j)}}],\SoPr[\bm{\theta^{(0)}}]) \leq \cst[(2)]{\psi, \delta} \log(M) K^{-2s-1}
  \end{equation}
  where $\cst[(1)]{\psi, \delta}>0$ and
  $\cst[(2)]{\psi, \delta}<\infty$ depend on $\psi$
  and $\delta$ only.  
  Selecting $K=\ceil{n^{1/(2s+1)}}$  follows
  \begin{align*}
    \frac{1}{M}\sum_{j=1}^M
    \text{KL}((\SoPr[\bm{\theta^{(j)}}])^{\otimes
    n},(\SoPr[\bm{\theta^{(0)}}])^{\otimes n})
    &= \frac{n}{M} \sum_{j=1}^M \text{KL}(
      \SoPr[\bm{\theta^{(j)}}],\SoPr[\bm{\theta^{(0)}}] )
      \leq \cst[(2)]{\psi, \delta}  \log(M)
  \end{align*}
  where $\cst[(2)]{\psi, \delta}< 1/8$ for all
  $\delta^2<\frac{\log(2)}{8c_{f_o}\VnormA{\psi}^2}$ and $M\geq 2$ for
  $n\geq n_0=8^{2s+1}$. Thereby, we can use Theorem 2.5 of
  \cite{Tsybakov2008}, which in turn for any estimator $\hSo$ of $\So$
  implies
  \begin{multline*}
    \sup_{\So\in\rwcSobD{\wSob,\rSob,\mSob}}
    \nVg\big(\VnormA{\hSo-\So}^2\geq
    \tfrac{\cst[(1)]{\psi, \delta}}{2}n^{-2s/(2s+1)} \big)\geq
    \tfrac{\sqrt{M}}{1+\sqrt{M}}\big(1-2\cst[(2)]{\psi, \delta}
    -\sqrt{\tfrac{2\cst[(2)]{\psi, \delta}}{\log(M)}} \big) \geq 0.07.
  \end{multline*}
  Note that the constant $\cst[(1)]{\psi, \delta}$ does only depend on
  $\psi$ and $\delta$, hence implicitly also on $\mSob$, but
  it is independent of the parameters $\wSob,\rSob$ and $n$. The claim
  of \cref{theorem:lower_bound} follows by using Markov's inequality,
  which completes the proof.\proEnd
\end{pro}

%% file: _3Aggreg.tex
%
%
%
\section{Data-driven aggregation}\label{ag}
\begin{te}
 Given a family $\Kset[\Di\in\nset{1,\maxDi}]{\hSoPr}$ of orthogonal
 series estimators as in \eqref{equation:proj_est_const} the optimal choice of the dimension parameter $\Di$ in an oracle or  minimax sense,
  does depend on characteristics of the unknown density.  Introducing \textit{aggregation weights}
  $\We[]=(\We[\Di])_{\Di\in \nset{1, M_n}} \in [0,1]^{\maxDi}$ with
  $\sum_{\Di=1}^{\maxDi} \We[\Di]=1$ we consider here and subsequently
  the aggregation $\hSoPr[{\We[]}]= \sum_{k=1}^{\maxDi}
  \We[k]\hSoPr$. Note that the aggregation weights  define a discrete probability measure
  $\FuVg{\We[]}$ on the set $\nset{1,\maxDi}$ by
  $\FuVg{\We[]}(\{\Di\})=\We$.  Clearly,  the random coefficients
  $\Ksuite[j\in\Nz_0]{\fou[j]{\hSoPr[{\We[]}]}}$
 of $\hSoPr[{\We[]}]$ satisfy $\fou[j]{\hSoPr[{\We[]}]}=0$ for
 $j\geq \maxDi$ and  for any  $j\in\nsetro{0,\maxDi}$ holds
 $\fou[j]{\hSoPr[{\We[]}]}=(\sum_{\Di=j+1}^{\maxDi}\We)\times\hfSo{j}=\FuVg{\We[]}(\nsetlo{j,\maxDi})\times\hfSo{j}$. Our aim is to prove an upper bound for its risk $\nEx\VnormA{\hSoPr[{\We[]}]-\So}^2$ and its  maximal risk
  $\nRi{\hSoPr[{\We[]}]}{\rwcSobD{\wSob,\rSob,\mSob}}{}$. For 
  arbitrary aggregation weights and penalty sequence the next lemma establishes an
  upper bound for the loss of the aggregated estimator. Selecting suitably
  the weights and penalties  this bound provides in the sequel our key
  argument. 
\end{te}
\begin{lem}\label{co:agg}
Consider an  aggregation $\hSoPr[{\We[]}]=\sum_{\Di=1}^{\maxDi}\We\hSoPr$
and   sequences
$(\pen)_{\Di\in\nset{1,\maxDi}}$ and $(\hpen)_{\Di\in\nset{1,\maxDi}}$ of non-negative penalty terms. For any $\mDi\in\nset{1,\maxDi}$ and $\pDi\in\nset{1,\maxDi}$  holds
  \begin{multline}\label{co:agg:e1}
    \VnormA{\hSoPr[{\We[]}]-\So}^2\leq \tfrac{2}{14}\pen[\pDi] +2\VnormA{\So}^2\bias[\mDi]^2(\So)\\\hfill
    +2\VnormA{\So}^2\FuVg{\We[]}(\nsetro{1,\mDi})+\tfrac{2}{14}\sum_{\Di=1+\pDi}^{\maxDi}\pen\We\Ind{\{\VnormA{\hSoPr-\SoPr}^2<\hpen/7\}}\\
+2\sum_{\Di=\pDi}^{\maxDi}\vectp{\VnormA{\hSoPr-\SoPr}^2-\pen/14}  
+\tfrac{2}{14}\sum_{\Di=1+\pDi}^{\maxDi}\We\pen\Ind{\{\VnormA{\hSoPr-\SoPr}^2\geq\hpen/7\}}
\end{multline}
\end{lem}
\begin{rem}\label{co:agg:rem}Keeping \cref{co:agg} in mind let us outline briefly
  the principal arguments of our aggregation strategy.
  Selecting the values of $\pDi$ and $\mDi$  
 close to the oracle dimension $\oDiSo$ the first two terms 
  in the upper bound of \eqref{co:agg:e1} are of the order of the oracle rate.  The
  weights are on the other hand selected such that
  the third and fourth term on the right hand side in \eqref{co:agg:e1} are
  negligible with respect to the oracle rate, while the choice of the
  penalties allows  as usual  to bound the deviation of the last two terms 
by concentration inequalities.\remEnd
\end{rem}
  \paragraph{Risk bounds.}
  We derive bounds for the risk of the aggregated estimator $\hSoPr[{\We[]}]$
  using either \bwn $\We[]:=\rWe[]$
  as in \eqref{ag:de:rWe} or \mwn $\We[]:=\msWe[]$
  as in \eqref{ag:de:msWe}. Until now we have
  not yet specified  the sequences of
  penalty terms. Keeping in mind that the oracle dimension $\oDiSo$
  belongs to $\nset{1,n^2}$ we set  $\maxDi:=\floor{n^2(600\log n)^{-4}}$ and $\maxDi:=n^2$ in case \LagCase and \HerCase, respectively.  We estimate $\vFu$ defined  in
  \eqref{de:bas:vFu} by its empirical counterpart $\hvFu:=1+\tfrac{1}{n}\sum_{i=1}^n
  |X_i|^{\aSob}$ with  $\aSob=-1/2$ in case
\LagCase and $\aSob=2/3$ in case \HerCase. 
For  each $\Di\in\Nz$ and a numerical
constant $\cpen>0$ we set
\begin{equation}\label{ag:de:pen}
\pen:=\vpen\quad\text{ and }\quad \hpen:=\hvpen.
\end{equation}
  Our theory necessitates a lower bound for the numerical constant
  $\Delta$ which is for practical application in general too large. In
  the simulations we use  preliminary experiments
  to determine a good choice for $\Delta$ (c.f. \cite{BaudryMaugisMichel2012}).
\begin{thm}\label{ag:ub:pnp} 
  Consider an aggregation $\hSoPr[{\We[]}]$   using either \bwn $\We[]:=\rWe[]$
  as in \eqref{ag:de:rWe} or \mwn $\We[]:=\msWe[]$
  as in \eqref{ag:de:msWe} with $\Delta \geq  84 \Vcst[]$ and
  $\Vcst[]$ as in \eqref{de:bas:cst}. Suppose the density
  $\So\in\LpA$ satisfies $\VnormInf{\So}<\infty$ and $\FuEx{\So}(|X_1|^{2\aSob})<\infty$. \begin{Liste}[]
\item[\mylabel{ag:ub:pnp:p}{\dgrau\upshape\bfseries{(p)}}]Assume there is $K\in\Nz$
  with   $1\geq \bias[{K-1 }](\So)>0$ and $\bias[K](\So)=0$. 
Then there is a finite constant $\cst{\So}$
given in \eqref{ag:ub:pnp:p6} depending only on $\So$ such that for all $n\geq3$ holds
\begin{equation}\label{ag:ub:pnp:e1}
  \oRi  \leq \cst{\So}n^{-1}.
\end{equation}
\item[\mylabel{ag:ub:pnp:np}{\dgrau\upshape\bfseries{(np)}}] If
  $\bias(\So)>0$ for all  $\Di\in\Nz$, then there is a  numerical constant $\cst{}$ 
such that for all
$n\geq3$ 
 \begin{multline}\label{ag:ub:pnp:e2}
  \oRi 
  \leq
  \cst{}\big((\VnormA{\So}^2+\vc)\mRaSo+([\VnormInf{\So}^3\vee1]\,\vc + [\VnormA{\So}^2\vee1]\,  \FuEx{\So}(|X_1|^{2\aSob}))\,n^{-1}\big)\\
    \text{ with }\mRaSo:=\min_{\Di\in\nset{1,\maxDi}}\setB{\dRaSo\vee\exp\big(-\tfrac{\Vcst[]\vc}{1\vee400\VnormInf{\So}}\Di^{1/2}\big)}.
\end{multline}
\end{Liste}  
\end{thm}
\begin{te}
Before we proof the main result. Let us state an immediate consequence.  
\end{te}
\begin{co}\label{ak:ag:ub2:pnp}
  Let the assumptions of \cref{ag:ub:pnp} be satisfied. If in case \ref{ag:ub:pnp:np} in addition
    \begin{inparaenum}\item[\mylabel{ag:ub2:pnp:npc}{{\dgrau\upshape\bfseries(A1)}}]
      there is  $n_{\So}\in\Nz$ such that $\oDiSo$ and $\oRaSo$  as in \eqref{de:oRa} satisfy
      $\maxDi^{1/2}\geq(\oDiSo)^{1/2}\geq \tfrac{1\vee400\VnormInf{\So}}{\Vcst[]\vc}|\log\oRaSo|$
      for all $n\geq n_{\So}$,
    \end{inparaenum}
    then there is a constant $\cst{\So}$ depending
    only on $\So$ such that $\oRi\leq\cst{\So}\oRaSo$ for all $n\in\Nz$ holds true.
  \end{co}
\begin{pro}[Proof of \cref{ak:ag:ub2:pnp}]If the additional assumption
  \ref{ag:ub2:pnp:npc} is satisfied, then the oracle dimension $\oDiSo$ 
     satisfies trivially $\oDiSo\in\nset{1,\maxDi}$ and
  $\exp\big(\tfrac{-\Vcst[]\vc}{1\vee400\VnormInf{\So}}\oDiSo^{1/2}\big)\leq
  \oRaSo$ while for $n\in\nset{1,n_{\So}}$ we have
  $\exp\big(\tfrac{-\Vcst[]\vc}{1\vee400\VnormInf{\So}}\oDiSo^{1/2}\big)\leq1\leq
  n\oRaSo\leq n_{\So}\oRaSo$. Thereby, from \eqref{ag:ub:pnp:e2} with
  $\oRaSo=\min_{\Di\in\nset{1,\maxDi}}\dRaSo$
  follows  the claim, which completes the proof.\proEnd\end{pro}
\begin{rem}\label{ag:ub:pnp:il}Let us briefly comment on the last
  results.
  In case \ref{ag:ub:pnp:p}  the data-driven aggregation leads
  to an estimator attaining the parametric oracle rate (see
  \cref{oo:rem:ora}).   On the other hand in case
  \ref{ag:ub:pnp:np} the data-driven aggregation leads
  to an estimator attaining the oracle  rate $\oRa$ (see
  \cref{oo:rem:ora}), if the additional assumption
  \ref{ag:ub2:pnp:npc} is satisfied. Otherwise, the upper bound $\mRaSo$ in \eqref{ag:ub:pnp:e2}
  faces a deterioration compared to the rate $\oRaSo$. Considering again the  particular
  specification  $\bias^2(\So)\sim \Di^{-\wSob}$, $\wSob>0$, we have seen that $\oDiSo\sim n^{2/(2\wSob+1)}$ and $\oRaSo\sim n^{-2\wSob/(2\wSob+1)}$ is 
 the oracle dimension and rate, respectively. Obviously, in this situation the additional assumption
  \ref{ag:ub2:pnp:npc} is satisfied for any $\wSob>0$. Thereby, the
  data-driven aggregated estimator attains the oracle rate for all
  $\wSob>0$ and thus  it is adaptive. However, if $\bias^2(\So)\sim
  \exp(-\Di^{2\wSob})$, $\wSob>0$, then $\oDiSo\sim (\log n)^{1/\wSob}$  and $\oRaSo\sim (\log n)^{1/2\wSob}n^{-1}$ is 
 the oracle dimension and rate, respectively. In this situation the additional assumption
  \ref{ag:ub2:pnp:npc} is satisfied only for $\wSob\in(0,1/2]$. Hence,
  for $\wSob\in(0,1/2]$ the  data-driven aggregated estimator attains the oracle rate.
  In case $\wSob>1/2$, however, with
  $(\sDi{n})^{1/2}:=\tfrac{1\vee400\VnormInf{\So}}{\Vcst[]\vc}|\log\oRaSo|\sim(\log
  n)$ the upper bound in \eqref{ag:ub:pnp:e2} satisfies
  $\mRaSo\leq\dRa[\sDi{n}]\sim(\log
  n)n^{-1}$. Thereby, the rate $\mRaSo$ of the  data-driven estimator
  $\hSoPr[{\We[]}]$ features a deterioration at most by a logarithmic factor
  $(\log n)^{(1-1/(2\wSob))}$ compared to the oracle rate $\oRaSo$, i.e.
  $(\log n)n^{-1}$ versus
  $(\log n)^{1/(2\wSob)}n^{-1}$.\ilEnd
\end{rem}
\begin{pro}[Proof of \cref{ag:ub:pnp}]We outline here the main steps
  of the proof, while more technical details are deferred to the \cref{a:ag}.  Given penalties as in \eqref{ag:de:pen} for $\Di\in\Nz$ holds by construction
  \begin{equation}\label{ag:ass:pen:oo:c}
    \dRaSo\geq\bias^2(\xdf)\quad\text{and}\quad\cpen[\hvc\vee\vc]\dRaSo\geq[\pen\vee\hpen]
    \quad\text{for all }\Di\in\Nz.
  \end{equation}
  For arbitrary $\pdDi,\mdDi\in\nset{1,\maxDi}$ (to be chosen
  suitable below) let us   define
  \begin{multline}\label{ag:de:*Di}
    \mDi:=\min\set{\Di\in\nset{1,\mdDi}: \VnormA{\So}^2\bias^2(\So)\leq
      \VnormA{\So}^2\bias[\mdDi]^2(\So)+6\pen[\mdDi]}\quad\text{and}\\
    \pDi:=\max\set{\Di\in\nset{\pdDi,\maxDi}:
      \hpen \leq 6\VnormA{\So}^2\bias[\pdDi]^2(\So)+ 4\hpen[\pdDi]}
  \end{multline}
  where the defining set obviously contains $\mdDi$ and $\pdDi$,
  respectively, and hence, it is not empty. Note that only $\pDi$ depends
  on the observations, and hence is random.  Consider further the event $\vEv:=\set{|\hvc-\vc|\leq\vc/2}$ and its
  complement $\vEvC$, where by construction  for all $\Di\in\Nz$ holds
  $\tfrac{1}{2}\pen\Ind{\vEv}\leq
  \hpen\Ind{\vEv}\leq\tfrac{3}{2}\pen$.
  Exploiting the last bounds and \eqref{ag:de:*Di} it follows
  \begin{multline*}
    \pen[\pDi]\Ind{\vEv}\leq2\hpen[\pDi]\Ind{\vEv}\leq
    2\big(6\VnormA{\So}^2\bias[\pdDi]^2(\So)+
    4\hpen[\pdDi]\big)\Ind{\vEv}\\\leq2\big(6\VnormA{\So}^2\bias[\pdDi]^2(\So)+
    4\tfrac{3}{2}\pen[\pdDi]\big)=12\big(\VnormA{\So}^2\bias[\pdDi]^2(\So)+
    \pen[\pdDi]\big)
  \end{multline*}
  and with $\pen\leq\pen[\maxDi]\leq \cpen\vc$ for all
  $\Di\in\nset{1,\maxDi}$ also
  \begin{multline*}
    \sum_{\Di=1+\pDi}^{\maxDi}\We\pen\Ind{\{\VnormA{\hSoPr-\SoPr}^2\geq\hpen/7\}}
    \leq \cpen\vc\Ind{\vEvC} +  
    \sum_{\Di=1}^{\maxDi}\pen\Ind{\{\VnormA{\hSoPr-\SoPr}^2\geq\pen/14\}}.
  \end{multline*}
  Combining the last bounds and \cref{co:agg} we obtain
  \begin{multline}\label{co:agg:e2}
    \VnormA{\hSoPr[{\We[]}]-\So}^2\leq
    \tfrac{2}{14}\big(\cpen\vc\Ind{\vEvC} +12\VnormA{\So}^2\bias[\pdDi]^2(\So)+
12      \pen[\pdDi]\big) +2\VnormA{\So}^2\bias[\mDi]^2(\So)\\\hfill
+2\sum_{\Di=1}^{\maxDi}\vectp{\VnormA{\hSoPr-\SoPr}^2-\pen/14}  
+\tfrac{2}{14}\big(\cpen\vc\Ind{\vEvC} +\sum_{\Di=1}^{\maxDi}\pen\Ind{\{\VnormA{\hSoPr-\SoPr}^2\geq\pen/14\}}\big)\\
+2\VnormA{\So}^2\big(\Ind{\vEvC}+\FuVg{\We[]}(\nsetro{1,\mDi})\Ind{\vEv}\big)+\tfrac{2}{14}\sum_{\Di=1+\pDi}^{\maxDi}\pen\We\Ind{\{\VnormA{\hSoPr-\SoPr}^2<\hpen/7\}}.
\end{multline}
We bound the last two terms on the right hand side 
considering \bwn  $\We[]:=\rWe[]$
  as in \eqref{ag:de:rWe} and \mwn $\We[]:=\msWe[]$
  as in \eqref{ag:de:msWe} in \cref{ag:re:SrWe} and
  \cref{ag:re:SrWe:ms}, respectively. 
  Combining
those upper bounds and \eqref{co:agg:e2} we obtain
    \begin{multline}\label{co:agg:e3}
      \VnormA{\hSoPr[{\We[]}]-\So}^2\leq
      \tfrac{12}{7}\VnormA{\So}^2\bias[\pdDi]^2(\So)+ \tfrac{12}{7}\pen[\pdDi]
      +2\VnormA{\So}^2\bias[\mDi]^2(\So)\\\hfill
      +\tfrac{8\VnormA{\So}^2}{\rWc^2\cpen^2\vc^2}\Ind{\{\mDi>1\}}
    \exp\big(-\tfrac{3\rWc\cpen\vc}{28}(\mdDi)^{1/2}\big)
    \big)
    + n^{-1}\tfrac{192\vc}{\rWc^{3}\cpen^2}   \\\hfill
      +2\sum_{\Di=1}^{\maxDi}\vectp{\VnormA{\hSoPr-\SoPr}^2-\pen/14}
      +\tfrac{2}{14}\sum_{\Di=1}^{\maxDi}\pen\Ind{\{\VnormA{\hSoPr-\SoPr}^2\geq\pen/14\}}\\
  + \big(\tfrac{2}{7}\cpen\vc+2\VnormA{\So}^2\big)\Ind{\vEvC}+ 2\VnormA{\So}^2\Ind{\{\mDi>1\}}\Ind{\setB{\VnormA{\hSoPr[\mdDi]-\SoPr[\mdDi]}^2\geq\pen[\mdDi]/14}}.
\end{multline}
The deviations of the last four terms on the right hand side in  \eqref{co:agg:e3}  we bound in  \cref{ag:re:nd:rest} by exploiting usual
concentration inequalities. Thereby, with $\cpen\geq84\Vcst[]\geq1$, $\Vcst[]$ as in \eqref{de:bas:cst},
$\tfrac{3\rWc\cpen}{28}\geq 1$ and $\rWc\cpen\vc\geq1$ combining
\eqref{co:agg:e3} and \cref{ag:re:nd:rest}  there is a
  finite numerical constant $\cst{}>0$ such that for any sample size
  $n\in\Nz$, $n\geq3$, any dimension parameter  $\mdDi,\pdDi\in\nset{1,\maxDi}$ and associated $\mDi\in\nset{1,\maxDi}$  as defined in \eqref{ag:de:*Di} hold
  \begin{multline}\label{co:agg:e4}
    \nEx\VnormA{\hSoPr[{\We[]}]-\So}^2\leq
          \tfrac{12}{7}\VnormA{\So}^2\bias[\pdDi]^2(\So)+ \tfrac{12}{7}\pen[\pdDi]
          +2\VnormA{\So}^2\bias[\mDi]^2(\So)\\\hfill
                 + \cst{}\VnormA{\So}^2\Ind{\{\mDi>1\}} \exp\big(\tfrac{-\Vcst\vc}{1\vee400\VnormInf{\So}}(\mdDi)^{1/2}\big)\\
       +\cst{}\big(\VnormInf{\So}^3\vee1\big)\,\vc \,n^{-1}
      + \cst{}\big(\VnormA{\So}^2\vee1\big)\,  \FuEx{\So}(|X_1|^{2\aSob})\, n^{-1} 
  \end{multline}
 We distinguish now the two cases \ref{ag:ub:pnp:p} and
\ref{ag:ub:pnp:np} given in \eqref{ag:ub:pnp}. The tedious case-by-case analysis for
\ref{ag:ub:pnp:p} we defer to \cref{ag:ub:p} in the appendix. Here  we consider
the case \ref{ag:ub:pnp:np}  only.  However, in both cases the proof
is based on  an evaluation of the   upper bound 
\eqref{co:agg:e4}  for a suitable selection of the parameters
$\mdDi,\pdDi\in\nset{1,\maxDi}$. Recall the definition of the
oracle dimension and rate
$\oDi:=\oDiSo\in\nset{1,n^2}$
and $\dRa:=\dRaSo$, respectively.  We   select
$\pdDi:=\argmin\Nset[\Di\in\nset{1,\maxDi}]{\dRa}$. Further 
the inequalities \eqref{ag:ass:pen:oo:c} and 
the definition
\eqref{ag:de:*Di} of  $\mDi$ implies
$\VnormA{\So}^2\bias[\mDi]^2(\So)\leq
       (\VnormA{\So}^2+6\cpen\vc)\dRa[\mdDi]$. Keeping the last bound  together with 
 $\dRa[\mdDi]\geq\dRa[\pdDi]\geq\oRa\geq
 n^{-1}$, which holds for all $\mdDi\in\nset{1,\maxDi}$, in mind, we evaluate the upper bound \eqref{co:agg:e4} and obtain
 the assertion
 \eqref{ag:ub:pnp:e2}, 
which  completes the
proof of \cref{ag:ub:pnp}.\proEnd
\end{pro}
  \paragraph{Maximal risk bounds.}
  The following assertion shows that the aggregation
  $\hSoPr[{\We[]}]$ using either \bwn $\We[]:=\rWe[]$
  as in \eqref{ag:de:rWe} or \mwn $\We[]:=\msWe[]$
  as in \eqref{ag:de:msWe} attains the minimax optimal rate over
  Sobolev-ellipsoids $\rwcSobD{\wSob,\rSob,2}$ as in \eqref{equation:dens_sobol}.
\begin{thm}\label{ag:mub:sob}  
  Consider an aggregation $\hSoPr[{\We[]}]$ 
  using either \bwn
  $\We[]:=\rWe[]$ as in \eqref{ag:de:rWe} or \mwn  $\We[]:=\msWe[]$ as in \eqref{ag:de:msWe} with $\Delta \geq 84\Vcst[]$ and $\Vcst[]$ as in \eqref{de:bas:cst}. For each
  $\wSob,\rSob\in\pRz$ with
  $\wSob>1$ there is a finite constant $\cst{\wSob,\rSob}$ depending only on the class
  $\srwcSo[\wSob]$  such that for
  all $n\geq3$ holds $ \mmRi[{\We[]}]{\wSob}{2}
  \leq \cst{\wSob,\rSob}\; n^{-2s/(2s+1)}.$
\end{thm}
\begin{te}
  The proof of \cref{ag:mub:sob} follows a long the lines of the proof
  of case \ref{ag:ub:pnp:np} in \cref{ag:ub:pnp} where we did not
  specify the asymptotic behaviour of the sequence
  $\Ksuite[\Di\in\Nz]{\bias^2(\So)}$. Therefore, rather imposing a
  specific polynomial decay as implied by a Sobolev ellipsoid we
  characterise it by a strictly positive sequence
  $\wSo=\Ksuite[\Di\in\Nz]{\wSoDi}$. Precisely, let $
  \rwcSob{\wSo}:=\{f \in \LpA[2] : \Vabs[\wSo]{f}^2 := \sum_{\Di\in
    \Nz} (|\fSo{\Di}|^2/ \wSoDi)< \infty\} \text{ and }
  \rwcSob{\wSo,\rSob}:= \{f \in \rwcSob{\wSo}: \Vabs[\wSo]{f}^2\leq
  \rSob\}$.  Obviously, the Sobolev ellipsoid
  $\rwcSob{\wSob,\rSob}$ corresponds to the special case
  $\wSo=\Nsuite[\Di]{\Di^{-\wSob}}$.  Keeping \eqref{equation:dens_sobol} in mind we denote
  further a corresponding subset of densities with finite
  $\aSob\mSob$-th moment,
  $\mSob\in\pRz$, by $\rwcSobD{\wSo,\rSob,\mSob}:= \{f\in
  \rwcSob{\wSo,\rSob}: f \text{ is a density and }
  \FuEx{\So}\big(|X|^{\aSob\mSob}\big) \leq
  L\}$. For $n,\Di\in\Nz$ we set
  \begin{multline}\label{de:mmRa}
    \dRaC:=[\wSoDi\vee n^{-1}\Di^{1/2}],\quad 
    \mmDiC:=\argmin\Nset[\Di\in\Nz]{\dRaC}\quad\text{and}\\
    \mmRaC :=\min\Nset[\Di\in\Nz]{\dRaC}.\hfill\end{multline}
 Here and subsequently, we impose the following minimal
  regularity conditions.
\end{te}
\setcounter{ass}{1}
\begin{ass}\label{ass:Cw}The sequence 
  $\wSo=\Ksuite[\Di\in\Nz]{\wSoDi}$ is strictly positive, monotonically non-increasing with
  $\wSoDi[1]\leq1$,
  $\lim_{\Di\to\infty}\wSoDi=0$ and there is
  $\Vcst[\wSo,\rSob]\in\pRz$ such that
  $\db1\leq\rSob\vee\vc\vee\VnormA{\So}^2\vee\VnormInf{\So}^2\leq\Vcst[\wSo,\rSob]$ for all
  $\So\in\rwcSobD{\wSo,\rSob,2}$.
  \end{ass}%
\begin{rem}\label{rem:ass:Cw}
  We shall emphasise that for any
  $\So\in\rwcSobD{\wSo,\rSob,\mSob}$ hold $\vc\leq 1+\rSob$ and 
  $\VnormA{\So}^2\bias^2(\So)\leq\rSob\wSoDi$ for all
  $\Di\in\Nz$. Keeping further $\sup_{j\in\Nz_0}\VnormInf{\bas_j}\leq
  \sqrt{2}$ in mind we have $\fSo{0}\leq \VnormInf{\bas_0}\leq
  \sqrt{2}$ and hence $\VnormA{\So}^2\leq
  2+\rSob$. Moreover, if 
  $\Vnormlp[{\lp[1]}]{\wSo}:=\sum_{\Di\in\Nz}\wSoDi<\infty$, then
  uniformly for all
  $\So\in\rwcSobD{\wSo,\rSob,\mSob}$ we have $\VnormInf{\So}^2\leq
  (2+\sum_{\Di\in\Nz}(|\fSo{\Di}|^2/\wSoDi))\VnormInf{\bas_{0}^2+
    \sum_{\Di\in\Nz}\wSoDi\bas_{\Di}^2}\leq(2+\rSo)2(1+\Vnormlp[{\lp[1]}]{\wSo})$
  by applying the Cauchy-Schwarz inequality. Consequently, if
  $\Vnormlp[{\lp[1]}]{\wSo}<\infty$ then $\rSob\vee\vc\vee\VnormA{\So}^2\vee\VnormInf{\So}^2\leq(2+\rSo)2(1+\Vnormlp[{\lp[1]}]{\wSo})
  =:\Vcst[\wSo,\rSob]$  for all
  $\So\in\rwcSobD{\wSo,\rSob,\mSob}$. In particular, the Sobolev-ellipsoid $\srwcSo[\wSob]$  satisfies \cref{ass:Cw} for all
  $\wSob>1$. Note that, under \cref{ass:Cw} hold
  $\mmRaC=\dRaC[\mmDiC]$ and
  $\mmDiC\in\nset{1,n^2}$. Moreover,  we have $\mmRaC\geq
  n^{-1}$, $\mmRaC=o(1)$ and $n\mmRaC\to\infty$ as
  $n\to\infty$. In this situation the rate
  $\mmRaC$ is non-parametric and for any
  $\So\in\rwcSobD{\wSo,\rSob,\mSob}$ holds by construction
  $(2+\rSo)\mmRaC\geq\VnormA{\So}^2\oRaSo$ for all $n\in\Nz$. \remEnd
\end{rem}
\begin{te}
  Exploiting again the identity \eqref{oo:rdec}, the upper bound \eqref{de:bas:cst} and the definition
  \eqref{de:mmRa} under \cref{ass:Cw} there is a numerical constant
  such that for all $n\in\Nz$
  \begin{equation}\label{mm:e1}
    \mmRi[\mmDiC]{\wSo}{2}\leq \cst{}\,\Vcst[\wSo,\rSob]\,\mmRaC.
  \end{equation}
  By applying \cref{co:agg} we derive next bounds for the
  maximal risk over ellipsoids
  $\rwcSobD{\wSo,\rSob,2}$ of the aggregation
  $\hSoPr[{\We[]}]$ using either \bwn
  $\We[]:=\rWe[]$ as in \eqref{ag:de:rWe} or \mwn
  $\We[]:=\msWe[]$ as in \eqref{ag:de:msWe} based on the penalties $(\hpen)_{\Di\in\maxDi}$
  given in \eqref{ag:de:pen}.
\end{te}
\begin{thm}\label{ag:ub:mm} Consider  an aggregation $\hSoPr[{\We[]}]$
  using either \bwn $\We[]:=\rWe[]$
  as in \eqref{ag:de:rWe} or \mwn $\We[]:=\msWe[]$
  as in \eqref{ag:de:msWe}. Under \cref{ass:Cw} there is a  numerical constant $\cst{}$ such that for all
$n\geq3$  holds 
\begin{multline}\label{ag:ub:mm:e1}
       \mmRi{\wSo}{2} \leq   \cst{}\;\big(\Vcst[\wSo,\rSob]\;\mRaC
       +\Vcst[\wSo,\rSob]^3\;n^{-1}\big)\\
    \text{ with }\mRaC:=\min_{\Di\in\nset{1,\maxDi}}\setB{\big[\dRaC\vee\exp\big(-\tfrac{\Vcst[]\Di^{1/2}}{400\Vcst[\wSo,\rSob]}\big)\big]}.
\end{multline}
\end{thm}
\begin{te}
Before we proof the main result let us state an immediate consequence.  
\end{te}
\begin{co}\label{ag:ub:mm2}
  Let the assumptions of \cref{ag:ub:mm} be satisfied. If  in addition
    \begin{inparaenum}\item[\mylabel{ag:ub:mm2:c}{{\dgrau\upshape\bfseries(A1')}}]
      there is  $n_{\wSo,\rSob}\in\Nz$ such that $\mmDiC$  and $\mmRaC$  as
      in \eqref{de:mmRa} satisfy  $\maxDi^{1/2}\geq(\mmDi)^{1/2}\geq\tfrac{400\Vcst[\wSo,\rSob]}{\Vcst}|\log\mmRaC|$ for all $n\geq n_{\wSo}$,
    \end{inparaenum}
    then there is a numerical constant $\cst{}>0$ such that  for all
    $n\in\Nz$ holds $\mmRi{\wSo}{2}\leq\cst{}\;\big([n_{\wSo,\rSob}\vee\Vcst[\wSo,\rSob]]\;\mmRaC
       +\Vcst[\wSo,\rSob]^3\;n^{-1}\big)$.
\end{co}
\begin{pro}[Proof of \cref{ag:ub:mm2}] follows in analogy to
  \cref{ak:ag:ub2:pnp} and we omit the details.\proEnd\end{pro}
\begin{rem}\label{ag:ub:mm:rem}
  Let us briefly comment on the last
  results. The data-driven aggregation leads
  to an estimator attaining the rate $\mmRaC$ due to \cref{ag:ub:mm2}, if the additional assumption
  \ref{ag:ub:mm2:c} is satisfied. Otherwise, the upper bound $\mRaC$  in \eqref{ag:ub:mm:e1}
  faces a deterioration compared to the rate $\mmRaC$. Considering the Sobolev ellipsoid $\rwcSobD{\wSo,\rSob,2}$, i.e., $\wSo=\Nsuite[\Di]{\Di^{-s}}$, $s\in\pRz$, where
  $\mmDiC\sim n^{2/(2s+1)}$ and 
  $\mmRaC\sim n^{-2s/(2s+1)}$, \cref{ass:Cw} and
  \ref{ag:ub:mm2:c} are satisfied for each $s>1$.  Consequently, \cref{ag:mub:sob}
  follows immediately from \cref{ag:ub:mm2}. On the other hand if 
  $\wSo=\Nsuite[\Di]{\exp(-\Di^{2s})}$, $s\in\pRz$, then $\mmDiC\sim (\log n)^{1/\wSob}$
  and $\mmRaC\sim (\log n)^{1/2\wSob}n^{-1}$.
  In this situation the additional assumption
  \ref{ag:ub:mm2:c} is satisfied only for $\wSob\in(0,1/2]$. Hence,
  for $\wSob\in(0,1/2]$ the  data-driven aggregation attains
  the rate $\mmRaC$.  In case $\wSob>1/2$, however, with
  $(\sDi{n})^{1/2}:=\tfrac{400\Vcst[\wSo,\rSob]}{\Vcst}|\log\mmRaC|\sim(\log
  n)$ we have 
  $\mRaC\leq\dRa[\sDi{n}]\sim(\log
  n)n^{-1}$. Thereby, the rate $\mRaC$ of the  aggregation
  $\hSoPr[{\We[]}]$ features a deterioration at most by a logarithmic factor
  $(\log n)^{(1-1/(2\wSob))}$ compared to the rate $\mmRaC$, i.e.
  $(\log n)n^{-1}$ versus
  $(\log n)^{1/(2\wSob)}n^{-1}$.\remEnd
\end{rem}
\begin{pro}[Proof of \cref{ag:ub:mm}] 
  We make use of the upper bound
  \eqref{co:agg:e3} derived in the proof of \cref{ag:ub:pnp}. We note
  that uniformly for all $\So\in\rwcSobD{\wSo,\rSob,2}$  under \cref{ass:Cw}  the definition \eqref{ag:de:*Di} of $\pDi$  and $\mDi$ implies
  $\VnormA{\So}^2\bias[\pdDi]^2(\So)+\pen[\pdDi]  \leq(1+\cpen)\Vcst[\wSo,\rSob]\dRaC[\pdDi]$ and
  $\VnormA{\So}^2\bias[\mDi]^2(\So)
  \leq(1+6\cpen)\Vcst[\wSo,\rSob]\dRaC[\mdDi]$. Combining
  \eqref{co:agg:e3}, the last bounds,
  $\VnormA{\So}^2\vee\vc\vee\tfrac{8\VnormA{\So}^2}{\rWc^2\cpen^2\vc^2}\vee\tfrac{192\vc}{\rWc^{3}\cpen^2}\leq
  \Vcst[\wSo,\rSob]$ and $\tfrac{3\rWc\cpen\vc}{28}\geq\tfrac{\Vcst}{400\Vcst[\wSo,\rSob]}$  we obtain
 \begin{multline*}
      \VnormA{\hSoPr[{\We[]}]-\So}^2\leq
      \tfrac{12}{7}(1+\cpen)\Vcst[\wSo,\rSob]\dRaC[\pdDi]
      +2(1+6\cpen)\Vcst[\wSo,\rSob]\dRaC[\mdDi]\\\hfill
      +\Vcst[\wSo,\rSob]\exp\big(-\tfrac{\Vcst}{400\Vcst[\wSo,\rSob]}(\mdDi)^{1/2}\big)
    \big)
    + \Vcst[\wSo,\rSob]n^{-1}   \\\hfill
      +2\sum_{\Di=1}^{\maxDi}\vectp{\VnormA{\hSoPr-\SoPr}^2-\pen/14}
      +\tfrac{2}{14}\sum_{\Di=1}^{\maxDi}\pen\Ind{\{\VnormA{\hSoPr-\SoPr}^2\geq\pen/14\}}\\
  + \big(\tfrac{2}{7}\cpen+2\big)\Vcst[\wSo,\rSob]\Ind{\vEvC}+ 2\Vcst[\wSo,\rSob]\Ind{\setB{\VnormA{\hSoPr[\mdDi]-\SoPr[\mdDi]}^2\geq\pen[\mdDi]/14}}.
\end{multline*}
 In \cref{ag:re:2nd:rest} in \cref{a:prel} we bound
  the last four  terms uniformly for all
  $\So\in\rwcSobD{\wSo,\rSob,2}$. Therewith, there exists a finite numerical constant
  $\cst{}>0$ such that for all $n\in\Nz$
 \begin{multline}\label{ag:ub:mm:p1}
    \mmRi{\wSo}{2} \leq  \cst{}\,\Vcst[\wSo,\rSob]\big(\dRaC[\pdDi]
      +\dRaC[\mdDi]+\exp\big(\tfrac{-\Vcst}{400\Vcst[\wSo,\rSob]}(\mdDi)^{1/2}\big)\big)\\\hfill
      +\cst{}\,\Vcst[\wSo,\rSob]^3\,n^{-1}.
\end{multline}
For
$\mmDi:=\mmDiC\in\nset{1,n^2}$
and $\dRaC$ as in \eqref{de:mmRa} set
$\pdDi:=\argmin\Nset[\Di\in\nset{1,\maxDi}]{\dRaC}$, then for any
$\mdDi\in\nset{1,\maxDi}$ we have
 $\dRaC[\mdDi]\geq\dRaC[\pdDi]\geq\dRaC[\mmDi]=\min\Nset[\Di\in\Nz]{\dRaC}$.
 The last bounds  together with \eqref{ag:ub:mm:p1} imply \eqref{ag:ub:mm:e1}, which completes the proof of
\cref{ag:ub:mm}.\proEnd
\end{pro}


%% file: _4Numeric.tex
%
%
%
\section{Numerical study}\label{si}

Let us illustrate the performance of the aggregated
estimator $\hSoPr[{\We[]}]$ using \wname   $\We[]=\rWe$ (see
\eqref{ag:de:rWe}) or model selection weights $\We[]=\msWe$ (see
\eqref{ag:de:msWe}) 
in both cases \HerCase and \LagCase. 
For the case of \LagCase we consider the densities
\begin{resListeN}
	\item\label{si:lag:i} \textit{Gamma Mixture:} $f(x)=0.4 \cdot 3.2^2 x\exp(-3.2x) +0.6 \cdot \frac{6.8^{16} x^{15}}{15!}\exp(-6.8x)$ ,
	\item\label{si:lag:ii}  \textit{Gamma Distribution}: $f(x)=\frac{x^{4}}{4!}\exp(-x)$,
	\item\label{si:lag:iii}  \textit{Beta Distribution}: $f(x)=\frac{1}{560} (0.5x)^3(1-0.5x)^4 \1_{[0,1]}(0.5x)$ and
	\item\label{si:lag:iv}  \textit{Weibull Distribution}: $f(x)=0.75 x^{-0.25} \exp(-x^{0.75})$.
\end{resListeN}
In case of \HerCase we investigate the densities
\begin{resListeN}
	\item \textit{Gaussian Mixture}: $f(x) = 0.6 \cdot \frac{2.5}{\sqrt{2\pi}}  \exp(-x^2/0.8))+ 0.4\cdot \frac{2.5}{\sqrt{2\pi}}  \exp(-(x-3)^2/0.8))$ , 
	\item \textit{Finite representation}: $f(x)= \frac{1}{\sqrt{2\pi}} x^2 \exp(-x^2/2)$, 
	\item \textit{Beta Distribution:} $f(x)=\frac{1}{560} (0.5x)^3(1-0.5x)^4 \1_{[0,1]}(0.5x)$ and
	\item \textit{Pareto Distribution:}  $f(x)=\frac{0.75}{x^{1.75}} \1_{[1,\infty)}(x)$.
\end{resListeN}
We consider these four cases for the following reasons. The bias of
both densities in \ref{si:lag:i} has an exponential decay as shown by
\cite{BelomestnyComteGenon-Catalot2019} and
\cite{ComteGenon-Catalotothers2015}. The densities in \ref{si:lag:ii}
have a finite representation in the Laguerre respectively the Hermite
basis. The case \ref{si:lag:iii} and \ref{si:lag:iv} illustrates the behaviour of the
estimators when firstly the density has a compact support and secondly $\hvFu$ 
does not have a finite second moment.
By minimising an integrated squared error over a family of histogram
densities with randomly  drawn partitions  and weights we select
$\Delta = 1.02$  and $\Delta = 0.95$ in case \LagCase and \HerCase, respectively.
Furthermore, we chose $\kappa= 9.8$  and $\kappa = 5.2$ in case
\LagCase and \HerCase, respectively.\\[2ex] 
\begin{minipage}{\textwidth}
\centerline{\begin{minipage}[t]{0.32\textwidth}
		\includegraphics[width=\textwidth,height=40mm]{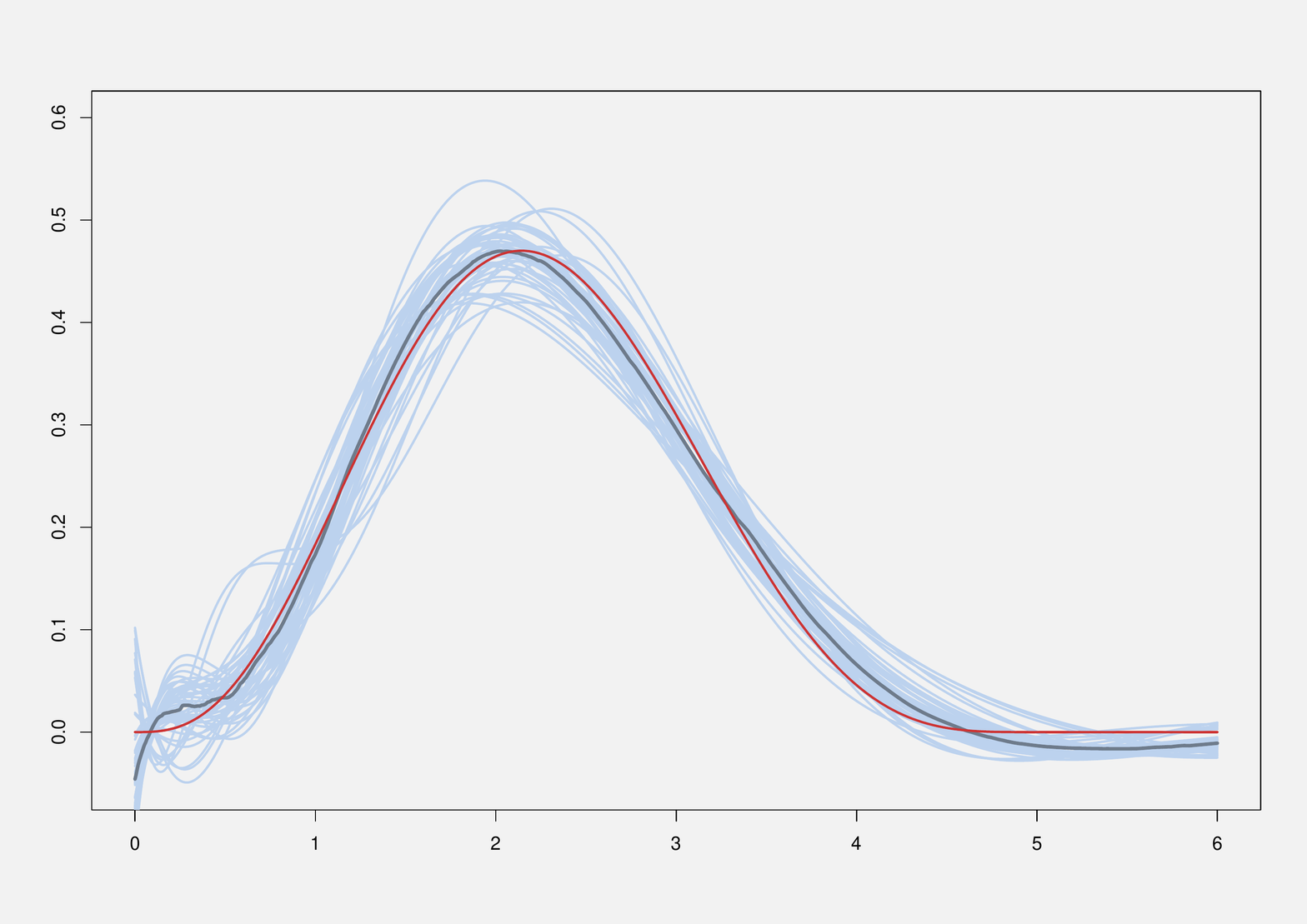}
	\end{minipage}
	\begin{minipage}[t]{0.32\textwidth}
		\includegraphics[width=\textwidth,height=40mm]{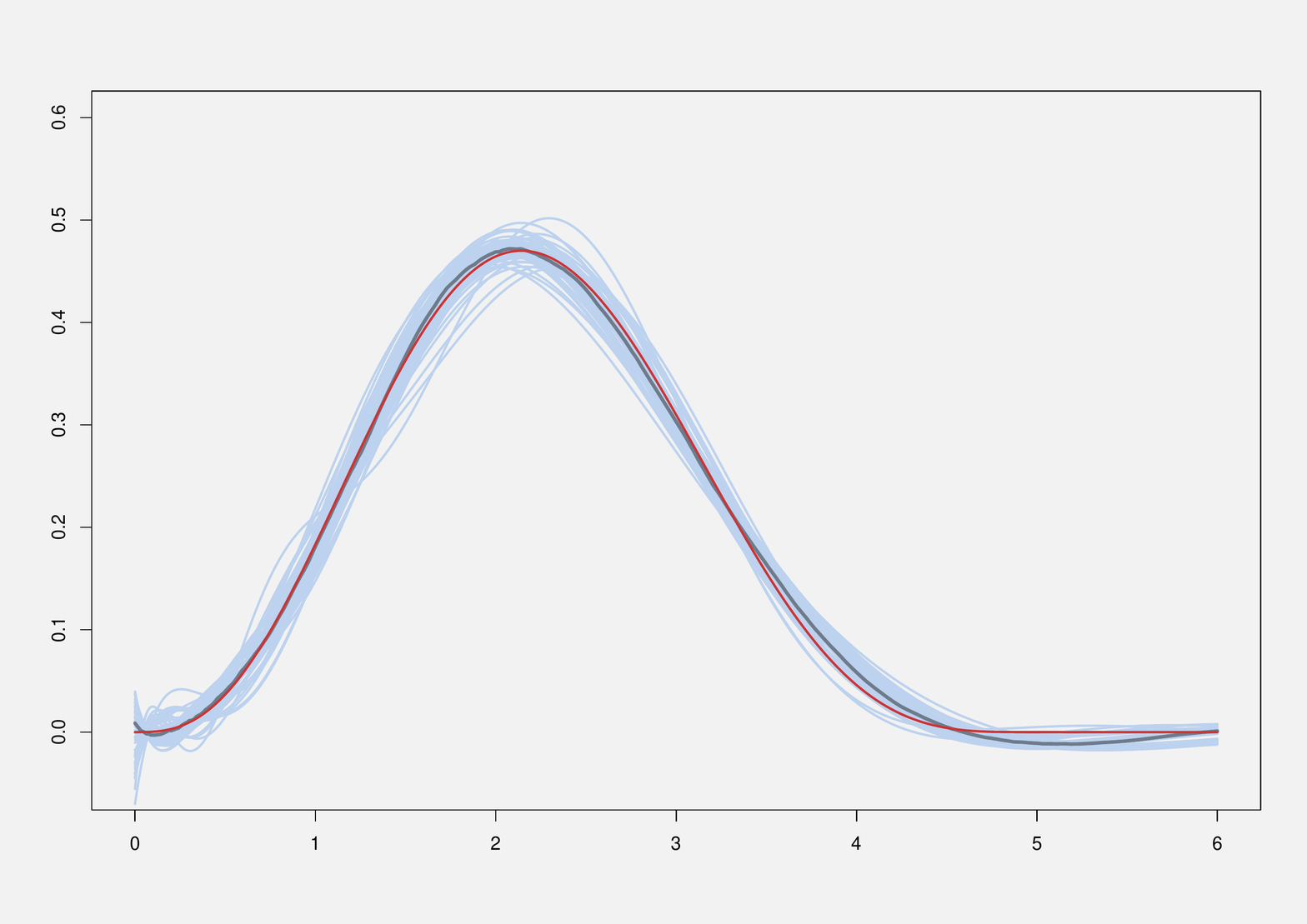}
	\end{minipage}
	\begin{minipage}[t]{0.32\textwidth}
		\includegraphics[width=\textwidth,height=40mm]{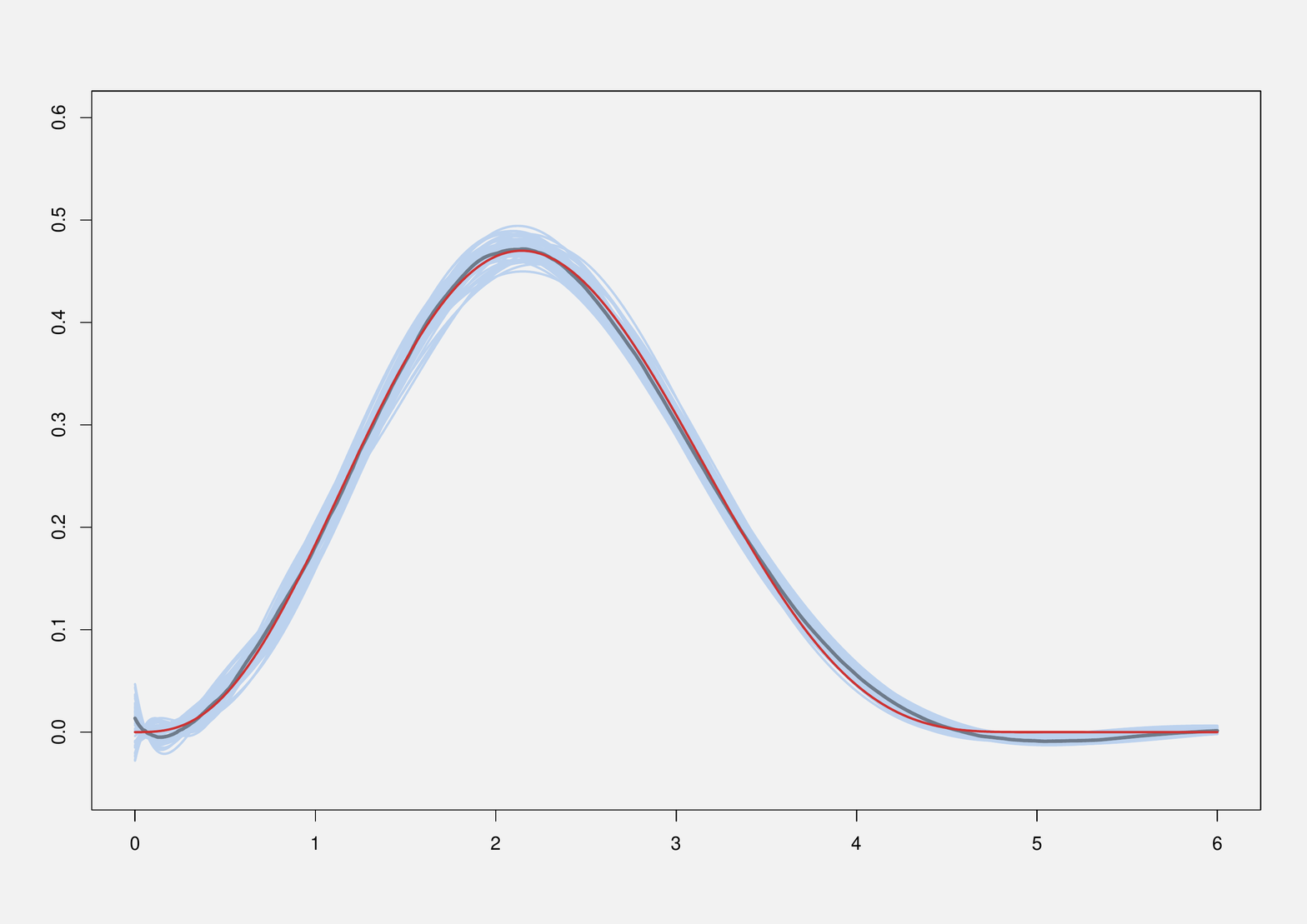}
	\end{minipage}}
\centerline{\begin{minipage}[t]{0.32\textwidth}
		\includegraphics[width=\textwidth,height=40mm]{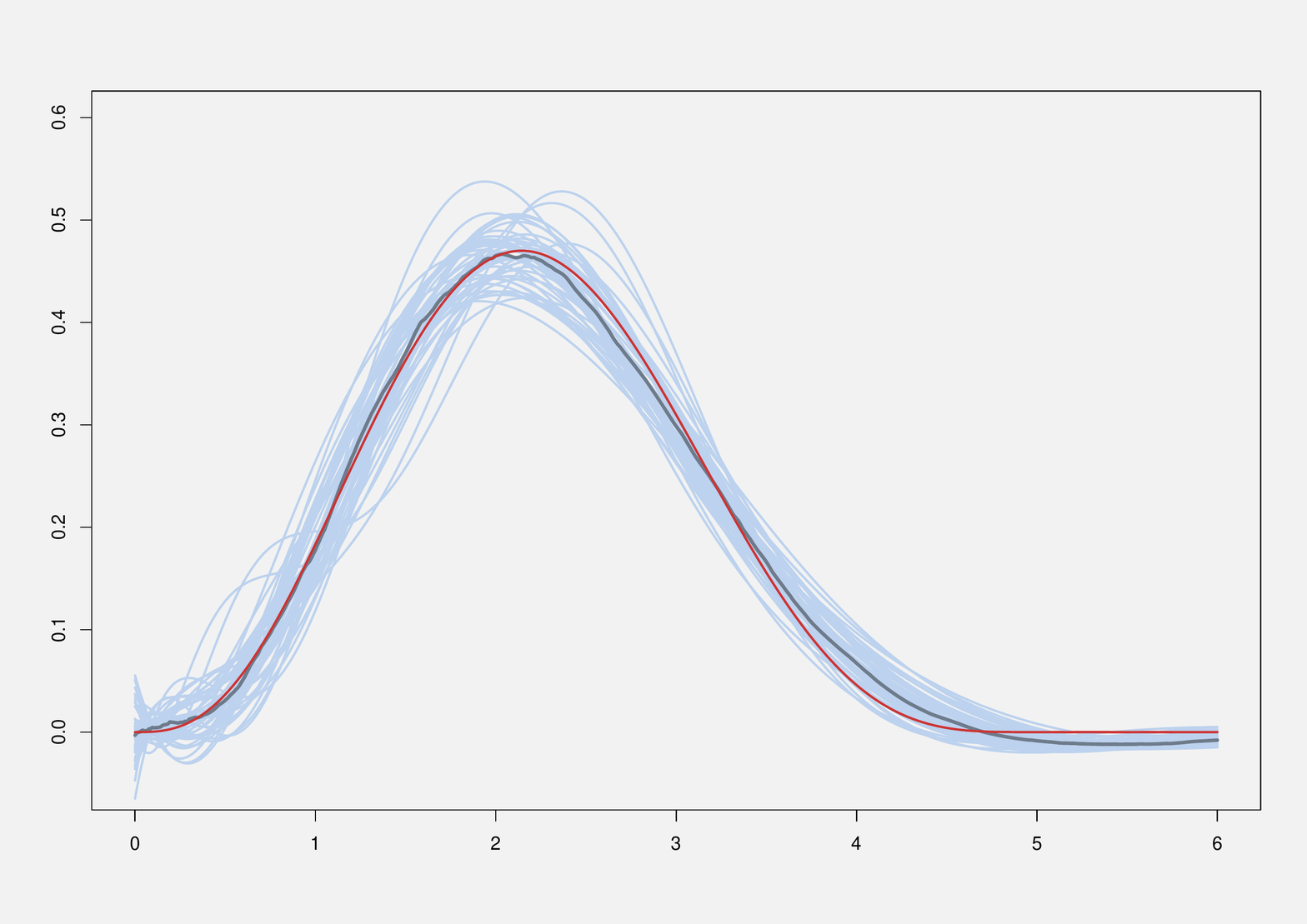}
	\end{minipage}
	\begin{minipage}[t]{0.32\textwidth}
		\includegraphics[width=\textwidth,height=40mm]{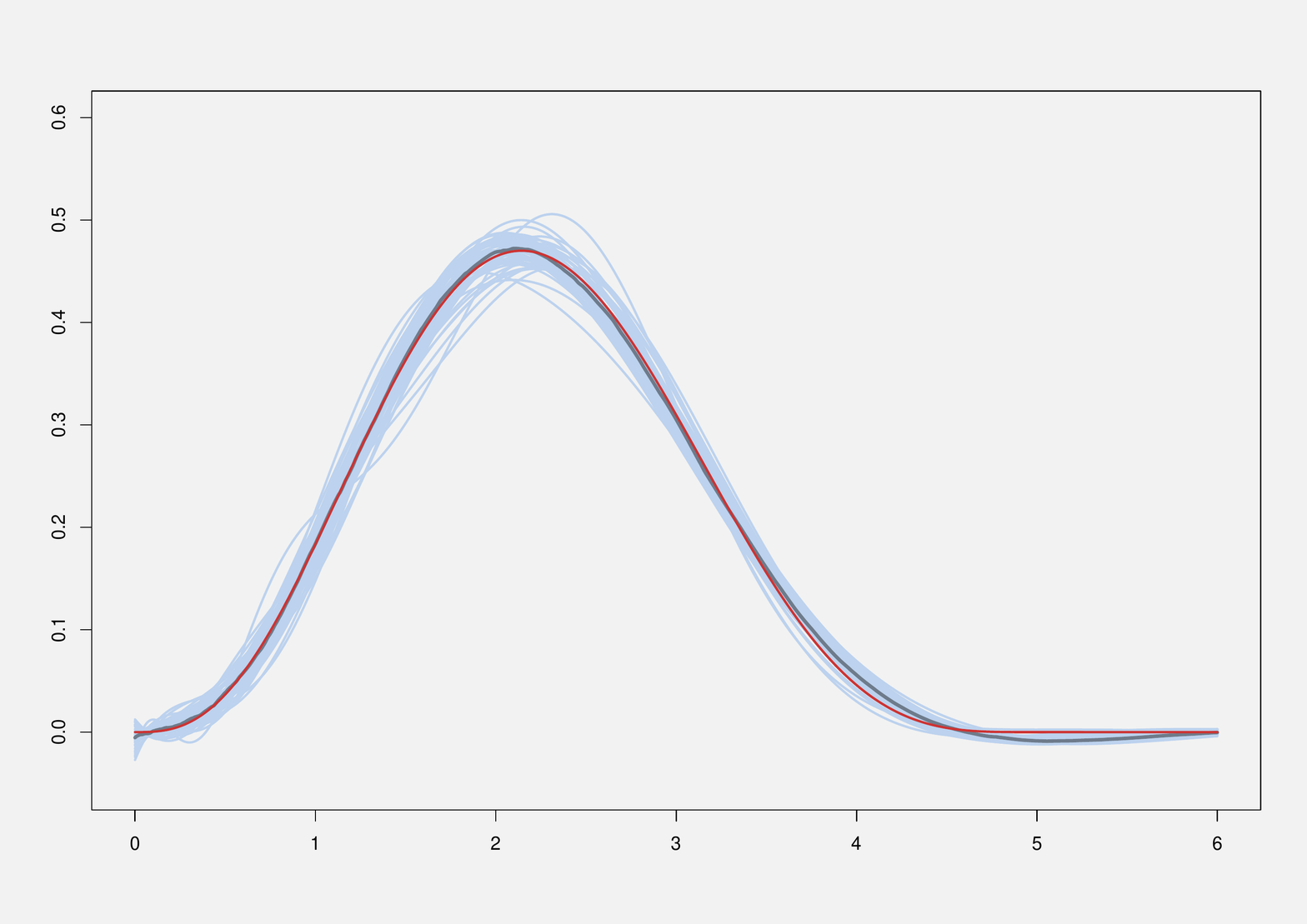}
	\end{minipage}
	\begin{minipage}[t]{0.32\textwidth}
		\includegraphics[width=\textwidth,height=40mm]{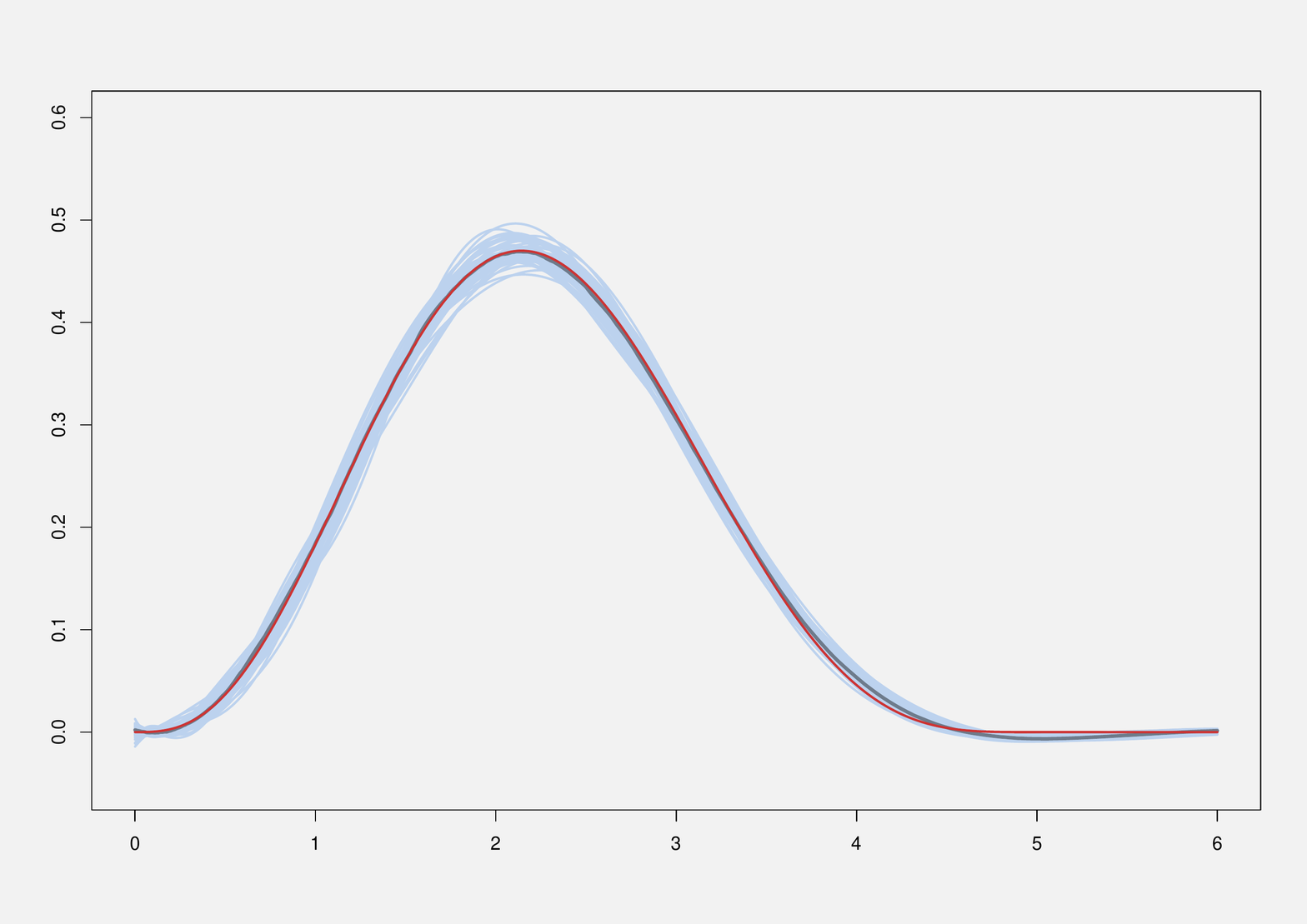}
	\end{minipage}}
      \captionof{figure}{Considering case \LagCase \ref{si:lag:iii}
        the estimators are depict for 
        50  Monte-Carlo simulations using  model selection weights
        (top) and  \wname (bottom) with varying sample size $n=200$ (left), $n=1000$ (middle) and $n=2000$ (right). The true density $\So$ is given by the red curve while the dark blue curve is the point-wise empirical median of the 50 estimates.}
\end{minipage}\\[2ex]
\begin{minipage}[t]{\textwidth}
\centerline{\begin{minipage}[t]{0.32\textwidth}
		\centering \includegraphics[width=\textwidth,height=40mm]{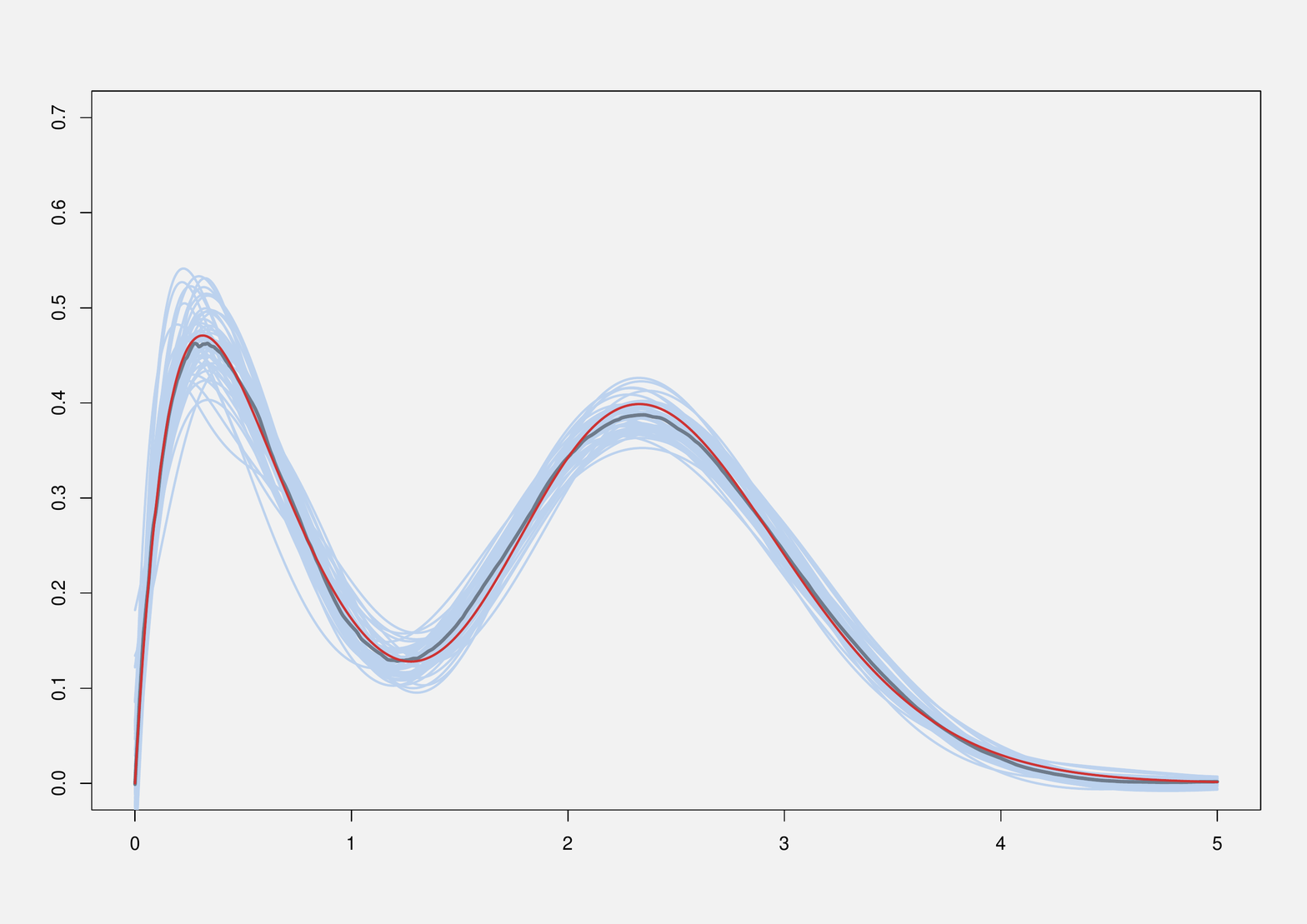}
	\end{minipage}
	\begin{minipage}[t]{0.32\textwidth}
	\centering	\includegraphics[width=\textwidth,height=40mm]{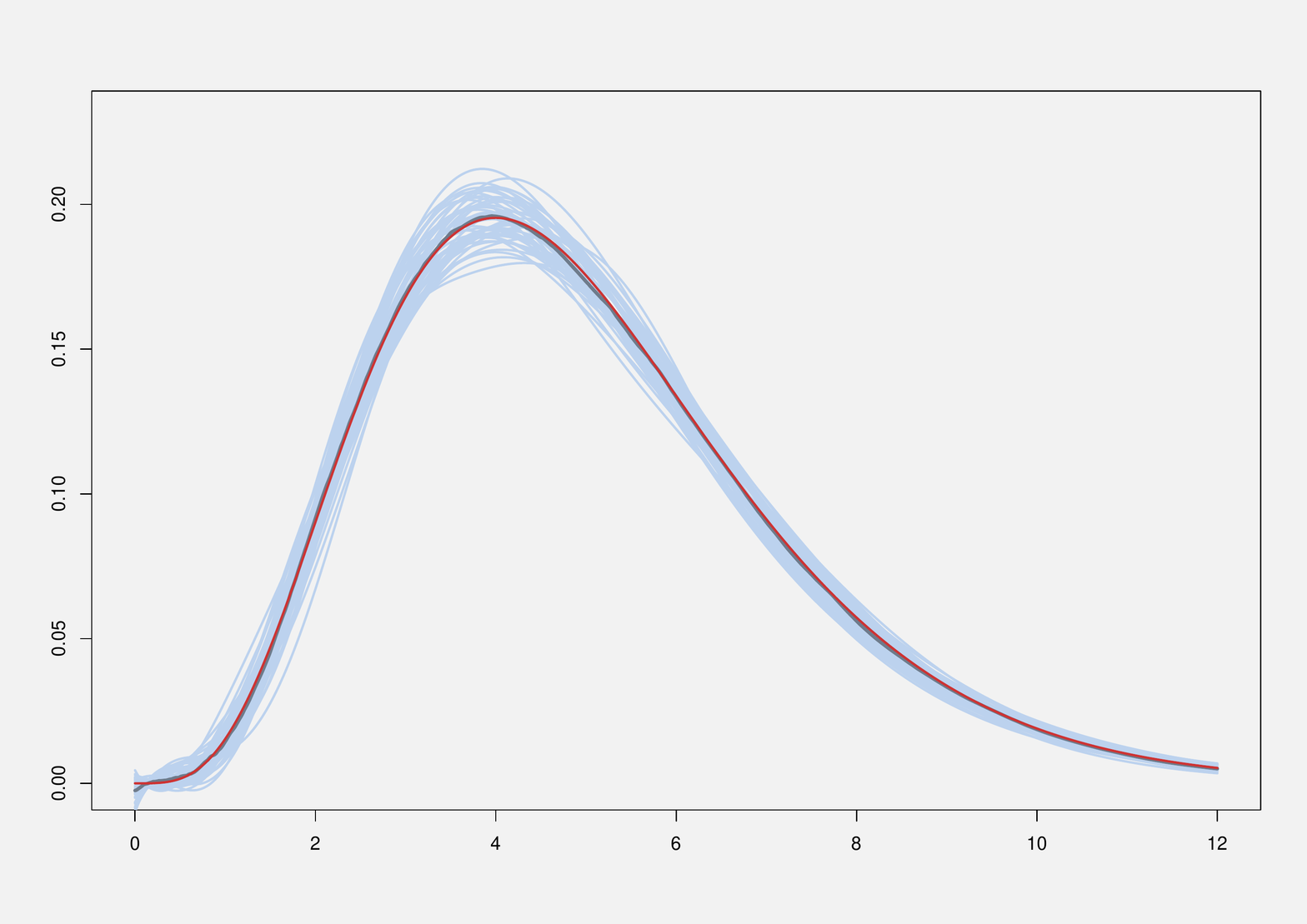}
	\end{minipage}
	\begin{minipage}[t]{0.32\textwidth}
		\centering	\includegraphics[width=\textwidth,height=40mm]{pics/Den5_BW_1000.pdf}
	\end{minipage}}
\centerline{\begin{minipage}[t]{0.32\textwidth}
	\centering		\includegraphics[width=\textwidth,height=40mm]{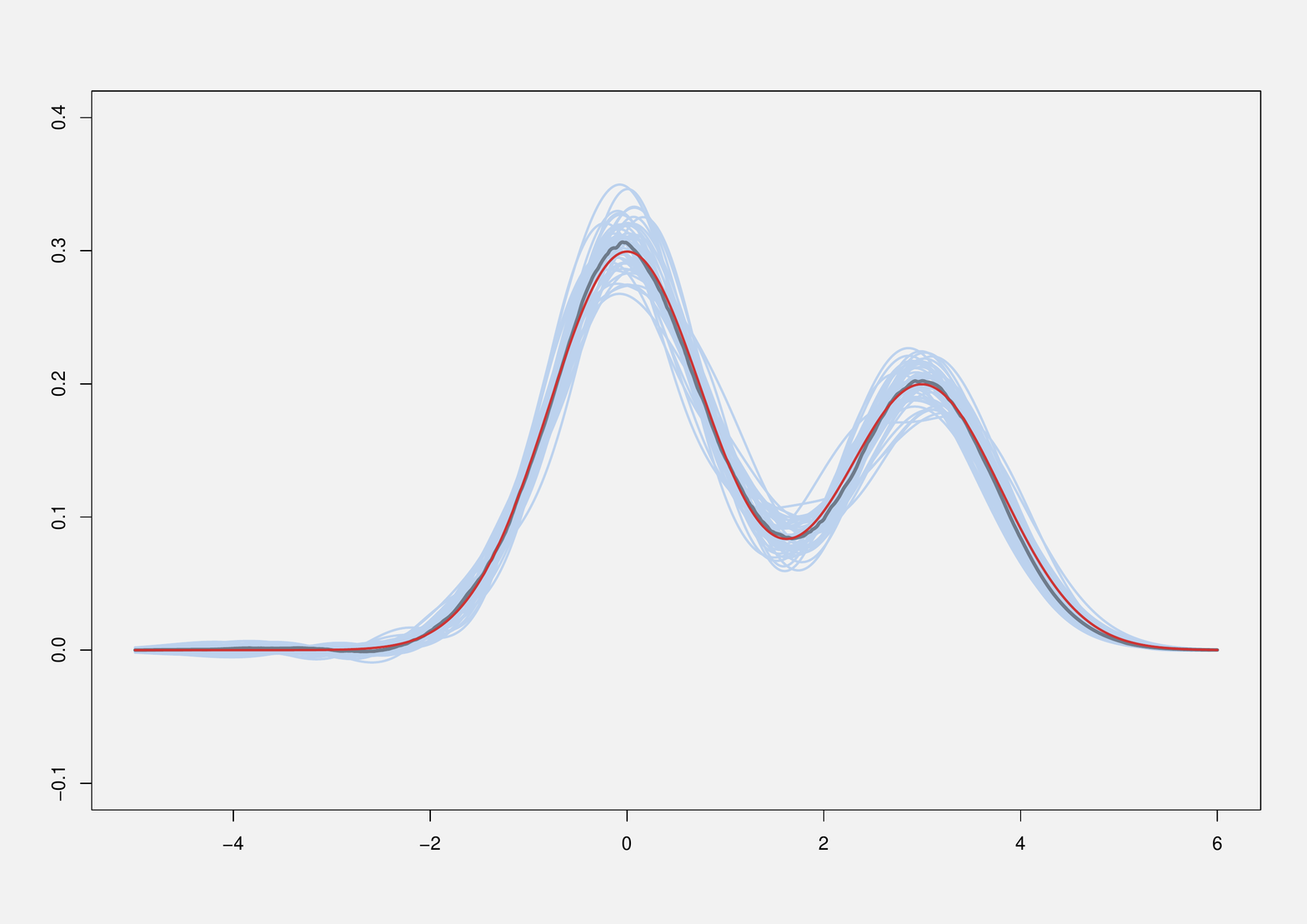}
	\end{minipage}
	\begin{minipage}[t]{0.32\textwidth}
	\centering		\includegraphics[width=\textwidth,height=40mm]{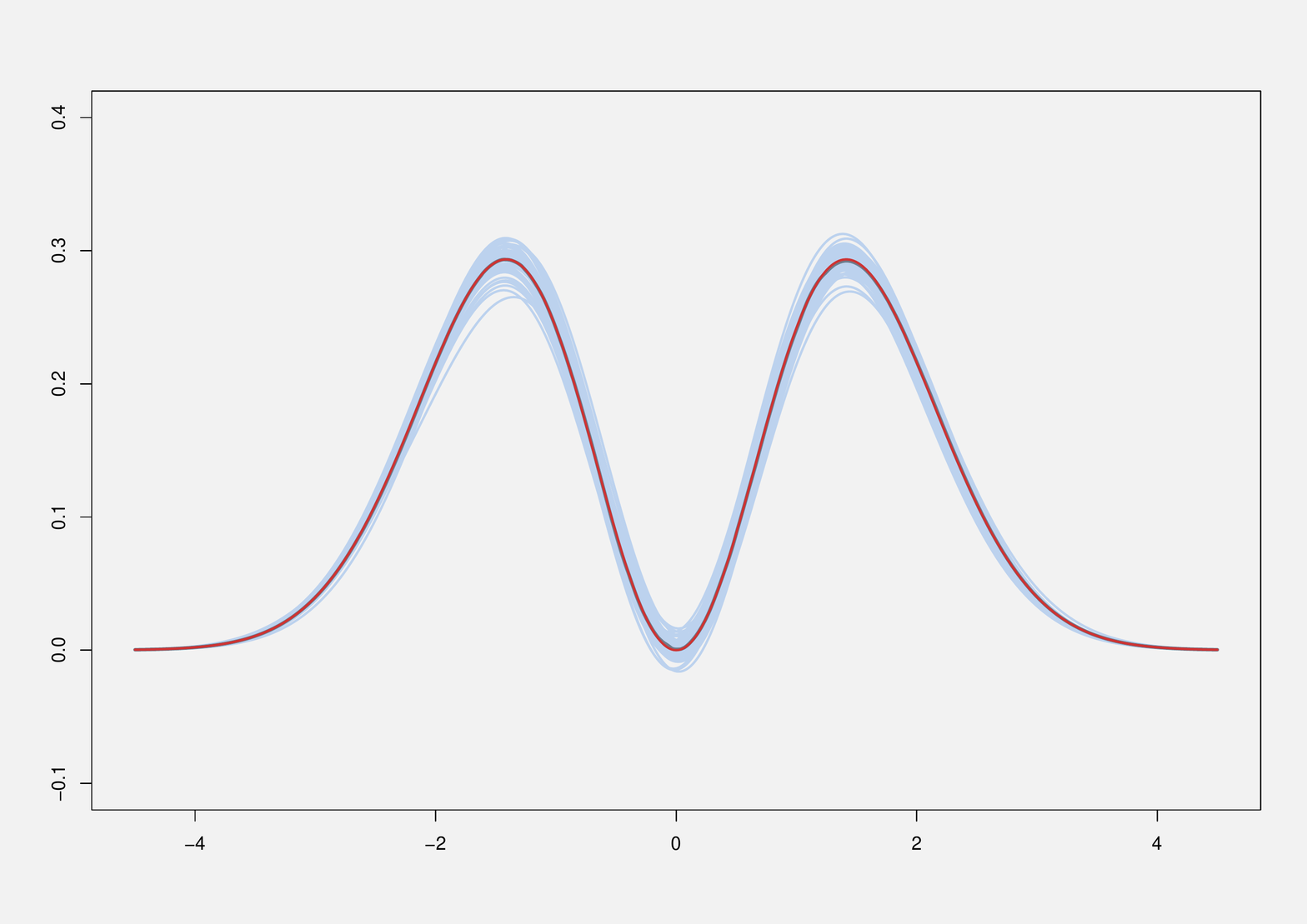}
	\end{minipage}
	\begin{minipage}[t]{0.32\textwidth}
	\centering		\includegraphics[width=\textwidth,height=40mm]{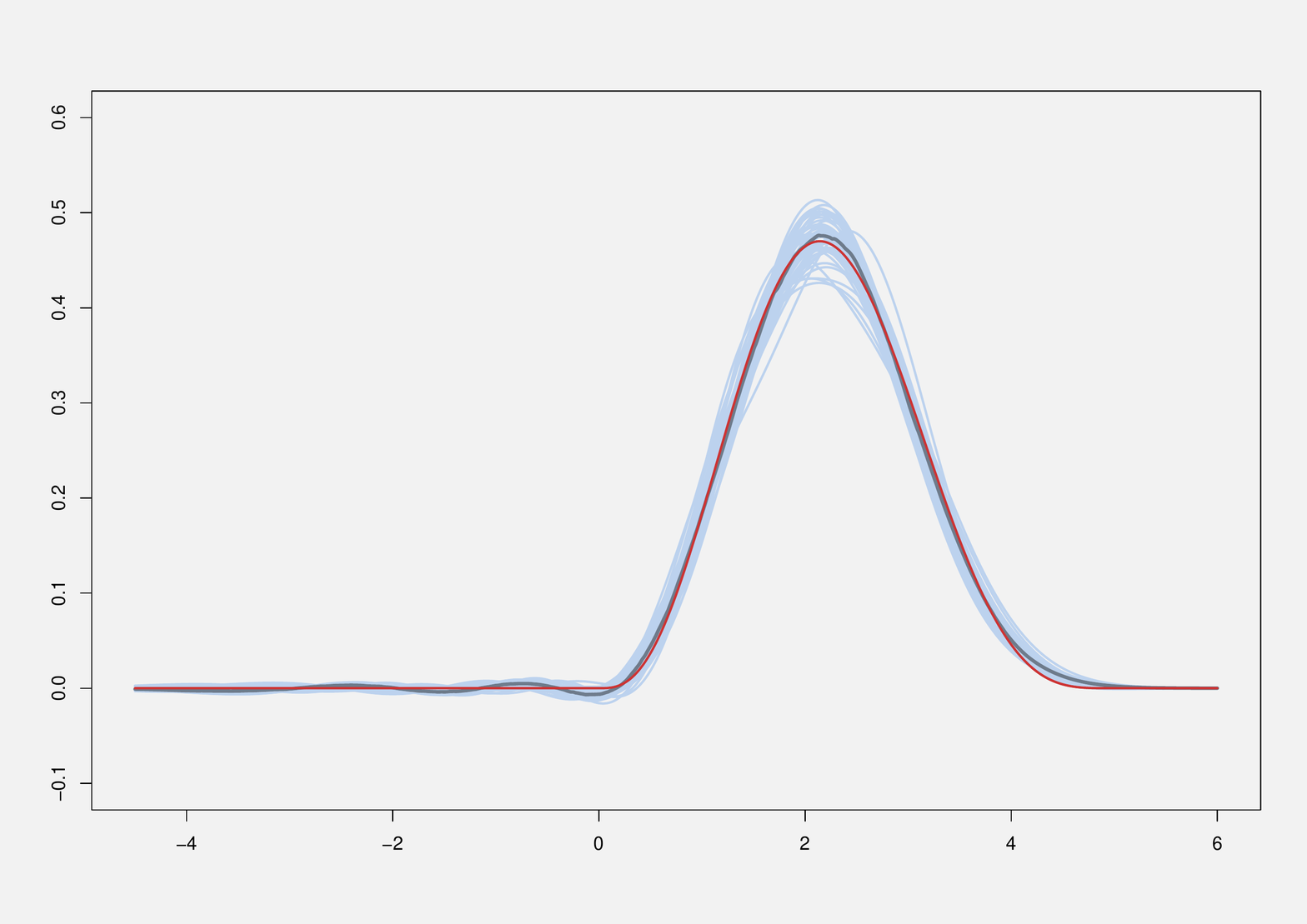}
	\end{minipage}}
	\captionof{figure}{Considering \bwn and a sample size
          $n=1000$ the aggregated estimators are depict for 50
          Monte-Carlo simulations using the Laguerre (top) and Hermite
          (bottom) basis in the cases \ref{si:lag:i} (left), \ref{si:lag:ii} (middle) and \ref{si:lag:iii} (right). The true density $\So$ is given by the red curve while the dark blue curve is the point-wise empirical median of the 50 estimates.}
\end{minipage}\\[2ex]
\begin{minipage}[t]{\textwidth}
\centerline{\begin{minipage}[t]{0.32\textwidth}
	\includegraphics[width=\textwidth,height=40mm]{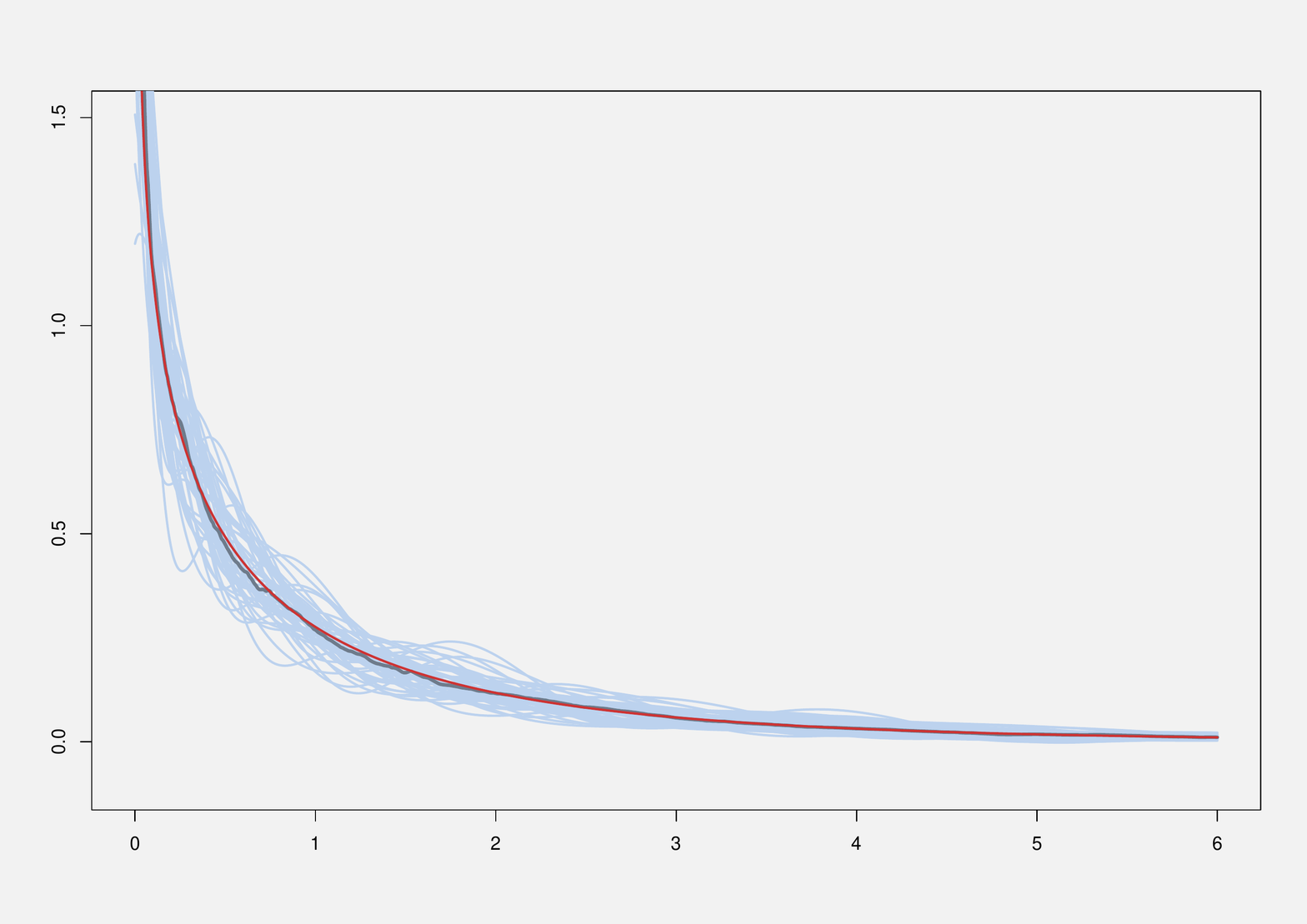}
	\end{minipage}
	\begin{minipage}[t]{0.32\textwidth}
		\includegraphics[width=\textwidth,height=40mm]{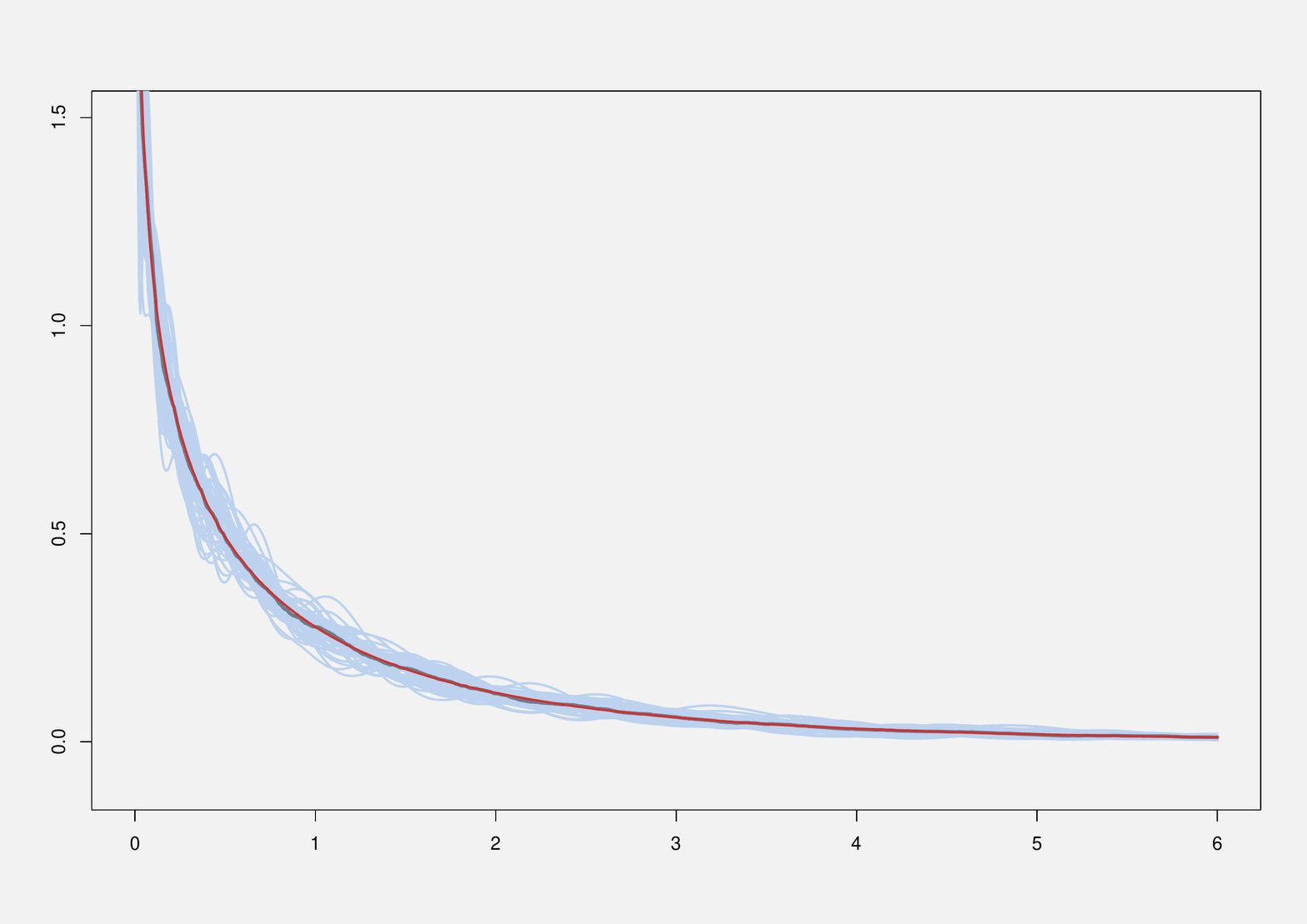}
	\end{minipage}
\begin{minipage}[t]{0.32\textwidth}
	 \includegraphics[width=\textwidth,height=40mm]{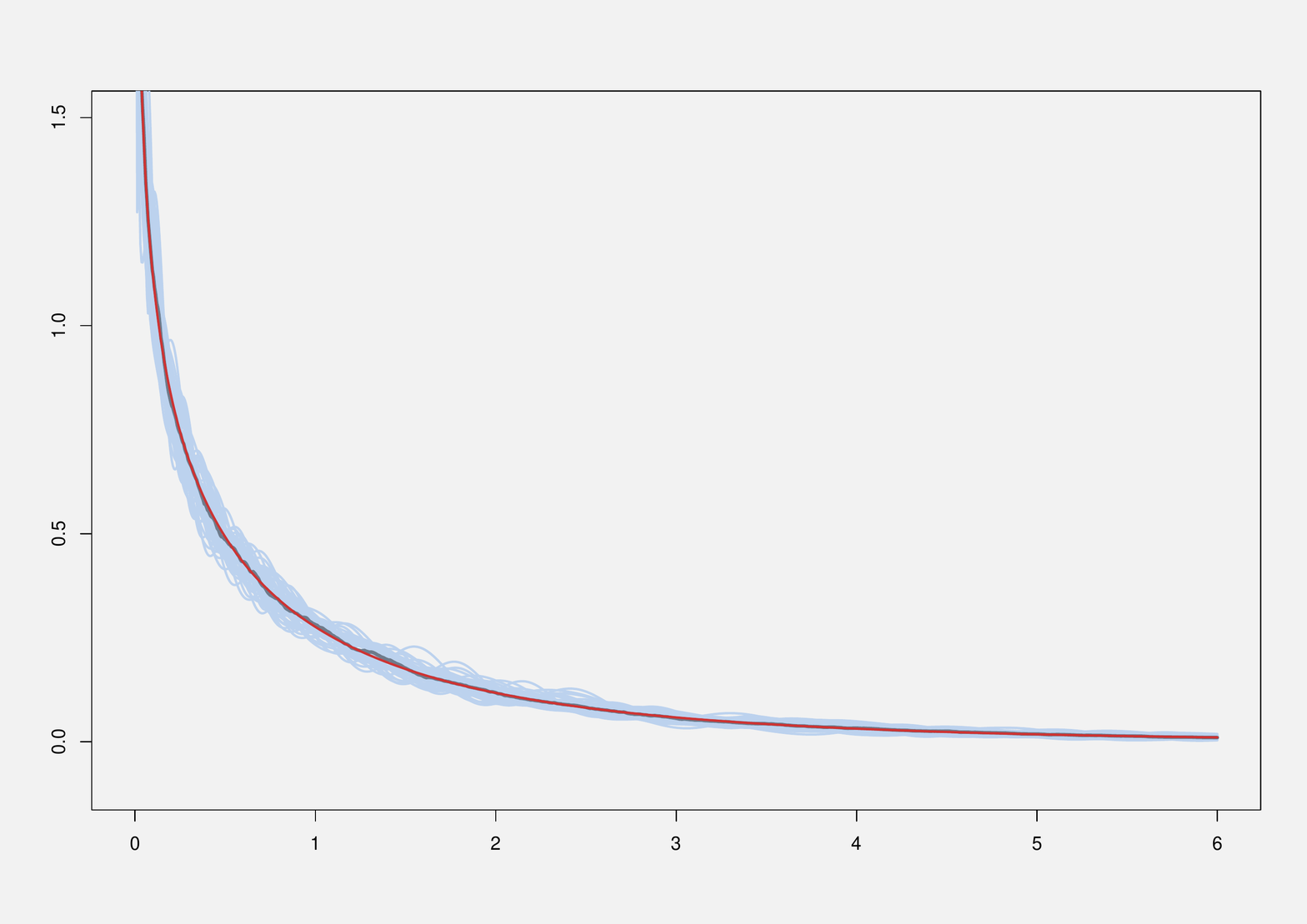}
\end{minipage}}
\centerline{\begin{minipage}[t]{0.32\textwidth}
		\includegraphics[width=\textwidth,height=40mm]{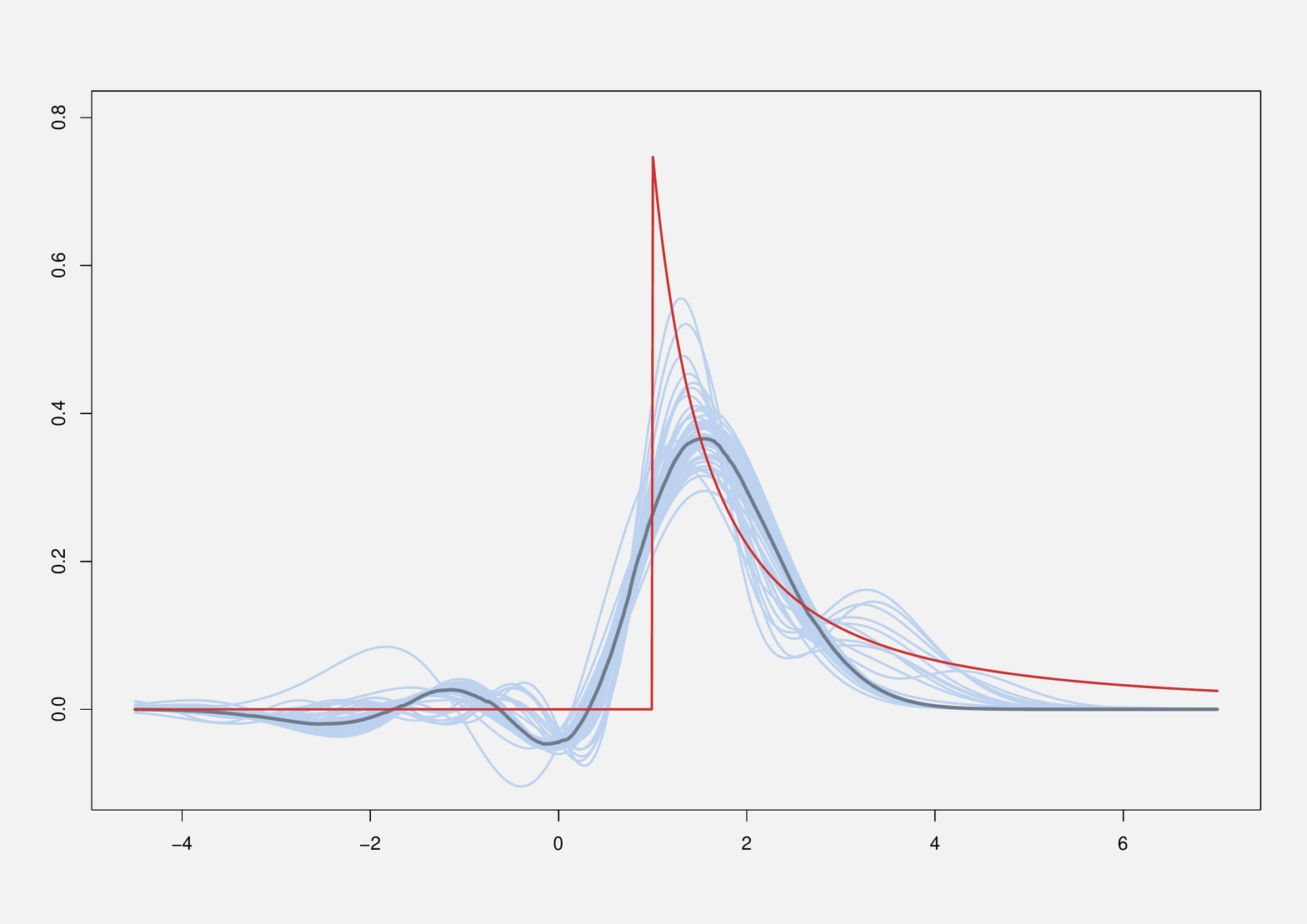}
	\end{minipage}
	\begin{minipage}[t]{0.32\textwidth}
	\includegraphics[width=\textwidth,height=40mm]{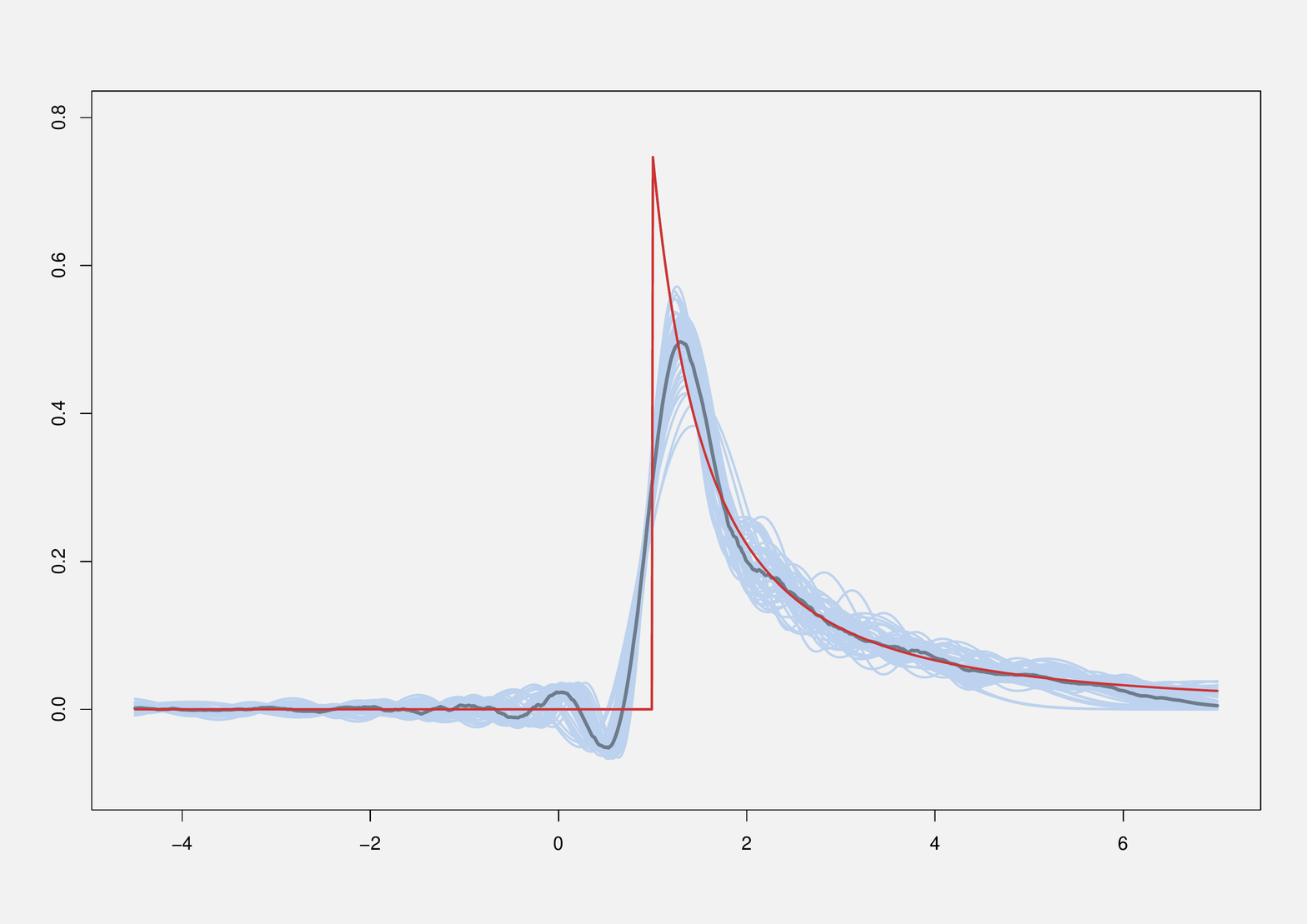}
	\end{minipage}
	\begin{minipage}[t]{0.32\textwidth}
	\includegraphics[width=\textwidth,height=40mm]{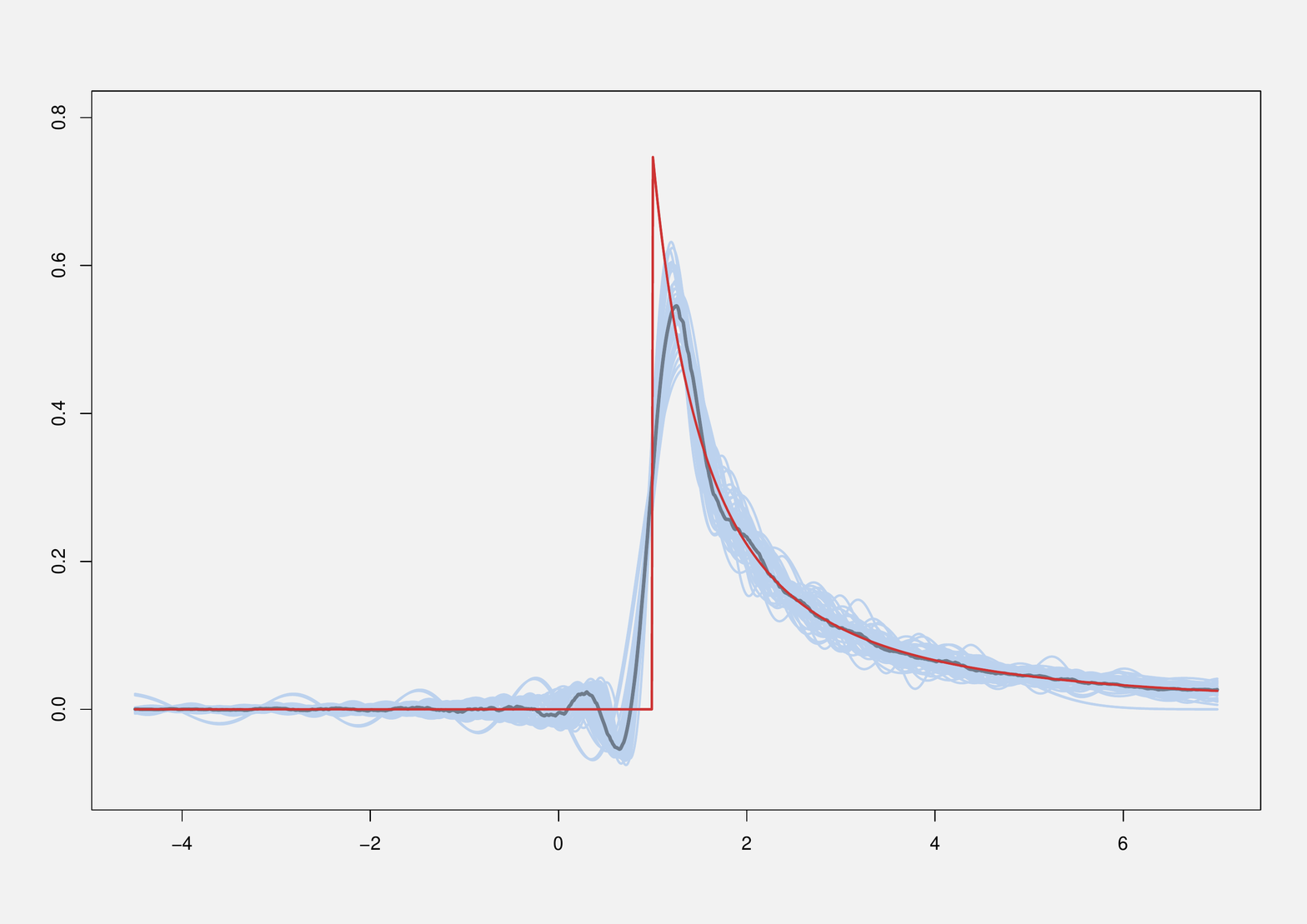}
	\end{minipage}}
\captionof{figure}{Considering \bwn in case \ref{si:lag:iv}
          the aggregated estimators are depict for 50     Monte-Carlo
          simulations using the Laguerre (top) and Hermite (bottom)
          basis   with varying sample size $n=200$ (left), $n=1000$ (middle) and $n=2000$ (right). The true density $\So$ is given by the red curve while the dark blue curve is the point-wise empirical median of the 50 estimates.}
\end{minipage}\\[2ex]
\begin{minipage}[t]{\textwidth}\label{table:num_res}
\centerline{\begin{tabular}{@{}cc|c|c|c|c|c|c|c|c|c@{}}\toprule
&&\multicolumn{3}{c|}{Model selection}&\multicolumn{3}{c|}{Bayesian weights}&\multicolumn{3}{c}{optimal OSE} \\
&$n=$ & 200 & 1000 & 2000 & 200 & 1000  & 2000 &200 & 1000& 2000\\\midrule
\LagCase &\ref{si:lag:i} & 0.899 & 0.408 & 0.282 & 0.529 & 0.226  & 0.160 & 0.456& 0.089& 0.055\\
&\ref{si:lag:ii}& 0.540 & 0.265 & 0.168 & 0.289 & 0.139 & 0.088 &0.096& 0.022&  0.009\\
&\ref{si:lag:iii}& 0.755 & 0.323 & 0.233 & 0.374 & 0.180  & 0.131 &0.265 & 0.075& 0.040\\
&\ref{si:lag:iv}& 1.129 & 0.359& 0.271 & 0.536 & 0.232  &  0.180 &0.315 & 0.086& 0.048\\\midrule
\HerCase&\ref{si:lag:i}& 0.595 & 0.135 & 0.067 & 0.544& 0.124  & 0.064&0.450 & 0.106& 0.060\\
&\ref{si:lag:ii}& 0.179& 0.0039 & 0.018& 0.174 & 0.038  & 0.017 &0.148 & 0.031& 0.015\\
&\ref{si:lag:iii}& 0.454 & 0.121& 0.061 & 0.415 & 0.114  & 0.054 &0.308 & 0.096& 0.043\\
&\ref{si:lag:iv}& 6.411 & 3.316 & 2.520 & 6.357 & 3.232  & 2.471 & 2.968 & 1.552& 1.230\\\bottomrule
\end{tabular}}
\captionof{table}{Over 50 Monte-Carlo  simulations  accumulated
          squared distances between the unknown density $\So$ and
          the aggregated estimator   $\hSoPr[{\msWe[]}]$
          with \mwn (left), $\hSoPr[{\rWe[]}]$ with \bwn
          (middle) and the OSE $\hSoPr[k_{\mathrm{opt}}]$
          (right) are presented where $k_{\mathrm{opt}}$ minimises in
          each iteration the squared distance $\hSoPr[k]$ and $\So$
          over $\nset{1,\maxDi}$.}
\end{minipage}


%% file: _app1.tex
%
%
%
\subsection{Preliminaries}\label{a:prel}
%


This section gathers preliminary technical results.
\input{_app1_A.tex}
For abbreviation,
we denote by $\Proj$ and $\ProjC$ the orthogonal projections on the
linear subspace $S_{\Di}$ and its orthogonal complement $S_{\Di}^{\perp}$
in $\LpA$, respectively.  
The next result
can be found in \cite{JohannesSimoniSchenk2015}. 
\begin{lem}\label{re:diff_fun}
Given $n\in\Nz$ and
  $\So,\bSo\in\LpA$ consider the
  families of  orthogonal projections
  $\setB{\SoPr:=\Proj\So,\Di\in\nset{1,n}}$ and
  $\setB{\bSoPr:=\Proj\bSo,\Di\in\nset{1,n}}$.  For any $l\in\nset{1,n}$ hold
\begin{resListeN}[]
\item\label{re:diff_fun:i}
$\VnormA{\bSoPr[k]}^2-\VnormA{\bSoPr[l]}^2\leq
\tfrac{11}{2}\VnormA{\bSoPr[l]-\SoPr[l]}^2-\tfrac{1}{2}\VnormA{\So}^2\{\bias[k]^2(\So)-\bias[l]^2(\So)\}$,
for all $k\in\nsetro{1,l}$;
\item\label{re:diff_fun:ii}
$\VnormA{\bSoPr[k]}^2-\VnormA{\bSoPr[l]}^2\leq \tfrac{7}{2}\VnormA{\bSoPr[k]-\SoPr[k]}^2+\tfrac{3}{2}\VnormA{\So}^2
\{\bias[l]^2(\So)-\bias[k]^2(\So)\}$, for all $k\in\nsetlo{l,n}$.
\end{resListeN}
\end{lem}

\begin{te}
The next assertion
provides our key arguments in order to control the deviations of the
reminder terms.  Both inequalities are due to
\cite{Talagrand1996}, the formulation of the first part  can be found for
example in \cite{KleinRio2005}, while the second part is  based on
equation (5.13) in Corollary 2 in \cite{BirgeMassart1995} and stated
in this form for example in  \cite{ComteMerlevede2002}.
\end{te}
\begin{lem}(Talagrand's inequalities)\label{tal:re} Let
  $\Ob_1,\dotsc,\Ob_n$ be independent $\cX$-valued random variables and let $\overline{\nu_{\He}}=n^{-1}\sum_{i=1}^n\left[\nu_{\He}(\Ob_i)-\Ex\left(\nu_{\He}(\Ob_i)\right) \right]$ for $\nu_{\He}$ belonging to a countable class $\{\nu_{\He},\He\in\cH\}$ of measurable functions. Then,
\begin{align}
	 &\Ex\vectp{\sup_{\He\in\cH}|\overline{\nu_{\He}}|^2-6\TcH^2}\leq C \left[\frac{\Tcv}{n}\exp\left(\frac{-n\TcH^2}{6\Tcv}\right)+\frac{\Tch^2}{n^2}\exp\left(\frac{- n \TcH}{100\Tch}\right) \right];\label{tal:re1} \\
	&\pM\big(\sup_{\He\in\cH}|\overline{\nu_{\He}}|^2\geq6\TcH^2\big)
   \leq 3\big[\exp\big(\frac{-n\TcH^2}{400\Tcv}\big)+\exp\big(\frac{-n\TcH}{200\Tch}\big)\big]\label{tal:re2}
\end{align}
 for some  numerical constant $C>0$ and where
\begin{equation}\label{tal:re3}
	\sup_{\He\in\cH}\sup_{z\in\cZ}|\nu_{\He}(z)|\leq \Tch,\qquad \Ex(\sup_{\He\in\cH}|\overline{\nu_{\He}}|)\leq \TcH,\qquad \sup_{\He\in\cH}\frac{1}{n}\sum_{i=1}^n \Var(\nu_{\He}(\Ob_i))\leq \Tcv.
\end{equation}
\end{lem}
\begin{rem}\label{remark:emp_proc}For $\Di\in\Nz$ define the unit ball
  $\mBaH:=\set{\He\in S_{\Di}:\VnormA{\He}\leq1}$ 
 contained in the linear subspace
 $S_{\Di}$ spanned by $\set{\bas_{j},j\in\nsetro{0,\Di}}$. For
 $\He\in\mBaH$ we define
 $\nu_{\He}:=\Proj\He$, 
 where
 $\Ex[\So](\nu_{\He}(\Ob_1)) = \sum_{j=0}^{\Di-1} \fHe{j}\fSo{j}$ and
 $\overline{\nu_{\He}}=\sum_{j=0}^{\Di-1}(\hfSo{j}-\fSo{j})\fHe{j}$, thus
\begin{equation*}
\VnormA{\hSoPr-\SoPr}^2=\sup_{\He\in\mBaH}|\VskalarA{\hSoPr-\SoPr,\He}|^2=\sup_{\He\in\mBaH}|\sum_{j=0}^{\Di-1}(\hfSo{j}-\fSo{j})\fHe{j}|^2=\sup_{\He\in\mBaH}|\overline{\nu_{\He}}|^2.
\end{equation*}
The last identity provides the necessary argument to apply Talagrand's inequality \cref{tal:re}
in the proof of \cref{re:conc}.  Note that, the unit ball $\mBaH$ is not a countable set of functions, however, it contains a countable dense subset, say $\cH$, since $\LpA$ is separable, and it is straightforward to see that $\sup_{\He\in\mBaH}|\overline{\nu_{\He}}|^2=\sup_{\He\in\cH}|\overline{\nu_{\He}}|^2$. \remEnd
\end{rem}
\begin{lem}\label{re:conc}
  Consider $\Vcst[]\in\pRz$  and $\vFu\geq1$ as in
  \eqref{de:bas:cst}. 
  There exists a numerical constant $\cst{}$
  such that for any density $\So\in\LpA$ with
  $\VnormInf{\So}<\infty$, for all $n\geq3$ and 
  $\Di\in\nset{1,\maxDi}$  hold
  \begin{resListeN}[]
  \item\label{re:conc:i}
    $\nEx \vectp{\VnormA{\hSoPr-\SoPr}^2 - 6\, \Vcst\, \vFu\, \Di^{1/2}\, n^{-1}}   
       \leq 
   \cst{} \big(\frac{\VnormInf{\So}}{n}\exp\big(\frac{-
         \Vcst\vc}{6\VnormInf{\So}}\Di^{1/2}\big)+\frac{1}{n\maxDi}\big)
$;
  \item\label{re:conc:ii}
    $\nVg\big(\VnormA{\hSoPr-\SoPr}^2 \geq 6\, \Vcst\, \vFu\, \Di^{1/2}\, n^{-1}\big)\leq  \cst{}\big(\exp\big(\frac{-
       \Vcst\vc}{400\VnormInf{\So}}\Di^{1/2}\big)+\frac{1}{n\maxDi}\big) $
  \end{resListeN}
where  in  case \LagCase  $\maxDi:=n^2(600\log n)^{-4}$ and in case \HerCase $\maxDi:=n^2$.
\end{lem}
\begin{pro}[Proof of \cref{re:conc}.]\label{pro:conc}
  For $\He\in\mBaH$ setting
  $\nu_{\He}=\Proj\He=\sum_{j=0}^{\Di-1}\fHe{j}\bas_j$
  we observe (see \cref{remark:emp_proc}) that
  $\VnormA{\hSoPr-\SoPr}^2=\sup_{\He\in\mBaH}|\overline{\nu_{\He}}|^2$. We
  intent to apply \cref{tal:re}. Therefore, we compute next quantities
  $\Tch$, $\TcH$, and $\Tcv$ verifying the three inequalities \eqref{tal:re3} required
  in \cref{tal:re}. First, making use of   
  \eqref{de:bas:cst} we have
    $\sup_{\He\in\mBaH}\VnormInf{\nu_{\He}}^2 
    =\VnormInf{\sum_{j=0}^{\Di-1}\bas_j^2}  \leq \Vcst[]\,\Di^{\dSob}=:\Tch^2$.
   Next, find $\TcH$. Notice that
  $\VnormA{\hSoPr-\SoPr}^2=
  \sum_{j=0}^{\Di-1}|\hfSo{j}-\fSo{j}|^2$.  As
  $\nEx|\hfSo{j}-\fSo{j}|^2\leq \tfrac{1}{n}
  \FuEx{\So}\big(\bas_j^2(\Ob_1)\big)$ and by employing
  \eqref{de:bas:cst}, we obtain 
    $\nEx\big(\sup_{\He\in\mBaH}|\overline{\nu_{\He}}|^2\big)
    \leq \tfrac{1}{n}
  \sum_{j=0}^{\Di-1}\FuEx{\So}\big(\bas_j^2(\Ob_1)\big)\leq
 \Vcst[]\vc\Di^{1/2}n^{-1} =: \TcH^2$.
   Finally, we set $\Tcv:=\VnormInf{\So}$. Indeed, for each  $\He\in \mBaH$
   it holds
$\tfrac{1}{n}\sum_{i=1}^n\Var(\nu_{\He}(\Ob_i)\leq
  \FuEx{\So}\big(|\nu_{\He}(\Ob_1)|^2\big)\leq \VnormInf{\So}\VnormA{\Proj\He}^2\leq\Tcv$. 
  Replacing in \eqref{tal:re1} the quantities $\Tch,\TcH$ and $\Tcv$
  there is a numerical constant $\cst{}$ such that (keep in mind that $1-2\dSob\leq 0$ and $\vc\geq1$) it holds
  \begin{multline}\label{pro:conc:e1}
   \Ex\vectp{\VnormA{\hSoPr-\SoPr}^2-6\Vcst\vc \Di^{1/2} n^{-1}}
   \\
   \leq 
   \cst{} \left(\frac{\VnormInf{\So}}{n}\exp\left(\tfrac{-
         \Vcst\vc}{6\VnormInf{\So}}\Di^{1/2}\right)+\tfrac{1}{n\maxDi}\frac{\maxDi^{\dSob+1}}{n}\exp\left(\tfrac{-1}{100}\frac{n^{1/2}}{\maxDi^{\dSob/2-1/4}}\right) \right) 
 \end{multline}
 for all
  $\Di\in\nset{1,\maxDi}$. 
  In case \LagCase with $\dSob=1$,  $\maxDi=n^2(600\log
  n)^{-4}$ and $n\geq3$ we have 
  \begin{equation*}
    \frac{\maxDi^{\dSob+1}}{n}\exp\left(\tfrac{-1}{100}(\tfrac{n^2}{\maxDi^{2\dSob-1}})^{1/4}\right)
     \leq
     \exp\left((3-\tfrac{600}{100})(\log n)\right)\leq1,  \end{equation*}
  and  in case \HerCase with $\dSob=10/12$ and $\maxDi=n^2$ we obtain
   \begin{equation*}
    \frac{\maxDi^{\dSob+1}}{n}\exp\left(\tfrac{-1}{100}(\tfrac{n^2}{\maxDi^{2\dSob-1}})^{1/4}\right)
     \leq n^{8/3}\exp\left(-\tfrac{1}{100} n^{1/6}\right)\leq
     \big(\tfrac{1600}{e}\big)^6.
  \end{equation*}
  Consequently, combining
   the bounds for the cases  \LagCase and \HerCase we obtain the
   assertion \ref{re:conc:i}. 
  Analogously, replacing in \eqref{tal:re2} the quantities $\Tch,\TcH$ and $\Tcv$
  we obtain \ref{re:conc:ii}, and we omit the details, 
 which completes the
  proof.\proEnd
\end{pro}
\begin{co}\label{ag:re:nd:rest}Under the assumptions of \cref{re:conc}
  for  $(\pen)_{\Di\in\Nz}$ as in \eqref{ag:de:pen} with $\cpen\geq 84\Vcst{}$ there is a finite numerical constant  $\cst{}>0$ such that for all
$n\geq3$ and  $\mdDi\in\nset{1,\maxDi}$    hold
\begin{resListeN}
\item\label{ag:re:nd:rest1}
$\sum_{\Di=1}^{\maxDi}\nEx\vectp{\VnormA{\hSoPr-\SoPr}^2-\pen/14}\leq
\cst{}\big(\VnormInf{\So}^3\vee1\big)n^{-1}$;
\item\label{ag:re:nd:rest2}
  $\sum_{\Di=1}^{\maxDi}\pen\nVg\big(\VnormA{\hSoPr-\SoPr}^2\geq\pen/14\big)\leq\cst{}\,\vc\,\big(\VnormInf{\So}^3\vee1\big)n^{-1}$;
\item\label{ag:re:nd:rest3}
  $\nVg\big(\VnormA{\hSoPr[\mdDi]-\SoPr[\mdDi]}^2\geq\pen[\mdDi]/14\big)\leq    \cst{} \big(\exp\big(\tfrac{-\Vcst\vc}{400\VnormInf{\So}}(\mdDi)^{1/2}\big)+n^{-1}\big)$.
\end{resListeN}
Consider $\vc$, $\hvc$ and $\vEvC=\set{|\hvc-\vc|\leq\vc/2}$, then for all $n\in\Nz$ holds
\begin{resListeN}[\addtocounter{ListeN}{3}]
\item\label{ag:re:nd:rest4} $\nVg\big(\vEvC\big)\leq 4
  \tfrac{\FuEx{\So}(|X_1|^{2\aSob})}{\vc^2} n^{-1}$ with $\aSob\in\Rz$ as in
  \eqref{de:bas:vFu}.
\end{resListeN}
\end{co}
\begin{pro}[Proof of \cref{ag:re:nd:rest}.]Consider
  \ref{ag:re:nd:rest1}.  Exploiting
   the elementary bounds $\sum_{\Di\in\Nz}\exp\big(-\lambda\Di^{1/2}\big)\leq
  \lambda^{-2}$ for all $\lambda>0$, $\pen/14\geq 6\, \Vcst\,\vc\,
  \Di^{1/2}\, n^{-1}$ for all $\Di\in\nset{1,\maxDi}$ and $\vc\geq1$ 
   from \cref{re:conc} \ref{re:conc:i} follows \ref{ag:re:nd:rest1}, that is,
  \begin{multline*}
    \sum_{\Di=1}^{\maxDi}\nEx \vectp{\VnormA{\hSoPr-\SoPr}^2 - \pen/14}   
       \leq 
   \cst{} \big(\frac{\VnormInf{\So}}{n}\sum_{\Di\in\Nz}\exp\big(-\tfrac{\Vcst\vc}{6\VnormInf{\So}}\Di^{1/2}\big)+\frac{1}{n}\big)\\\leq
   \cst{} \big(\frac{\VnormInf{\So}}{n}
   (\tfrac{6\VnormInf{\So}}{\Vcst\vc})^2 +n^{-1}\big)\leq  \cst{} \big(\VnormInf{\So}^3+1\big)n^{-1}.
  \end{multline*} 
  Analogously, from \cref{re:conc} \ref{re:conc:ii} together
  with $\sum_{\Di\in\Nz}\Di^{1/2}\exp\big(-\lambda\Di^{1/2}\big)\leq3
  \lambda^{-3}$ for all $\lambda>0$ we obtain \ref{ag:re:nd:rest2},
  and we omit the details. 
The assertion  \ref{ag:re:nd:rest3}
follows immediately
  from \cref{re:conc} \ref{re:conc:ii}. It remains to show \ref{ag:re:nd:rest4}.
 Recall that
 $\hvc-\vc=\tfrac{1}{n}\sum_{i=1}^n(|X_i|^{\aSob}-\Ex(|X_i|^{\aSob}))$
 with ${\aSob}$  as in \eqref{de:bas:vFu}, where $Y_i:=|X_i|^{\aSob}-\Ex(|X_i|^{\aSob})$,
$i\in\nset{1,n}$, are independent, identically distributed and
centred. Thereby, we have  $\nEx(|\hvc-\vc|^2) =
n^{-1}\FuEx{\So}(|Y_1|^{2})\leq n^{-1}\FuEx{\So}(|X_1|^{2\aSob})$. By
using Tchebychev’s inequality  we deduce the assertion
\ref{ag:re:nd:rest4}, 
which completes the proof.\proEnd
\end{pro}
\begin{co}\label{ag:re:2nd:rest} For  $(\pen)_{\Di\in\Nz}$ as in \eqref{ag:de:pen} with $\cpen\geq 84\Vcst{}$ there is under \cref{ass:Cw} a finite numerical constant  $\cst{}>0$ such that for all $n\geq3$ and $\mdDi\in\nset{1,\maxDi}$  hold
\begin{resListeN}
\item\label{ag:re:2nd:rest1}
$\sup_{\So\in\rwcSo}\sum_{\Di=1}^{\maxDi}\nEx\vectp{\VnormA{\hSoPr-\SoPr}^2-\pen/14}\leq
\cst{}\, \Vcst[\wSo,\rSob]^{2}\,n^{-1}$;
\item\label{ag:re:2nd:rest2}
  $\sup_{\So\in\rwcSo}\sum_{\Di=1}^{\maxDi}\pen\nVg\big(\VnormA{\hSoPr-\SoPr}^2\geq\pen/14\big)\leq\cst{}\,\Vcst[\wSo,\rSob]^3\,n^{-1}$;
\item\label{ag:re:2nd:rest3}
  $\sup_{\So\in\rwcSo}\nVg\big(\VnormA{\hSoPr[\mdDi]-\SoPr[\mdDi]}^2\geq\pen[\mdDi]/14\big)\leq    \cst{} \big(\exp\big(\tfrac{-\Vcst}{400\Vcst[\wSo,\rSob]}(\mdDi)^{1/2}\big)+n^{-1}\big)$.
\end{resListeN}
Consider $\vc$, $\hvc$ and $\vEvC=\set{|\hvc-\vc|\leq\vc/2}$, then for all $n\in\Nz$ holds
\begin{resListeN}[\addtocounter{ListeN}{3}]
\item\label{ag:re:2nd:rest4} $\sup_{\So\in\rwcSo}\nVg\big(\vEvC\big)\leq 4
  \rSo n^{-1}$.
\end{resListeN}
\end{co}
\begin{pro}[Proof of \cref{ag:re:2nd:rest}.]The assertions follow
  immediately from \ref{ag:re:nd:rest1}-\ref{ag:re:nd:rest4} in
  \cref{ag:re:nd:rest}, respectively,  by using that for all
  $\So\in\rwcSo$ and $\Di\in\Nz$ hold
  $1\leq\vc\vee\VnormA{\So}^2\vee\VnormInf{\So}^2\leq\Vcst[\wSo,\rSob]$
  and  $\FuEx{\So}(|X_1|^{2\aSob})\leq\rSob$, and
   we
  omit the details.\proEnd
\end{pro}


%% file: _app1_A.tex
%
%
%

Let us firstly introduce the Laguerre and Hermite basis and secondly
briefly argue that the inequalities in \eqref{de:bas:cst} are fulfilled by both basis. The \textit{Laguerre functions} are defined by
\begin{align*}
\begin{matrix}
\ell_j(x):=\sqrt{2}L_j(2x)\exp(-x)\Ind{\pRz}(x), &L_j(x):=
\sum_{k=0}^{j}(-1)^k \binom{j}{k} \frac{x^k}{k!}, & x\in\pRz
\end{matrix}
\end{align*} 
where $L_j$ is the \textit{Laguerre polynomial} of order
$j\in\Nz_0$. 
As proven in \cite{Szego1918} the Laguerre polynomials are bounded by
$\exp(x/2)$ and therefore for all $j \in\Nz_0$ the function $\ell_j$
is bounded by $\sqrt{2}$, in equal $\VnormInf{ \ell_j} \leq \sqrt{2}$.
The \textit{Hermite polynomial} $\Ksuite[j\in\Nz_0]{H_j}$   and the \textit{Hermite function} $h_j$ of order $j\in\Nz_0$ is given by
\begin{align*}
\begin{matrix}
h_j(x): = \tfrac{1}{(2^j j! \sqrt{\pi}))^{1/2}}H_j(x)\exp(-\tfrac{x^2}{2}), &
H_j(x):=(-1)^j \exp(x^2) \frac{d^j}{dx^j}(\exp(-x^2)), x\in\Rz,
\end{matrix}
\end{align*}
where for each $j\in\Nz_0$ holds $\VnormInf{ h_j}\leq 1$ (see
\cite{OlverLozierBoisvertClark2010} p.450).  Moreover, due to
\cite{Szego1939} p. 242 there is 
$C_{\infty}\in\pRz$ such that  $  \|h_j\|_{
  \infty} \leq C_{\infty} (j+1)^{-1/12}$ for all $j\in\mathbb{N}_0$  which implies the first part of \eqref{de:bas:cst}.
For the second part of \eqref{de:bas:cst} we refer to
\cite{ComteGenon-Catalot2018}  in case of the  \textit{Hermite functions} while
for the \textit{Laguerre functions} we slightly alternate their proof.  Here we change the upper bound for the integral $I_1$ as follows
\begin{align*}
I_1\leq \frac{c_1^2}{2^{p+1}}\int_0^{1/\nu} u^p(u\nu)^{\delta} f(u/2)du \leq \frac{c_1^2}{\sqrt{2\nu}}\int_0^{1/2\nu} u^{p-1/2} f(u)du \leq \frac{c_1^2 }{\sqrt{2k}} \E_f[X^{p-1/2}].
\end{align*}
Now following the steps as in \cite{ComteGenon-Catalot2018} there exists $C_{p,\delta}\in\pRz$ such that
$\E_f[X^p|\ell_k^{(\delta)}(X)|^2] \leq C_{p, \delta}
(\E_f[X^{p-1/2}]+\E_f[X^p]) k^{-1/2}$ which with $p=\delta=0$ shows the
first part of \eqref{de:bas:cst}.  In the sequel we  stick to the
unified notation of 
an orthonormal basis $\Ksuite[j\in\Nz_0]{\bas_j}$ in 
$\LpA$ where for each $j\in\Nz_0$ in case \LagCase and  \HerCase
$\bas_j := \ell_j$ and $\bas_j := h_j$, respectively. 

%% file: _app2.tex
%
%
%
\subsection{Proof of \cref{theorem:lower_bound}}\label{a:mm}
 Due to the construction \eqref{equation:lobodens} of the functions $\psi_{k,K}$ we easily see that the function  $\psi_{k,K}$ has support on $[1+k/K,1+(k+1)/K]$ which lead to  $\psi_{k,K}$ and $ \psi_{l,K}$ having disjoint support if $k\neq l$.  Further we will choose $\psi\in C_c^{\infty}(\Rz)$, the set of all smooth functions with compact support in $\Rz$, which implies  that $\|\psi\|_{\infty}< \infty$. For instance we can choose $\psi(x):= \sin(2\pi x) g(x)$ where $g(x):= \exp\left(- \frac{1}{1-(2x-1)^2}\right) \1_{(0,1)}(x)$. The function $g\in C_c^{\infty}(\Rz)$ is a often use bump function and it holds for all $x\in \Rz$ that $g(1/2+x)=g(1/2-x)$ which implies that $\int_0^1 \psi(x)dx=0$.
 Keep in mind that  $\SoPr[\bm\theta]$ as in \eqref{equation:lobodens} is a density for each $\bm\theta\in\{0,1\}^K$ and $\delta \leq \delta_{\SoPr[o], \psi}$ and the density $\SoPr[o]$ satisfies $c_{\SoPr[o]}:=\inf_{x \in [1,2]} \SoPr[o] >0$.
\begin{lem}\label{lemma:tsyb_vorb}For $K\geq 8$ there is a subset $\{\bm \theta^{(0)}, \dots, \bm \theta^{(M)}\}$ of $\{0,1\}^K$  with $\bm \theta^{(0)}=(0, \dots, 0)$ such that $M\geq 2^{K/8}$ and for all $j, l \in \llbracket 0, M \rrbracket, j \neq l$ holds $\| \SoPr[\bm\theta^{(j)}]-\SoPr[\bm\theta^{(l)}]\|^2_A \geq \frac{\VnormA{\psi}^2}{8} \delta^2 K^{-2\wSob}$ and $ \text{KL}(\SoPr[\bm\theta^{(j)}], \SoPr[\bm\theta^{(0)}]) \leq \frac{\VnormA{\psi}^2}{c_{f_o}\log(2)} \delta^2 \log(M) K^{-2\wSob-1}$ where $\text{KL}$ is the Kullback-Leibler-divergence.
\end{lem}
\begin{pro}[Proof of \cref{lemma:tsyb_vorb}] 
Since  $(\psi_{k,K})_{k\in\nsetro{0,K}}$ have disjoint support for
$\bm\theta, \bm \theta' \in \{0,1\}^K$ holds
	\begin{align*}
	\VnormA{\SoPr[\bm{\theta}]-\SoPr[\bm{\theta}']}^2 &= \delta^2 K^{-2\wSob} \VnormA{\sum_{k=0}^{K-1} ( \theta_{k+1}-\theta'_{k+1}) \psi_{k,K}}^2 = \delta^2  K^{-2\wSob}\rho(\bm \theta, \bm \theta')  \VnormA{\psi^2_{k,K}}^2
	\end{align*}
where $\rho(\bm \theta, \bm \theta'):= \sum_{j=0}^{K-1} \1_{\{\bm
  \theta_{j+1} = \bm \theta'_{j+1}\}}$ is the usual \textit{Hamming
  distance}. Applying a change of variables $v=xK-K-k$ we conclude $
\VnormA{ \SoPr[\bm{\theta}]-\SoPr[\bm{\theta}']}^2 =\delta^2
\VnormA{\psi}^2 K^{-2\wSob-1}\rho(\bm{\theta},\bm{\theta}')$. Due to
the \textsc{Varshamov-Gilbert Lemma} (see \cite{Tsybakov2008})  for $K
\geq 8$ there is a subset $\{\bm \theta^{(0)}, \dots,
\bm\theta^{(M)}\}$ of $\{0,1\}^K$ with $\bm \theta^{(0)}=(0, \dots,
0)$ such that $\rho(\bm \theta^{(j)}, \bm \theta^{(k)}) \geq K/8$ for
all $j ,k\in \llbracket 0, M \rrbracket, j\neq k $ and $M\geq 2^{K/8}$
 implying the first claim $\VnormA{ \SoPr[\bm\theta^{(j)}]-
  \SoPr[\bm\theta^{(l)}]}^2 \geq \frac{\VnormA{\psi}^2\delta^2}{8}
K^{-2\wSob}.$ 
	For the second part, sincee $\SoPr[o]=f_{\bm\theta^{(0)}}$ and
        $\text{KL}(\SoPr[\bm \theta], \SoPr[o])) \leq
        \chi^2(\SoPr[\bm\theta], \SoPr[o])= \int_{1}^{2}
        (\SoPr[\bm\theta](x) -\SoPr[o](x))^2/\SoPr[o](x) dx$ it is
        sufficient to bound the $\chi$-squared divergence where by construction
	\begin{align*}
	\chi^2(\SoPr[\bm\theta],\SoPr[o])&
	\leq c_{\SoPr[o]}^{-1} \VnormA{\SoPr[\bm \theta]- \SoPr[o]}^2 =c_{\SoPr[o]}^{-1}\VnormA{\psi}^2 \delta^2 K^{-2\wSob-1} \rho(\bm \theta, \bm \theta^{(0)})\leq c_{\SoPr[o]}^{-1}\VnormA{\psi}^2 \delta^2 K^{-2\wSob}.
	\end{align*}
Using  $M \geq 2^K$ follows the second claim $\text{KL}(\SoPr[\bm\theta^{(j)}],\SoPr[\bm\theta^{(0)}])\leq
        \frac{\VnormA{\psi}^2}{c_{\SoPr[o]}\log(2)} \delta^2 \log(M)
        K^{-2\wSob-1}$ which completes the proof. \proEnd
\end{pro}
It remains to show  that $\SoPr[\bm \theta]$, $\bm\theta\in\{0,1\}^K$,  as in \eqref{equation:lobodens}
are elements of the classes $\rwcSobD{\wSob,\rSob,\mSob}$. We will consider the cases \LagCase and \HerCase separately starting with the case \LagCase. A similar result was proven by \cite{BelomestnyComteGenon-Catalotothers2017}  without the additionally moment condition. For the sake of simplicity we denote by $\psi_{k,K}^{(j)}$ the $j$-th derivative of $\psi_{k,K}$ and define the finite constant $C_{j,\infty}:= \max(\|\psi^{(l)}\|_{\infty}, l\in \nset{0,j})$. Here we remark that due to the definition of $\psi_{k,K}$ for any $j\in\Nz$ the functions $\psi_{k,K}^{(j)}$ have also disjoint support for different values of the index $k$.
\begin{lem}\label{lemma:Lag_SobDen}Let
  $\wSob,\mSob\in\Nz$ and $\SoPr[o](x):=
  \frac{x^{\mSob}}{\mSob!}\exp(-x)$, $x\in\pRz$. Then,
  there is $L_{\wSob,\mSob,\delta}>0$ such that $\SoPr[o]$
  and any $\SoPr[\bm{\theta}]$ as in  \eqref{equation:lobodens} with
  $\bm\theta\in\{0,1\}^K$, $K\in\Nz$,  belong to $\rwcSobD[\Rz^+]{\wSob,L_{\wSob,\mSob,\delta},\mSob}$.
\end{lem}
\begin{pro}[Proof of \cref{lemma:Lag_SobDen}]Our proof starts with the
  observation that $a_j(\SoPr[o])=0$ for $j\geq \mSob+1$ and hence
  $|\SoPr[o]|_{\wSob}^2 =\sum_{j=0}^{\mSob} j^{\wSob} a_j(\SoPr[o])^2
  \leq 2 \mSob^{\wSob+1}$. On the other hand we use Lemma 7.1. in
  \cite{Catalotothers2016} to bound
  $|\SoPr[\bm \theta]-\SoPr[o]|_{\wSob}$. Precisely, there exists a
  constant $A(\wSob)>0$ such that
  $|\SoPr[\bm{\theta}] - \SoPr[o]|_{\wSob}^2 \leq A(\wSob) \lefttri
  \SoPr[\bm \theta] -\SoPr[o] \righttri_{\wSob}^2$ with
  $\lefttri \SoPr[\bm \theta] -\SoPr[o] \righttri_{\wSob}^2:=
  \sum_{j=0}^{\wSob} \|\SoPr[\bm \theta]-\SoPr[o] \|_j^2$ and for
  $j\in \llbracket 0, \wSob \rrbracket$
  \begin{equation*}
    \Vnorm[j]{\SoPr[\bm{\theta}]-\SoPr[o]}^2  := \delta^2
    K^{-2\wSob}\int_{\Rz^+}\left(
      x^{j/2}
      \sum_{l=0}^j \begin{pmatrix}
        j \\l \end{pmatrix}
      \sum_{k=0}^{K-1}\theta_{k+1} K^l\psi^{(l)}(xK-K-k)\right)^2
    dx .
  \end{equation*}
  Applying Jensen inequality and using disjoint support and boundness
  of the derivatives implies
  \begin{multline*} \Vnorm[j]{ \SoPr[\bm{\theta}]-\SoPr[o]}^2\leq
    \delta^2 2^j K^{-2\wSob} \sum_{l=0}^j \begin{pmatrix} j \\
      l \end{pmatrix} K^{2l} \sum_{k=0}^{K-1} \int_{1+k/K}^{1+(k+1)/K}
    x^j \left( \psi^{(l)}(xK-K-k)\right)^2 dx. \\\hfill \leq \delta^2
    2^j C_{\infty,\wSob}^2 \sum_{l=0}^j \begin{pmatrix} j \\
      l \end{pmatrix} \sum_{k=0}^{K-1} \int_{1+k/K}^{1+(k+1)/K} 2^jdx
    \leq \delta^2 2^{3\wSob} C_{\infty,\wSob}^2.
  \end{multline*}
  It follows
  $|\SoPr[\bm{\theta}] - \SoPr[o]|_{\wSob}^2 \leq A(\wSob)
  \lefttri \SoPr[\bm \theta] -\SoPr[o] \righttri_{\wSob}^2 \leq
  (\wSob+1) C_{\infty,\wSob}^2 A(\wSob) \delta^2 2^{3\wSob} $
  and hence $|\SoPr[\bm{\theta}]|^2_{\wSob} 
  \leq  2(| \SoPr[\bm{\theta}]-\SoPr[o]|_{\wSob}^2 +
  |\SoPr[o]|_{\wSob}^2)  \leq(\wSob+1)C_{\infty,\wSob}^2\delta^2 2^{3\wSob+1} + 4\mSob^{\wSob+1}$.
  Since $\FuEx[]{\SoPr[\bm\theta]}[X^{-\mSob/2}] =2+\delta K^{-\wSob+ 1} C_{\infty, 0} \leq 2+\delta C_{\infty,0}$
  \cref{lemma:Lag_SobDen} holds true with $L_{\wSob,\mSob,\delta }:=((\wSob+1)C_{\infty,
   \wSob}2\delta^2 2^{3\wSob+1} + 4\mSob^{\wSob+1})\vee( 2+\delta
 C_{\infty, 0})$.
        \proEnd
\end{pro}

For the case \HerCase we exploit an alternative characterisation of the Hermite-Sobolev spaces 
$\rwcSob[\Rz]{\wSob}$, $\wSob\geq 0$ providing a criteria for a
function $f\in \Lp[2]{\Rz}$ to be an element of
$\rwcSob[\Rz]{\wSob}$. Precisely, define the operator
$\mathcal{H}^{\wSob/2}$ mapping $f\in\Lp[2]{\Rz}$ with
$\sum_{j=0}^{\infty} (2j+1)^{\wSob} a_j(f)^2< \infty$ to
$\mathcal{H}^{\wSob/2}(f)= \sum_{j=0}^{\infty} (2j+1)^{\wSob/2}
a_j(f)h_j$ where $|f|_{\wSob}^2 \leq \|\mathcal
H^{\wSob/2}(f)\|_{\Rz}^2$. Due to \cite{BongioanniTorrea2006} (based on \cite{Thangavelu1993})  for any choice of $\iota \in C_c^{\infty}(\Rz)$ there exists a constant $C_{\iota}>0$ such that for all $f\in \Lp[2]{\Rz}$ 
\begin{align}\label{equation:Thanga}
\VnormA[\Rz]{\mathcal H^{s/2}(f \iota) }^2 \leq C_{\iota} \VnormA[\Rz]{(- \Delta + \text{Id})^{s/2} f}^2
\end{align}
where $\text{Id}$ denotes the identity and $\Delta$ the Laplacian
operator. In what follows we will consider only functions $f\in
C_c^{\infty}(\Rz)$ with $\text{supp}(f) \subset [1,2]$. We will fix
$\iota$ such that $\iota \in C_c^{\infty}(\Rz)$ and $\forall x\in
[1,2]: \iota(x)=1$. In this situation we have $f\iota = f$ and
$C_{\iota}$ is a universal constant which allows us to bound
$\VnormA[\Rz]{\mathcal H^{\wSob/2}(f)}^2$ for $f\in
C_{c}^{\infty}(\Rz)$ with $\text{supp}(f)\subset [1,2]$. 
Note that \cite{ComteDuvalSacko2019} used
$C_c^{\infty}(\Rz) \subset \rwcSob[\Rz]{\wSob}$ for all $\wSob \geq 0$
which is due
to 
\cite{StempakTorrea2003}.  However, for our purpose we exploit
\eqref{equation:Thanga} to make explicitly  the dependencies of $|\psi_{k,K}|_{\wSob}$ on its index $K$.
\begin{lem}\label{lemma:Her_SobDen}Let $\wSob,\mSob\in\Nz$ and
  $\SoPr[o](x):=(2\pi)^{-1/2} \exp(-0.5x^2)$, $x\in\Rz$. Then, there is $ L_{\wSob,\mSob,\delta}>0$ such that $\SoPr[o]$ and any $\SoPr[\bm{\theta}]$ as in  \eqref{equation:lobodens} with
  $\bm\theta\in\{0,1\}^K$, $K\in\Nz$,  belong to $\rwcSobD[\Rz]{\wSob,L_{\wSob,\mSob,\delta},\mSob}$. 
\end{lem}
\begin{pro}[Proof of \cref{lemma:Her_SobDen}]We start our proof with
  the observation that  $a_j(\SoPr[o])=0$  for every $j\geq 1$
  implying $|\SoPr[o]|_{\wSob} = 0$ and thus $\SoPr[o] \in
  \rwcSob[\Rz]{\wSob, \rSob}$ for all $L>0$. Setting $\Psi_K:=
  \sum_{k=0}^{K-1} \theta_{k+1}\psi_{k,K}$ and keeping in mind that
  $|\cdot|_{\wSob}$ defines a semi-norm  from \eqref{equation:Thanga} follows
  \begin{align*}
    |\SoPr[\bm\theta]-\SoPr[o]|_{\wSob}^2 &\leq \delta^2 K^{-2\wSob}| \Psi_K|_{\wSob}^2 \leq  \delta^2 K^{-2\wSob}\VnormA[\Rz]{\mathcal H^{\wSob/2} \Psi_K}^2 \leq C_{\iota}\delta^2 K^{-2\wSob} \VnormA[\Rz]{(-\Delta+ \text{Id})^{\wSob/2}\Psi_K}^2.
  \end{align*}
  Let $\mathscr S(\Rz):=\{ f \in C^{\infty}(\Rz) | \forall \alpha,
  \beta\in \Nz: \sup_{x\in \Rz} |x^{\alpha} f^{(\beta)}(x)|< \infty\}$
  be the Schwartz class. Introducing the Bessel potential operator
  $(-\Delta+\text{Id})^{\wSob/2}$ we use  for any $f \in \mathscr
  S(\Rz)$ the identity  $(-\Delta + \text{Id})^{\wSob/2} f=
  G_{\wSob/2} * f$ where $\mathcal F(G_{\wSob/2})(\xi)=
  (1+|\xi|^2)^{\wSob/2}$ for $\xi \in \Rz$. Here  $\mathcal
  F(f)(\xi):=\int_{\Rz} f(x) \exp(-ix\xi)dx$ denotes the usual Fourier
  transform of $f$ evaluated at $\xi \in \Rz$
  (e.g. \cite{AdamsHedberg2012}). In the sequel we use  $\mathcal
  F(G_{\wSob/2}) \mathcal F(G_{\wSob/2})= \mathcal F(G_{\wSob})$ and
  $\Psi_{K}, \psi_{k,K} \in C_c^{\infty}(\Rz) \subset \mathscr S
  (\Rz)$ which together with Plancherel's and the
  convolution theorem imply $\VnormA[\Rz]{(-\Delta +
    \text{Id})^{\wSob/2}\Psi_K}^2=\VnormA[\Rz]{\mathcal F((-\Delta +
    \text{Id})^{\wSob/2}\Psi_K)}^2= \VnormA[\Rz]{\mathcal
    F(G_{\wSob/2})\mathcal F(\Psi_k)}^2$ and
  \begin{equation*}
    \VnormA[\Rz]{(-\Delta + \text{Id})^{\wSob/2}\Psi_K}^2 
    =\langle \mathcal F(\Psi_K), \mathcal F(G_{\wSob})\mathcal F(\Psi_K)\rangle_{\Rz}= \int_{\Rz} \Psi_K(x) (-\Delta+\text{Id})^{\wSob} \Psi_K(x)dx.
  \end{equation*}
  Keeping $\wSob\in \Nz$ in mind for any $f \in C_c^{\infty}$ we have $(-\Delta + \text{Id})^{\wSob}f = \sum_{j=0}^{\wSob}\binom{\wSob}{j} f^{(2j)}$. Since the derivatives $\psi_{k,K}^{(j)}$ have disjoint supports for different values of the index $k$ follows
  \begin{multline*}
    \VnormA[\Rz]{(-\Delta + \text{Id})^{\wSob/2}\Psi_K}^2 = \sum_{k=0}^{K-1} \theta_{k+1}\int_{1+k/K}^{1+(k+1)/K} \psi_{k,K}(x) (-\Delta +\text{Id})^{\wSob} \psi_{k,K}(x)dx \\
    \leq  \sum_{k=0}^{K-1} \int_{1+k/K}^{1+(k+1)/K} |\psi_{k,K}(x) (-\Delta +\text{Id})^{\wSob} \psi_{k,K}(x)|dx.
  \end{multline*}
  Note that
		 $\|(-\Delta+\text{Id})^{\wSob} \psi_{k,K}\|_{\infty} \leq \sum_{j=0}^{\wSob} \binom{\wSob}{j} K^{2j} \| \psi^{(2j)}\|_{\infty} \leq 2^{\wSob}K^{2\wSob} C_{\infty,2\wSob}=C_{\wSob} K^{2\wSob}$. From $\int_{-\infty}^{\infty} |\psi_{k,K}(x)|dx  \leq C_{\infty, 0} K^{-1} $ follows  $\VnormA[\Rz]{(-\Delta+\text{Id})^{s/2} \Psi_K}^2 \leq C_{\wSob} K^{2\wSob} $ and whence $|\SoPr[\bm \theta]|_{\wSob}^2 =  |\SoPr[\bm\theta]-\SoPr[o] |_{\wSob}^2 \leq C_{\wSob} \delta^2.$ 
Since  $\FuEx[]{\SoPr[\bm\theta]}[|X|^{2\mSob/3}]\leq  \FuEx{\SoPr[\bm\theta]}[|X|^{2\mSob}]^{1/3}$ and
	$\FuEx{\SoPr[\bm\theta]}[|X|^{2\mSob}] 	= \int_{-\infty}^{\infty} x^{2\mSob} \SoPr[o](x)dx + \delta K^{-\wSob} \sum_{k=0}^{K-1} \int_1^2 x^{2\mSob} \psi_{k,K }(x)dx\leq (2\mSob)!+4^{\mSob}\delta C_{\infty, 0}$
\cref{lemma:Her_SobDen} is satisfied with $L_{\wSob, \mSob, \delta}:=
(C_{\wSob}\delta^2)\vee ((2\mSob)!+4^{\mSob}\delta
C_{\infty,0})^{1/3}$, which completes the proof. \proEnd
\end{pro}


%% file: _app3.tex
%
%
%
\subsection{Proofs of \cref{ag}}\label{a:ag}
\begin{pro}[Proof of \cref{co:agg}.]
We start the proof with the observation that
\begin{multline*}
  \fou[j]{\hSoPr[{\We[]}]}-\fSo{j}=(\hfSo{j}-\fSo{j})\FuVg{\We[]}(\nsetlo{j,\maxDi})-\fSo{j}\FuVg{\We[]}(\nset{1,j})\text{ for all }j\in\nsetro{0,\maxDi}\\\text{ and }\fou[j]{\hSoPr[{\We[]}]}-\fSo{j}=-\fSo{j}\text{ for all }j\geq\maxDi.
\end{multline*}
Consequently, we  have
  \begin{multline}\label{co:agg:pro1}
    \VnormA{\hSoPr[{\We[]}]-\So}^2
\leq
   2\sum_{j\in\nsetro{0,\maxDi}}\big(|\hfSo{j}-\fSo{j}|^2 \FuVg{\We[]}(\nsetlo{j,\maxDi})\big) \\+ 2\sum_{j\in\nsetro{1,\maxDi}}|\fSo{j}|^2\FuVg{\We[]}(\nset{1,j})+\sum_{j\geq\maxDi}|\fSo{j}|^2,
 \end{multline}
where we bound the first  and the two other terms on the right hand
side separately. Considering the first term we split the sum into two parts.  Precisely,
for $\pDi\in\nset{1,\maxDi}$ holds
\begin{multline}\label{co:agg:pro2}
\sum_{j\in\nsetro{0,\maxDi}}(\hfSo{j}-\fSo{j})^2
\FuVg{\We[]}(\nsetlo{j,\maxDi})
\leq
\VnormA{\hSoPr[\pDi]-\SoPr[\pDi]}^2
+\sum_{l\in\nsetlo{\pDi,\maxDi}}\We[l]\VnormA{\hSoPr[l]-\SoPr[l]}^2\\
\hfill\leq \tfrac{1}{14}\pen[\pDi]
+\sum_{l\in\nset{\pDi,\maxDi}}\vectp{\VnormA{\hSoPr[l]-\SoPr[l]}^2-\pen[l]/14}\\
+\tfrac{1}{14}\sum_{l\in\nsetlo{\pDi,\maxDi}}\We[l]\pen[l]\Ind{\{\VnormA{\hSoPr[l]-\SoPr[l]}^2\geq\hpen[l]/7\}}
+\tfrac{1}{14}\sum_{l\in\nsetlo{\pDi,\maxDi}}\We[l]\pen[l]\Ind{\{\VnormA{\hSoPr[l]-\SoPr[l]}^2<\hpen[l]/7\}}
\end{multline}
Consider the second and third term in \eqref{co:agg:pro1} we split the first sum into two parts and obtain
\begin{multline}\label{co:agg:pro3}
\sum_{j\in\nsetro{1,\maxDi}}|\fSo{j}|^2\FuVg{\We[]}(\nset{1,j})+\sum_{j\geq\maxDi}|\fSo{j}|^2\\
\hspace*{5ex}\leq  \sum_{j\in\nsetro{1,\mDi}}|\fSo{j}|^2\FuVg{\We[]}(\nset{1,j})+ \sum_{j\in\nsetro{\mDi,\maxDi}}|\fSo{j}|^2+
  \sum_{j\geq\maxDi}|\fSo{j}|^2\\\hfill
\leq \VnormA{\So}^2\{\FuVg{\We[]}(\nsetro{1,\mDi})+\bias[\mDi]^2(\So)\}
\end{multline}
Combining  \eqref{co:agg:pro1} and  \eqref{co:agg:pro2}, \eqref{co:agg:pro3} we obtain   the assertion, which completes the proof.\proEnd
\end{pro}
\subsubsection{Technical assertions used in the proof of \cref{ag:ub:pnp}}\label{a:bm:agg}
\begin{te}
  Below we state and proof the technical \cref{ag:re:SrWe,ag:re:SrWe:ms,ag:ub:p} used
  in the proof of \cref{ag:ub:pnp}. The
  proof of \cref{ag:re:SrWe} is based on \cref{re:rWe} given first.
\end{te}
\begin{lem}\label{re:rWe} Considering \bwn
  $\rWe[]$ as in \eqref{ag:de:rWe} for any $l\in\nset{1,\maxDi}$  hold
  \begin{resListeN}[]
  \item\label{re:rWe:i} for all $k\in\nsetro{1,l}$ we have\\
    $\rWe\Ind{\setB{\VnormA{\hSoPr[l]-\SoPr[l]}^2<\hpen[l]/7}}
    \leq\exp\big(\rWn\big\{\tfrac{25}{14}\hpen[l]+\tfrac{1}{2}\VnormA{\So}^2\bias[l]^2(\So)-\tfrac{1}{2}\VnormA{\So}^2\bias^2(\So)
    -\hpen\big\}\big)$
   \item\label{re:rWe:ii} for all $\Di\in\nsetlo{l,\maxDi}$ we have\\
    $\rWe\Ind{\setB{\VnormA{\hSoPr-\SoPr}^2<\hpen/7}}\leq\exp\big(\rWn\big\{-\tfrac{1}{2}\hpen
    +\tfrac{3}{2}\VnormA{\So}^2\bias[l]^2(\So)+\hpen[l]\big\}\big)$.
  \end{resListeN}
\end{lem}
\begin{pro}[Proof of \cref{re:rWe}.]
  Given $\Di,l\in\nset{1,\maxDi}$ and an event $\dmEv{\Di}{l}$ (to be
  specified below) it follows
  \begin{multline}\label{re:rWe:pro1}
    \rWe\Ind{\dmEv{\Di}{l}}
    =\frac{\exp(-\rWn\{-\VnormLp{\hSoPr}^2+\hpen\})}
    {\sum_{l\in\nset{1,n}}\exp(-\rWn\{-\VnormA{\hSoPr[l]}^2+\hpen[l]\})}
    \Ind{\dmEv{\Di}{l}}\\
    \leq
    \exp\big(\rWn\big\{\VnormA{\hSoPr}^2-\VnormA{\hSoPr[l]}^2
    +(\hpen[l]-\hpen)\big\}\big)\Ind{\dmEv{\Di}{l}}
  \end{multline}
  We distinguish the two cases \ref{re:rWe:i} $\Di\in\nsetro{1,l}$ and
  \ref{re:rWe:ii} $\Di\in\nsetlo{l,\maxDi}$.  Consider first \ref{re:rWe:i}
  $\Di\in\nsetro{1,l}$. From \eqref{re:rWe:pro1} and \ref{re:diff_fun:i} in \cref{re:diff_fun} 
  (with  $\bSo:=\hSoPr[n]$) follows
  \begin{multline*}
    \rWe\Ind{\dmEv{\Di}{l}}
    \leq \exp\big(\rWn\big\{\tfrac{11}{2}\VnormA{\hSoPr[l]-\SoPr[l]}^2-\tfrac{1}{2}\VnormA{\So}^2(\bias[k]^2(\So)-\bias[l]^2(\So))+(\hpen[l]-\hpen[k])\big\}\big)\Ind{\dmEv{k}{l}}.
  \end{multline*}
  Setting  $\dmEv{\Di}{l}:=\setB{\VnormA{\hSoPr[l]-\SoPr[l]}^2<\hpen[l]/7}$
  the last bound implies the
  assertion \ref{re:rWe:i}. 
  Consider secondly \ref{re:rWe:ii} $\Di\in\nsetlo{l,\maxDi}$. From \ref{re:diff_fun:ii}
  in \cref{re:diff_fun}  (with  $\bSo:=\hSoPr[n]$) and
  \eqref{re:rWe:pro1} follows
  \begin{multline*}
    \rWe[k]\Ind{\dmEv{l}{k}}
    \leq
    \exp\big(\rWn\big\{\tfrac{7}{2}\VnormA{\hSoPr[k]-\SoPr[k]}^2
    +\tfrac{3}{2}\VnormA{\So}^2(\bias[l]^2(\So)-\bias^2(\So))
    +(\hpen[l]-\hpen)\big\}\big)\Ind{\dmEv{l}{k}}.
  \end{multline*}
  Setting $\dmEv{l}{\Di}:=\{\VnormA{\hSoPr-\SoPr}^2<\hpen/7\}$ and
  exploiting $\bias^2(\So)\geq0$ we obtain \ref{re:rWe:ii}, 
  which completes the proof.\proEnd
\end{pro}
\begin{lem}\label{ag:re:SrWe}Consider \bwn $\rWe[]$
  as in \eqref{ag:de:rWe}.  For any
  $\mdDi,\pdDi\in\nset{1,\maxDi}$ and associated $\pDi,\mDi\in\nset{1,\maxDi}$
  as defined in \eqref{ag:de:*Di} hold
  \begin{resListeN}
  \item\label{ag:re:SrWe:i}
    $\FuVg{\rWe[]}(\nsetro{1,\mDi})\Ind{\vEv}\leq
    \tfrac{4}{\rWc^2\cpen^2\vc^2}\Ind{\{\mDi>1\}}
    \exp\big(-\tfrac{3\rWc\cpen\vc}{28}(\mdDi)^{1/2}\big)
    +\Ind{\{\mDi>1\}}\Ind{\setB{\VnormA{\hSoPr[\mdDi]-\SoPr[\mdDi]}^2\geq\pen[\mdDi]/14}}$;
  \item\label{ag:re:SrWe:ii}
    $\sum_{\Di\in\nsetlo{\pDi,\maxDi}}\pen\rWe
    \Ind{\{\VnormA{\hSoPr-\SoPr}^2<\hpen/7\}}
    \leq n^{-1}\tfrac{192\vc}{\rWc^{3}\cpen^2\hvc^3}\leq
    n^{-1}\tfrac{192\vc}{\rWc^{3}\cpen^2} $ (using $\hvc\geq1$).
  \end{resListeN}
\end{lem}
\begin{pro}[Proof of \cref{ag:re:SrWe}.]
  Consider \ref{ag:re:SrWe:i}. Let $\mDi\in\nset{1,\mdDi}$ as in
  \eqref{ag:de:*Di}. For the non trivial case $\mDi>1$ from
  \cref{re:rWe} \ref{re:rWe:i} with $l=\mdDi$ follows for all
  $\Di<\mDi\leq \mdDi$
  \begin{multline*}
    \rWe\Ind{\setB{\VnormA{\hSoPr[\mdDi]-\SoPr[\mdDi]}^2<\hpen[\mdDi]/7}}
    \\\leq
    \exp\big(\rWn\big\{-\tfrac{1}{2}\VnormA{\So}^2\bias^2(\So)
    +(\tfrac{25}{14}\hpen[\mdDi]
    +\tfrac{1}{2}\VnormA{\So}^2\bias[\mdDi]^2(\So))
    -\hpen\big\}\big).
  \end{multline*}
  By using
  $\tfrac{1}{2}\pen\Ind{\vEv}\leq\hpen\Ind{\vEv}\leq\tfrac{3}{2}\pen$
  and the definition \eqref{ag:de:*Di}  $\mDi$ satisfies
  $\VnormA{\So}^2\bias^2\geq \VnormA{\So}^2\bias[(\mDi-1)]^2>
  \VnormA{\So}^2\bias[\mdDi]^2(\So)+6\pen[\mdDi]\geq
  \VnormA{\So}^2\bias[\mdDi]^2(\So)+4\hpen[\mdDi]\Ind{\vEv}$, which 
  implies
  \begin{multline*}
    \rWe\Ind{\setB{\VnormA{\hSoPr[\mdDi]-\SoPr[\mdDi]}^2<\pen[\mdDi]/7}}\Ind{\vEv}
    \leq\exp\big(-\tfrac{3}{28}\rWn\pen[\mdDi]-\tfrac{1}{2}\rWn\pen\big),\quad\forall \Di\in\nsetro{1,\mDi}.
  \end{multline*}
  The last  bound,  
  $\db\pen=\hvpen$,  $\rWc\cpen\vc>0$,  and
  $\dr\sum_{\Di\in\Nz}\exp(-\lambda\Di^{1/2})\leq \lambda^{-2}$ for any $\lambda>0$
  imply together \ref{ag:re:SrWe:i}, that is,
  \begin{multline*}
    \FuVg{\rWe[]}(\nsetro{1,\mDi})\Ind{\vEv}
    \leq\exp\big(-\tfrac{3\rWc}{28}n\pen[\mdDi]\big)
    \sum_{k=1}^{\mDi-1}\exp(-\tfrac{\rWc\cpen\vc}{2}\Di^{1/2})
    +\Ind{\setB{\VnormA{\hSoPr[\mdDi]-\SoPr[\mdDi]}^2\geq\pen[\mdDi]/14}}\\
    \leq \tfrac{4}{\rWc^2\cpen^2\vc^2}\exp\big(-\tfrac{3\rWc\cpen\vc}{28}(\mdDi)^{1/2}\big)
    +\Ind{\setB{\VnormA{\hSoPr[\mdDi]-\SoPr[\mdDi]}^2\geq\pen[\mdDi]/14}}.
  \end{multline*} 
  Consider \ref{ag:re:SrWe:ii}. Let $\pDi\in\nset{\pdDi,\maxDi}$ as in \eqref{ag:de:*Di}. For the non trivial case $\pDi<\maxDi$ from \cref{re:rWe} \ref{re:rWe:ii}
  with $l=\pdDi$ follows for all $\Di>\pDi\geq \pdDi$
  \begin{equation*}
    \rWe\Ind{\setB{\VnormA{\hSoPr-\SoPr}^2<\hpen/7}}
    \leq \exp\big(\rWn\big\{-\tfrac{1}{2}\hpen
    +\tfrac{3}{2}\VnormA{\So}^2\bias[\pdDi]^2(\So)
    +\hpen[\pdDi]\big\}\big).
  \end{equation*}
  Thereby, for $\pDi$ as in \eqref{ag:de:*Di} satisfying
  $\tfrac{1}{4}\hpen\geq \tfrac{1}{4}\hpen[(\pDi+1)] >
  \tfrac{3}{2}\VnormA{\So}^2\bias[\pdDi]^2(\So)+\hpen[\pdDi]$ holds
  \begin{equation*}
    \rWe\Ind{\setB{\VnormA{\hSoPr-\SoPr}^2<\hpen/7}}
    \leq \exp\big(\rWn\big\{-\tfrac{1}{4} \hpen\big\}\big),\quad\forall\;\Di\in\nsetlo{\pDi,\maxDi}.
  \end{equation*}
   The last upper bound, $\hpen=\hvpen$ and $\pen=\vpen$ imply
  \begin{equation*}
    \sum_{\Di\in\nsetlo{\pDi,\maxDi}}\pen\rWe\Ind{\{\VnormA{\hSoPr-\SoPr}^2<\pen/7\}}
    \leq \cpen\vc n^{-1}\sum_{\Di\in\nsetlo{\pDi,\maxDi}} \Di^{1/2}\exp\big(-\tfrac{\rWc\cpen\hvc}{4}\Di^{1/2}\big),
  \end{equation*}
  which together with
  $\dr\sum_{\Di\in\Nz}\Di^{1/2}\exp(-\lambda\Di^{1/2})\leq 3\lambda^{-3}$  for any $\lambda>0$ implies the
  assertion \ref{ag:re:SrWe:ii} 
  and completes the proof.\proEnd
\end{pro}
\begin{te}
 The next result can be directly deduced from \cref{ag:re:SrWe} by letting
  $\rWc\to\infty$. However, we think the following direct proof  provides  an interesting illustration  of the values
  $\pDi,\mDi\in\nset{1,\maxDi}$ as defined in \eqref{ag:de:*Di}.
\end{te}
\begin{lem}\label{ag:re:SrWe:ms}Consider \mwn $\msWe[]$
  as in \eqref{ag:de:msWe}. For any $\mdDi,\pdDi\in\nset{1,\maxDi}$ and associated
  $\pDi,\mDi\in\nset{1,\maxDi}$ as in \eqref{ag:de:*Di} hold
  \begin{resListeN}[]
  \item\label{ag:re:SrWe:ms:i}
    $\FuVg{\msWe[]}(\nsetro{1,\mDi})\Ind{\vEv}\leq\Ind{\{\mDi>1\}}\Ind{\{\VnormA{\hSoPr[\mdDi]-\SoPr[\mdDi]}^2 \geq\pen[\mdDi]/14\}}$;
  \item\label{ag:re:SrWe:ms:ii}
    $\sum_{\Di\in\nsetlo{\pDi,\maxDi}}\pen\msWe\Ind{\{\VnormA{\hSoPr-\SoPr}^2<\hpen/7\}}=0$.
  \end{resListeN}
\end{lem}
\begin{pro}[Proof of \cref{ag:re:SrWe:ms}.]
  By definition of $\hDi$ it holds
  $-\VnormA{\SoPr[\hDi]}^2+\hpen[\hDi]\leq
  -\VnormA{\SoPr}^2+\hpen$ for all $\Di\in\nset{1,\maxDi}$, and
  hence
  \begin{equation}\label{ag:re:SrWe:ms:pr:e1}
    \VnormA{\SoPr[\hDi]}^2-\VnormA{\SoPr}^2\geq
    \hpen[\hDi]-\hpen\text{ for all }\Di\in\nset{1,\maxDi}.
  \end{equation}
  Consider \ref{ag:re:SrWe:ms:i}. Let $\mDi\in\nset{1,\mdDi}$ as in
  \eqref{ag:de:*Di}. For the non trivial case $\mDi>1$ it is
  sufficient to show, that on  the event
  $\vEv=\{|\hvc-\vc|\leq\vc/2\}$, where $\tfrac{1}{2}\pen\leq\hpen\leq\tfrac{3}{2}\pen$, holds
  $\{\hDi\in\nsetro{1,\mDi}\} \subseteq
  \{\VnormA{\SoPr[\mdDi]-\SoPr[\mdDi]}^2\geq\pen[\mdDi]/14\}$. Indeed,
  on $\vEv$ if  $\hDi\in\nsetro{1,\mDi}$, then the definition
  \eqref{ag:de:*Di} of $\mDi$  implies
  \begin{equation}\label{ag:re:SrWe:ms:pr:e2}
    \VnormA{\So}^2\bias[\hDi]^2(\So)\geq
    \VnormA{\So}^2\bias[(\mDi-1)]^2(\So)>
    \VnormA{\So}^2\bias[\mdDi]^2(\So)+6\pen[\mdDi]\geq\VnormA{\So}^2\bias[\mdDi]^2(\So)+4\hpen[\mdDi].
  \end{equation}
  On the other hand from \cref{re:diff_fun} \ref{re:diff_fun:i}  (with  $\bSoPr[]:=\hSoPr[n]$) follows
  \begin{equation}\label{ag:re:SrWe:ms:pr:e3}
    \VnormA{\hSoPr[\hDi]}^2-\VnormA{\hSoPr[\mdDi]}^2\leq
    \tfrac{11}{2}\VnormA{\hSoPr[\mdDi]-\SoPr[\mdDi]}^2
    -\tfrac{1}{2}\VnormA{\So}^2\{\bias[\hDi]^2(\So)-\bias[\mdDi]^2(\So)\}.
  \end{equation}
  Combining, \eqref{ag:re:SrWe:ms:pr:e1} and
  \eqref{ag:re:SrWe:ms:pr:e3} we conclude
  \begin{equation*}
    \tfrac{11}{2}\VnormA{\hSoPr[\mdDi]-\SoPr[\mdDi]}^2\geq
    \hpen[\hDi]-\hpen[\mdDi]
    +\tfrac{1}{2}\VnormA{\So}^2\{\bias[\hDi]^2(\So)-\bias[\mdDi]^2(\So)\},
  \end{equation*}
  which together with \eqref{ag:re:SrWe:ms:pr:e2} and 
  $\hpen[\hDi]\geq0$ implies
  \begin{multline*}
    \tfrac{11}{2}\VnormA{\hSoPr[\mdDi]-\SoPr[\mdDi]}^2\geq
    \tfrac{1}{2}\VnormA{\So}^2\bias[\hDi]^2(\So)-
    \tfrac{1}{2}\VnormA{\So}^2\bias[\mdDi]^2(\So)
    -\hpen[\mdDi]\\
    >\tfrac{1}{2}(\VnormA{\So}^2\bias[\mdDi]^2(\So)+4\hpen[\mdDi])
    -\tfrac{1}{2}\VnormA{\So}^2\bias[\mdDi]^2(\So)-\hpen[\mdDi]
    \geq\tfrac{11}{14}\hpen[\mdDi].
  \end{multline*}
  Consequently, on $\vEv$ holds
 $\{\hDi\in\nsetro{1,\mDi}\}\subseteq
  \{\VnormA{\hSoPr-\SoPr}^2\geq\hpen[\mdDi]/7\}\subseteq
  \{\VnormA{\hSoPr-\SoPr}^2\geq\pen[\mdDi]/14\}$, which shows
  \ref{ag:re:SrWe:ms:i}.  Consider \ref{ag:re:SrWe:ms:ii}.  Let $\pDi\in\nset{\pdDi,\maxDi}$ as
  in \eqref{ag:de:*Di}. For the non trivial case $\pDi<\maxDi$  it is sufficient to show that,
  $\{\hDi\in\nsetlo{\pDi,\maxDi}\}\subseteq
  \{\VnormA{\hSoPr[\hDi]-\SoPr[\hDi]}^2\geq\hpen[\hDi]/7\}$.  If
  $\hDi\in\nsetlo{\pDi,\maxDi}$ then the definition \eqref{ag:de:*Di}  of $\pDi$  implies
  \begin{equation}\label{ag:re:SrWe:ms:pr:e4}
    \pen[\hDi]\geq \pen[(\pDi+1)] > 6\VnormA{\So}^2\bias[\pdDi]^2(\So)+ 4\hpen[\pdDi]
  \end{equation}
  and due to \cref{re:diff_fun} \ref{re:diff_fun:ii}  (with  $\bSoPr[]:=\hSoPr[n]$) also
  \begin{equation}\label{ag:re:SrWe:ms:pr:e5}
    \VnormA{\hSoPr[\hDi]}^2-\VnormA{\hSoPr[\pdDi]}^2\leq
    \tfrac{7}{2}\VnormA{\hSoPr[\hDi]-\SoPr[\hDi]}^2+\tfrac{3}{2}\VnormA{\So}^2
    \{\bias[\pdDi]^2(\So)-\bias[\hDi]^2(\So)\}.
  \end{equation}
  Combining, \eqref{ag:re:SrWe:ms:pr:e1} and \eqref{ag:re:SrWe:ms:pr:e5} it
  follows that
  \begin{multline*}
    \tfrac{7}{2}\VnormA{\hSoPr[\hDi]-\SoPr[\hDi]}^2\geq
    \hpen[\hDi]-\hpen[\pdDi]  -\tfrac{3}{2}\VnormA{\So}^2
    \{\bias[\pdDi]^2(\So)-\bias[\hDi]^2(\So)\}\hfill
  \end{multline*}
  which  with $\bias[\hDi]^2(\So)\geq0$ and
  \eqref{ag:re:SrWe:ms:pr:e4} implies $\{\hDi\in\nsetlo{\pDi,\maxDi}\}\subseteq
  \{7\VnormA{\hSoPr[\hDi]-\SoPr[\hDi]}^2\geq\hpen[\hDi]\}$,
  that is
  \begin{multline*}
    \tfrac{7}{2}\VnormA{\hSoPr[\hDi]-\SoPr[\hDi]}^2\geq
    (\tfrac{1}{2}+\tfrac{1}{2})\hpen[\hDi]-\hpen[\pdDi]  -\tfrac{3}{2}\VnormA{\So}^2
    \bias[\pdDi]^2(\So)\\
    >\tfrac{1}{2}\hpen[\hDi]+\tfrac{1}{2}\big( 6\VnormA{\So}^2\bias[\pdDi]^2(\So)+ 4\hpen[\pdDi]\big)-\hpen[\pdDi]-\tfrac{3}{2}\VnormA{\So}^2
    \bias[\pdDi]^2(\So)
    \geq\tfrac{1}{2}\hpen[\hDi].
  \end{multline*}
  Thereby, we have shown \ref{ag:re:SrWe:ms:ii} and completed the proof.\proEnd
\end{pro}
\begin{lem}\label{ag:ub:p}
  Let the assumptions of \cref{ag:ub:pnp} be satisfied. Considering
  \ref{ag:ub:pnp:p} there is a finite constant $\cst{\So}$
given in \eqref{ag:ub:pnp:p6} such that for all $n\in\Nz$ holds
$  \oRi  \leq \cst{\So}n^{-1}$.
  \end{lem}
  \begin{pro}[Proof of \cref{ag:ub:p}.]The proof is based on
    an evaluation of the   upper bound 
\eqref{co:agg:e4}  for a suitable selection of the parameters $\mdDi,\pdDi\in\nset{1,\maxDi}$.
Considering \ref{ag:ub:pnp:p} there is $K\in\Nz$   with   $1\geq \bias[{[K-1] }](\xdf)>0$ and
$\bias(\xdf)=0$ for all $\Di\geq K$.
Let $c_{\So}:=\tfrac{6\cpen\vc}{\VnormA{\So}^2\bias[{[K-1]}]^2(\So)}>0$
and $n_{\So}:=\min\set{n\in\Nz:n>c_{\So}K^{1/2}\wedge M_n\geq K}\in\Nz$. We distinguish for $n\in\Nz$ the following two
 cases, \begin{inparaenum}[i]\renewcommand{\theenumi}{\dgrau\rm(\alph{enumi})}\item\label{ag:ub:pnp:p:c1}
$n\in\nset{1,n_{\So}}$ and \item\label{ag:ub:pnp:p:c2}
$n> n_{\So}$. \end{inparaenum}
Firstly, consider
\ref{ag:ub:pnp:p:c1} with $n\in\nset{1,n_{\So}}$, then setting $\mdDi:=1$, $\pdDi:=1$ we have
$\mDi=1$, $\bias[1](\So)\leq 1\leq n_{\So}n^{-1}$ and $\pen[1]\leq \cpen\vc n^{-1}$. Thereby,  from \eqref{co:agg:e4}
follows
 \begin{multline}\label{ag:ub:p:p1} 
   \nEx\VnormA{\hSoPr[{\We[]}]-\So}^2\leq\cst{}\big(\VnormA{\So}^2\,n_{\So}+
     [\VnormInf{\So}^3\vee1]\,\vc +[\VnormA{\So}^2\vee1]\,  \FuEx{\So}(|X|^{2\aSob})\big) n^{-1} 
\end{multline}
Secondly, consider \ref{ag:ub:pnp:p:c2}, i.e., $n>
n_{\So}$ and thus $K\in\nset{1,\maxDi}$. Setting $\pdDi:=K$, it
follows $\bias[\pdDi](\So)=0$ and $\pen[\pdDi]\leq \cpen\vc(K)^{1/2}n^{-1}$. From
\eqref{co:agg:e4} follows for all $n> n_{\So}$ thus
\begin{multline}\label{ag:ub:pnp:p2}
   \nEx\VnormA{\hSoPr[{\We[]}]-\So}^2\leq  
          2\VnormA{\So}^2\bias[\mDi]^2(\So)
                 + \cst{}\VnormA{\So}^2\Ind{\{\mDi>1\}} \exp\big(\tfrac{-\Vcst\vc}{1\vee400\VnormInf{\So}}(\mdDi)^{1/2}\big)\\
+\cst{}\big(\vc K^{1/2}+  [\VnormInf{\So}^3\vee1]\,\vc
+[\VnormA{\So}^2\vee1]\,  \FuEx{\So}(|X|^{2\aSob}) \big)\, n^{-1} 
\end{multline}
The defining set of
$\sDi{n}:=\max\{\Di\in\nset{K,\maxDi}:n>c_{\So}\Di^{1/2} \}$ is  not
empty, since it contains $K$ by construction for all $n> n_{\So}$. Consequently,  $\sDi{n}\geq
K$ and, hence 
$\bias[\sDi{n}](\So)=0$, and
$(\sDi{n})^{1/2}n^{-1}<c_{\So}^{-1}=\tfrac{\VnormA{\So}^2\bias[{[K-1]}]^2(\So)}{6\cpen\vc}$.
It follows
$\VnormA{\So}^2\bias[{[K-1]}]^2(\So)>6\cpen\vc (\sDi{n})^{1/2}n^{-1}=6\pen[\sDi{n}]+\VnormA{\So}^2\bias[\sDi{n}]^2(\So)$
and trivially
$\VnormA{\So}^2\bias[{K}]^2(\So)=0<6\pen[\sDi{n}]+\VnormA{\So}^2\bias[\sDi{n}]^2(\So)$. Therefore,
setting $\mdDi:=\sDi{n}$ the definition \eqref{ag:de:*Di} of
$\mDi$ implies $\mDi=K$ and hence
$\bias[\mDi]^2(\So)=\bias[K]^2(\So)=0$. From \eqref{ag:ub:pnp:p2} for all $n> n_{\So}$  follows
now 
\begin{multline}\label{ag:ub:pnp:p4}
     \nEx\VnormA{\hSoPr[{\We[]}]-\So}^2\leq  
                  \cst{}\VnormA{\So}^2\exp\big(\tfrac{-\Vcst\vc}{1\vee400\VnormInf{\So}}(\sDi{n})^{1/2}\big)\\
+\cst{}\big(\vc K^{1/2}+  [\VnormInf{\So}^3\vee1]\,\vc
+[\VnormA{\So}^2\vee1]\,  \FuEx{\So}(|X|^{2\aSob}) \big)\, n^{-1} 
\end{multline}
Setting $C_{\So}:=\tfrac{1\vee400\VnormInf{\So}}{\Vcst\vc}\vee1$, 
$\aDi{n}:=\floor{C_{\So}^2n^2(\log n)^{-6}}$ for all
  $n>\exp(c_{\So}^{1/3}C_{\So}^{1/3})\vee \exp(C_{\So})$ holds
  $(\aDi{n})^{1/2}n^{-1}\leq C_{\So}(\log n)^{-3}<c_{\So}^{-1}$ and
  $C_{\So}^2(\log n)^{-2}\leq 1$, thus $\aDi{n}\leq \floor{n^2(\log
    n)^{-4}}\leq \maxDi$. Thereby,  $\sDi{n}\geq\aDi{n}>C_{\So}^2n^2(\log n)^{-6}-1=(C_{\So}n(\log
  n)^{-3}-1)(C_{\So}n(\log n)^{-3}+1)\geq(C_{\So}n(\log n)^{-3}-1)^2$
  and hence
  \begin{multline*}
    \exp\big(\tfrac{-\Vcst\vc}{1\vee400\VnormInf{\So}}(\sDi{n})^{1/2}\big)\leq\exp\big(-C_{\So}^{-1}(\sDi{n})^{1/2}\big)\leq
    \exp\big(-C_{\So}^{-1}(C_{\So}n(\log
    n)^{-3}-1)\big)\\=\exp\big(C_{\So}^{-1}\big)\exp\big(-\tfrac{n}{(\log
      n)^{3}}\big)\leq e n^{-1}\exp\big(-(\log n)(n(\log
    n)^{-4}-1)\big) \end{multline*}
 where 
 $n(\log n)^{-4}>1$  for all $n\geq 5550$. Thereby, for all $n>[5550\vee
  \exp(c_{\So}^{1/3}C_{\So}^{1/3})\vee \exp(C_{\So})]$ holds
  $\exp\big(\tfrac{-\Vcst\vc}{1\vee400\VnormInf{\So}}(\sDi{n})^{1/2}\big)\leq
  en^{-1}$ while for $n\leq [5550\vee
  \exp(c_{\So}^{1/3}C_{\So}^{1/3})\vee \exp(C_{\So})]$ holds
  $\exp\big(\tfrac{-\Vcst\vc}{1\vee400\VnormInf{\So}}(\sDi{n})^{1/2}\big)\leq [5550\vee
  \exp(c_{\So}^{1/3}C_{\So}^{1/3})\vee
  \exp(C_{\So})]n^{-1}$. Combining both bounds and the definition of
  $c_{\So}$ and $C_{\So}$ there is a numerical
  constant $\cst{}$ such that for  all $n\in\Nz$ holds
  \begin{equation*}
   \exp\big(\tfrac{-\Vcst\vc}{1\vee400\VnormInf{\So}}(\sDi{n})^{1/2}\big)
\leq \cst{}\exp\big(8\tfrac{[1\vee\VnormInf{\So}]}{[1\vee\VnormA{\So}^2\bias[{[K-1]}]^2(\So)]}\big)n^{-1}
  \end{equation*}
The last bound together with
  \eqref{ag:ub:pnp:p4} implies
\begin{multline}\label{ag:ub:pnp:p5}
     \nEx\VnormA{\hSoPr[{\We[]}]-\So}^2\leq  
     \cst{}\big(\VnormA{\So}^2\exp\big(8\tfrac{[1\vee\VnormInf{\So}]}{[1\vee\VnormA{\So}^2\bias[{[K-1]}]^2(\So)]}\big)\\
+\vc K^{1/2}+  [\VnormInf{\So}^3\vee1]\,\vc
+[\VnormA{\So}^2\vee1]\,  \FuEx{\So}(|X|^{2\aSob}) \big)\, n^{-1} 
\end{multline}
Combining \eqref{ag:ub:pnp:p2} and
    \eqref{ag:ub:pnp:p5} for \ref{ag:ub:pnp:p:c1}
$n\in\nsetro{1,n_{\So}}$ and \ref{ag:ub:pnp:p:c2}
$n\geq n_{\So}$, respectively,  for all $K\in\Nz$ and for all
$n\in\Nz$ follows the claim of \cref{ag:ub:pnp}, that is
\begin{multline}\label{ag:ub:pnp:p6}
  \nEx\VnormA{\hSoPr[{\We[]}]-\So}^2\leq   
     \cst{}\big(\VnormA{\So}^2\exp\big(8\tfrac{[1\vee\VnormInf{\So}]}{[1\vee\VnormA{\So}^2\bias[{[K-1]}]^2(\So)]}\big)+ \VnormA{\So}^2   n_{\So}\\
+\vc K^{1/2}+  [\VnormInf{\So}^3\vee1]\,\vc
+[\VnormA{\So}^2\vee1]\,  \FuEx{\So}(|X|^{2\aSob}) \big)\, n^{-1},
\end{multline}
which  completes the
proof.\proEnd\end{pro}
